\documentclass[10pt,a4paper]{article}
\usepackage{amssymb,latexsym,times,standard,graphicx,amsmath}
\usepackage{pdfsync}
\usepackage{mathtools,bm,mathrsfs}
\usepackage[title]{appendix}
\usepackage{xcolor}

\textheight 
            675pt
\textwidth 
           460pt

\newcommand{\ee}{\mathrm{e}}
\newcommand{\ii}{\mathrm{i}}

\newcommand{\dbfk}{\, \mathrm{d}\bfk}

\newcommand{\dx}{\, \mathrm{d}x}

\newcommand{\du}{\, \mathrm{d}u}
\newcommand{\dy}{\, \mathrm{d}y}

\newcommand{\bfA}{{\bm{A}}}

\newcommand{\bfB}{{\bm{B}}}
\newcommand{\bfc}{{\bm{c}}}
\newcommand{\bfC}{{\bm{C}}}

\newcommand{\bfD}{{\bm{D}}}
\newcommand{\bfe}{{\bm{e}}}

\newcommand{\bfF}{{\bm{F}}}
\newcommand{\bfg}{{\bm{g}}}
\newcommand{\bfG}{{\bm{G}}}
\newcommand{\bfh}{{\bm{h}}}
\newcommand{\bfH}{{\bm{H}}}

\newcommand{\bfk}{{\bm{k}}}
\newcommand{\bfK}{{\bm{K}}}

\newcommand{\bfL}{{\bm{L}}}
\newcommand{\bfM}{{\bm{M}}}
\newcommand{\bfn}{{\bm{n}}}

\newcommand{\bfS}{{\bm{S}}}

\newcommand{\bfT}{{\bm{T}}}
\newcommand{\bfu}{{\bm{u}}}

\newcommand{\bfv}{{\bm{v}}}

\newcommand{\bfw}{{\bm{w}}}

\newcommand{\bfx}{{\bm{x}}}

\newcommand{\bfzero}{{\bm{0}}}

\newcommand{\sdfrac}[2]{\mbox{\small$\displaystyle\frac{#1}{#2}$}}

\renewcommand{\AA}{\mathcal A}

\newcommand{\FF}{\mathcal F}
\newcommand{\GG}{\mathcal G}
\newcommand{\HH}{\mathcal H}

\newcommand{\JJ}{\mathcal J}

\newcommand{\OO}{\mathcal O}

\newcommand{\RR}{\mathcal R}

\newcommand{\WW}{\mathcal W}
\newcommand{\XX}{\mathcal X}

\newcommand{\ZZ}{\mathcal Z}

\renewcommand{\i}{\mathrm{i}}

\newcommand{\nablap}{\nabla^\perp}
\newcommand{\nablac}{\nabla \cdot}
\newcommand{\nablapc}{\nablap \cdot}

\DeclareMathOperator{\curl}{curl}
\DeclareMathOperator{\Div}{div}
\DeclareMathOperator{\Grad}{grad}

\DeclareMathOperator{\cosec}{cosec}

\newcommand{\qed}{\hfill$\Box$\bigskip}

\newcommand{\alignqed}{\tag*{$\Box$}}

\newcommand{\nn}{|{\mskip-2mu}|{\mskip-2mu}|}

\newcommand{\udl}[1]{\underline{\smash{#1}}}

\newtheorem{theorem}{Theorem}[section]
\newtheorem{lemma}[theorem]{Lemma}

\newtheorem{proposition}[theorem]{Proposition}
\newtheorem{corollary}[theorem]{Corollary}
\newtheorem{remark}[theorem]{Remark}

\allowdisplaybreaks

\begin{document}
\newcounter{count}

\title{Fully localised three-dimensional solitary water waves on Beltrami flows with strong surface tension}

\author{M. D. Groves\footnote{Fachrichtung Mathematik, Universit\"{a}t des Saarlandes,
Postfach 151150, 66041 Saarbr\"{u}cken, Germany (groves@math.uni-sb.de)}
\and
E. Wahl\'{e}n\footnote{Centre for Mathematical Sciences, Lund University, PO Box 118, 22100 Lund, Sweden (erik.wahlen@math.lu.se)}}
\date{}
\maketitle{}

\begin{abstract}
\emph{Fully localised three-dimensional solitary waves} are
steady water waves which are evan\-escent in every horizontal direction.
This paper presents an existence theory for such waves under the
assumptions that the relative vorticity and velocity fields are parallel
(`Beltrami flows'), that the free surface of the water takes the form $\{z=\eta(x,y)\}$
for some function $\eta\colon {\mathbb R}^2\rightarrow{\mathbb R}$, and that the influence of surface tension is sufficiently
strong. The governing equations are formulated as a single equation
for $\eta$, which is then reduced to a perturbation of the KP-I equation. This equation has recently
been shown to have a family of nondegenerate localised solutions, and an application of a suitable variant of the
implicit-function theorem shows that they persist under perturbations.
\end{abstract}

\section{Introduction}

\subsection{The hydrodynamic problem}

\begin{figure}
\centering
\includegraphics[scale=0.75]{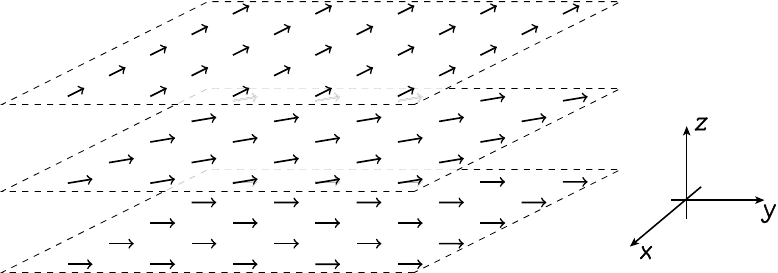}
\caption{\emph{The trivial flow \eqref{ABC flow}}} \label{Trivial flow}
\end{figure}

Consider an incompressible perfect fluid of unit density
occupying a three-dimensional domain bounded below by a rigid horizontal plane and above
by a free surface. A \emph{steady water wave} is a fluid flow of this kind in which both the velocity field and free-surface profile
are stationary with respect to a uniformly (horizontally) translating frame of reference; a \emph{(fully localised) solitary wave} is
a nontrivial steady wave whose free surface decays to the height of the fluid at rest in every horizontal direction. Working in a frame of reference moving with the wave
and in dimensionless coordinates, we suppose that the fluid domain is
$D_\eta=\{(x,y,z)\colon -1 < z < \eta(x,y)\}$ (so that the free surface is the graph $S_\eta$
of an unknown function $\eta\colon {\mathbb R}^2\rightarrow{\mathbb R}$), and
the flow is a \emph{(strong) Beltrami flow} whose velocity and vorticity fields $\bfu$ and $\curl \bfu$ are parallel, so
that $\curl \bfu = \alpha \bfu$ for some fixed constant $\alpha$. \emph{Irrotational} flows (with $\curl \bfu=\mathbf{0}$) are included as the special case $\alpha=0$.
The hydrodynamic problem is to solve the equations
\begin{align}
& \parbox{11.5cm}{$\curl \bfu = \alpha \bfu$}\mbox{in $D_\eta$,} \label{pre-Beltrami1}\\
& \parbox{11.5cm}{$\Div \bfu = 0$}\mbox{in $D_\eta$,} \label{pre-Beltrami2}\\
& \parbox{11.5cm}{$\bfu\cdot\bfe_3 = 0$}\mbox{at $z=-1$,} \label{pre-Beltrami3}\\
& \parbox{11.5cm}{$\bfu\cdot\bfn = 0$}\mbox{at $z=\eta$,} \label{pre-Beltrami4}\\
& \parbox{11.5cm}{$\displaystyle\tfrac{1}{2}|\bfu|^2 + \eta - \beta
\!\left(\frac{\eta_x}{(1+|\nabla\eta|^2)^\frac{1}{2}}\right)_{\!\!\!x}-\beta\!\left(\frac{\eta_y}{(1+|\nabla\eta|^2)^\frac{1}{2}}\right)_{\!\!\!y}=\tfrac{1}{2}|\bfc|^2$}\mbox{at $z=\eta$} \label{pre-Beltrami5},
\end{align}
where $\nabla =(\partial_x,\partial_y)^T$, $\nabla^\perp=(\partial_y,-\partial_x)^T$, $\bfc\coloneqq(c_1,c_2)^T$ is the dimensionless wave velocity, $\bfe_3=(0,0,1)^T$ and
$\bfn\coloneqq(-\eta_x,-\eta_y,1)^T$
is the outward normal vector at $S_\eta$; we have also introduced the Bond number\linebreak
$\beta=\sigma/gh^2$, where $h$ is the depth of the fluid at rest, $g$ is the acceleration due to gravity and
$\sigma>0$ is the coefficient of surface tension. (The pressure $p$ in the fluid
is recovered using the formula $p(x,y,z)=-\tfrac{1}{2}|\bfu(x,y,z)|^2-y$, and the variables $\bfu$ and $p$ automatically solve
the stationary Euler equation in $D_\eta$.)
 Equations \eqref{pre-Beltrami4} and \eqref{pre-Beltrami5} are referred to as respectively the \emph{kinematic} and
 \emph{dynamic} boundary conditions at the free surface.
It is natural to write $\eta$ and $\bfu$ as perturbations of the trivial solution
\begin{equation}
\eta^\star = 0, \qquad \bfu^\star=c_1\begin{pmatrix} \cos \alpha z \\ -\sin \alpha z \\ 0\end{pmatrix}
+c_2 \begin{pmatrix} \sin \alpha z \\ \cos \alpha z \\ 0 \end{pmatrix}
\label{ABC flow}
\end{equation}
 of \eqref{pre-Beltrami1}--\eqref{pre-Beltrami5} (see Figure \ref{Trivial flow}), so that $\bfv=\bfu-\bfu^\star$ satisfies the equations
 \begin{align}
 & \parbox{12.5cm}{$\curl \bfv = \alpha \bfv$}\mbox{in $D_\eta$,} \label{Beltrami 1}\\
& \parbox{12.5cm}{$\Div \bfv = 0$}\mbox{in $D_\eta$,} \label{Beltrami 2}\\
& \parbox{12.5cm}{$\bfv\cdot\bfe_3 = 0$}\mbox{at $z=-1$,} \label{Beltrami 3}\\
& \parbox{12.5cm}{$\bfv\cdot\bfn + \bfu^\star \cdot \bfn= 0$}\mbox{at $z=\eta$,} \label{Beltrami 4}\\
& \parbox{12.5cm}{$\displaystyle\tfrac{1}{2}|\bfv|^2 + \bfv \cdot \bfu^\star + \eta - \beta
\!\left(\frac{\eta_x}{(1+|\nabla\eta|^2)^\frac{1}{2}}\right)_{\!\!\!x}-\beta\!\left(\frac{\eta_y}{(1+|\nabla\eta|^2)^\frac{1}{2}}\right)_{\!\!\!y}=0$}\mbox{at $z=\eta$.}\label{Beltrami 5}
\end{align}
Our task is to find solutions $(\eta,\bfv)$ of \eqref{Beltrami 1}--\eqref{Beltrami 5} which are evanescent as
$|(x,y)| \rightarrow \infty$ and therefore represent fully localised solitary waves `riding' the trivial flow \eqref{ABC flow}.
Note that these equations are invariant under
$$\eta(x,y)\!\mapsto\!\eta(-x,-y), \qquad (v_1(x,y,z),v_2(x,y,z),v_3(x,y,z)\!) \!\mapsto\! (v_1(-x,-y,z),v_2(-x,-y,z),-v_3(-x,-y,z)\!),$$
and we in fact seek solutions which are themselves invariant under this transformation.

Irrotational fully localised solitary waves have been found by Groves \& Sun \cite{GrovesSun08} and Buffoni \emph{et al.} \cite{BuffoniGrovesSunWahlen13}
for $\beta>\tfrac{1}{3}$ using variational methods. Their result has recently been made more precise by Gui \emph{et al.} \cite{GuiLaiLiuWeiYang25a,GuiLaiLiuWeiYang25b}, who
obtained waves  which are perturbations and scalings of localised `lump' solutions of the KP-I equation and discussed their stability. In this paper we use a related method to
obtain the same family of waves for equations \eqref{Beltrami 1}--\eqref{Beltrami 4} for sufficiently large values of $\beta$ (depending upon $\alpha$);
the existence result of Gui \emph{et al.}  is included as a special case. Other types of three-dimensional steady water waves have also been studied, in particular
\emph{doubly periodic} steady waves, that is waves which are periodic in two different horizontal directions. Their existence was established for irrotational
flows with surface tension ($\beta>0$) by Craig \& Nicholls \cite{CraigNicholls00,CraigNicholls02}, for irrotational waves without surface tension $(\beta=0$) by
Iooss \& Plotnikov \cite{IoossPlotnikov11} and for Beltrami flows with surface tension by Lokharu, Seth \& Wahl\'{e}n \cite{LokharuSethWahlen20}
(see also Groves \emph{et al.} \cite{GrovesNilssonPasqualiWahlen24} for an existence theory in a framework similar to that used in the present paper).
Doubly periodic gravity-capillary waves with more general (but small) vorticity were recently constructed by Seth, Varholm \& Wahl\'{e}n \cite{SethVarholmWahlen24}.
 
\subsection{Heuristics}

The KP-I equation arises in weakly nonlinear theory as a universal model equation for two-dimensional nonlinear dispersive systems
whose linearisation has a distinguished wave speed attained only by long waves. Suppose that $|\alpha| \in [0,\tfrac{1}{2}\pi)$ and $\beta>\beta_0$, where
$$\beta_0=\beta_0=\frac{1}{2\alpha^2}(-\cos\alpha+\alpha\cosec\alpha).$$
A straightforward calculation shows that the linearised version of \eqref{Beltrami 1}--\eqref{Beltrami 5} has solutions of the form\linebreak $\eta(x,y)\sim\cos (\bfk\cdot\bfx)$, where $\bfx=(x,y)^\mathrm{T}$ and $\bfk=(k_1,k_2)^\mathrm{T}$, if
\begin{equation}
g(\bfk)=0, \label{dispersion relation - intro}
\end{equation}
where
$$
g(\bfk)=-\frac{1}{|\bfk|^2}(\alpha(\bfc\cdot\bfk^\perp)(\bfc\cdot\bfk)+\mathtt{c}(|\bfk|^2)(\bfc\cdot\bfk)^2)+1+\beta |\bfk|^2
$$
and
$$
\mathtt{c}(\mu) =
\begin{cases}
\sqrt{\alpha^2-\mu} \, \cot(\sqrt{\alpha^2-\mu}),   & \mbox{if $\mu < \alpha^2$,} \\[2mm]
\sqrt{\mu-\alpha^2} \, \coth(\sqrt{\mu-\alpha^2}), & \mbox{if $\mu \geq \alpha^2$.}
\end{cases}$$
The calculation
$$g(\bfk)=-\frac{1}{1+\tfrac{k_2^2}{k_1^2}}\alpha\big(c_1\tfrac{k_2}{k_1}-c_2\big)\big(c_1+c_2\tfrac{k_2}{k_1}\big)
-\frac{1}{1+\tfrac{k_2^2}{k_1^2}}\mathtt{c}\big(k_1^2(1+\tfrac{k_2^2}{k_1^2})\big)\big(c_1+c_2\tfrac{k_2}{k_1}\big)^2
+1+\beta k_1^2(1+\tfrac{k_2^2}{k_1^2})$$
shows that $g$ is an analytic function $\tilde{g}$ of $k_1$ and $\frac{k_2}{k_1}$. We find that $\tilde{g}(0,0)=0$ if $\bfc=\bfc_0$, where
$$\bfc_0=\begin{pmatrix}c_0\cos\tfrac{1}{2}\alpha \\[1mm] -c_0\sin\tfrac{1}{2}\alpha\end{pmatrix}, \qquad c_0^2=\tfrac{2}{\alpha} \tan \tfrac{1}{2}\alpha,$$
and it is shown in Appendix \ref{disprel} that the dispersion relation \eqref{dispersion relation - intro} has no further solutions for sufficiently large values of $\beta$. Substituting the Ansatz
$$\bfc=(1-\varepsilon^2)\bfc_0$$
and
\begin{equation}
\eta(x,y)=\varepsilon^2 \zeta (X,Y), \qquad X=\varepsilon x,\ Y=\varepsilon^2 y \label{KP scaling}
\end{equation}
into equations \eqref{Beltrami 1}--\eqref{Beltrami 5}, one duly finds that to leading order $\zeta$ satisfies the KP-I equation
\begin{equation}
- (\beta-\beta_0) \zeta_{xx} +2 \zeta +\sec^2\!\tfrac{1}{2}\alpha \frac{D_2^2}{D_1^2}\zeta
+ d_\alpha \zeta^2=0, \label{Physical KP}
\end{equation}
where
$$d_\alpha=\alpha \cosec \alpha + \tfrac{1}{2}\alpha \cot \alpha,$$
we have replaced $(X,Y)$ with $(x,y)$ for notational simplicity and $D_1=-\i \partial_x$, $D_2=-\i \partial_y$.

\begin{figure}[h]
\begin{center}
\includegraphics[scale=0.4]{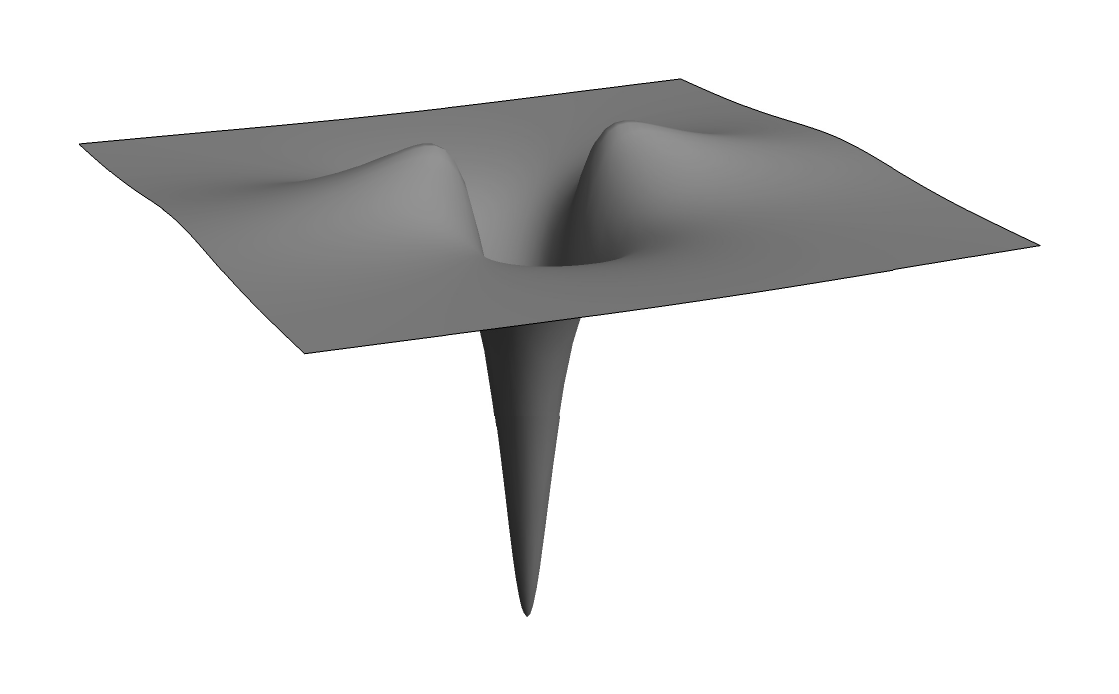}\hspace{5mm}
\includegraphics[scale=0.4]{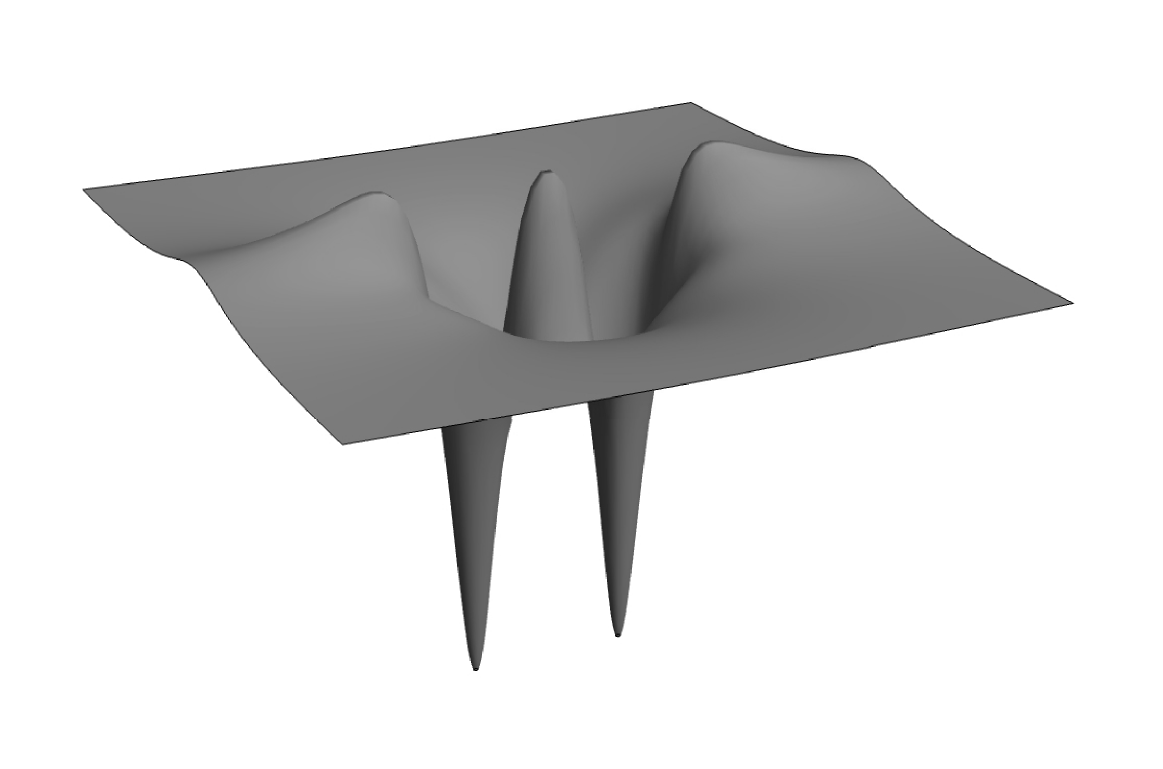}
\caption{The KP lumps $\zeta_1^\star$ (left) and $\zeta_2^\star$ (right).}\label{lumps}
\end{center}
\end{figure}

Using the change of variable
$$u(\tilde{x},\tilde{y})=\frac{1}{6}d_\alpha \zeta(x,y), \quad x=\left(\frac{\beta-\beta_0}{2}\right)^{\!\!\frac{1}{2}}\tilde{x},\quad y=\frac{1}{2}\big((\beta-\beta_0)\sec^2\tfrac{1}{2}\alpha\big)^\frac{1}{2} \tilde{y}$$
we can write \eqref{Physical KP} in the normalised form
\begin{equation}
\partial_x^2 (-\partial_x^2 u +u +3 u^2)+\partial_y^2 u =0, \label{Normalised KP}
\end{equation}
where we have dropped the tildes for notational simplicity.
Equation \eqref{Normalised KP} has a family of explicit symmetric `lump' solutions of the form
\begin{equation}
u_k^\star(x,y)=-2\partial_x^2 \log \tau_k^\star(x,y),\qquad k=1,2,\ldots, \label{log lump}
\end{equation}
where $\tau_k^\star$ is a polynomial of degree $k(k+1)$ with $\tau_k^\star(x,y)=\tau_k(-x,y)=\tau_k^\star(x,-y)$ for all $(x,y) \in {\mathbb R}^2$;
the first two members of the family are
\begin{align*}
\tau_1^\star(x,y) &=x^2+y^2+3, \\
\tau_2^\star(x,y) &=x^6+3x^4y^2+3x^2y^4+y^6+25x^4+90x^2y^2+17y^4-125x^2+475y^2+1875.
\end{align*}
Note that the lump solutions $u_k^\star$ are smooth, decaying rational functions, so that the same is true of their
derivatives of all orders.
The functions $\zeta_1^\star$ and $\zeta_2^\star$ (where $\zeta_k^\star$ is obtained from $u_k^\star$ by reversing the normalisation) are sketched in Figure \ref{lumps}.

The following result was established by Liu \& Wei \cite{LiuWei19} and
Liu, Wei \& Yang \cite{LiuWeiYang24a,LiuWeiYang24b}.\pagebreak

\begin{lemma}$ $ \label{LWY results}
\begin{itemize}
\item[(i)] Every smooth, algebraically decaying lump solution of \eqref{Normalised KP} has the form $u(x,y)=-2\partial_x^2 \log \tau(x,y)$
for some polynomial $\tau$ of degree $k(k+1)$ with $k \in {\mathbb N}$ and satisfies $|u(x,y)| \lesssim (1+x^2+y^2)^{-1}$ for all $(x,y) \in {\mathbb R}^2$.
\item[(ii)] There is a unique symmetric lump solution of the form \eqref{log lump} for each $k \in {\mathbb N}$ with $k(k+1) \leq 600$
(and it is conjectured that this result holds for all $k \in {\mathbb N}$).
\item[(iii)] The lump solutions $\zeta_k^\star$, $k=1,2$ of \eqref{Normalised KP} are nondegenerate in the sense that the only smooth solutions
of the linearised equation
$$\partial_x^2(-\partial_x^2 u +u +6 u_k^\star u)+\partial_y^2u=0$$
with $u(x,y) \rightarrow 0$ as $x^2+y^2 \rightarrow \infty$
are linear combinations of $\partial_x u_k^\star$ as $\partial_y u_k^\star$ (and it is conjectured that this result holds for all $k \in {\mathbb N}$; see Remark \ref{More lumps} below).
\end{itemize}
\end{lemma}

The KP lump solution $\zeta_k$ formally corresponds to a fully localised solitary water wave $\eta(x,z)=\varepsilon^2 \zeta_k(\varepsilon x, \varepsilon^2 y)$.
In this article we rigorously reduce the hydrodynamic equations \eqref{Beltrami 1}--\eqref{Beltrami 5} to a perturbation of the KP-I equation and combine the nondegeneracy result in Lemma \ref{LWY results}(iii) with an implicit-function theorem argument to establish the following result.

\begin{theorem} \label{Main result}
Suppose that $|\alpha| \in [0,\tfrac{1}{2}\pi)$ and
$$c_1=c_0(1-\varepsilon^2)\cos\tfrac{1}{2}\alpha, \qquad c_2=-c_0(1-\varepsilon^2)\sin\tfrac{1}{2}\alpha$$
where
$$c_0^2=\tfrac{2}{\alpha} \tan \tfrac{1}{2}\alpha.$$
For each sufficiently large value of $\beta>0$ and each sufficiently small value of
$\varepsilon>0$ equations \eqref{Beltrami 1}--\eqref{Beltrami 5} possess fully localised solitary-wave solutions
$\eta_1^\star$, $\eta_2^\star \in H^3({\mathbb R}^2)$ which satisfy
$\eta_k^\star(x,y)=\eta_k^\star(-x,-y)$ for all $(x,y) \in {\mathbb R}^2$ and
\begin{equation}
\eta_k^\star(x,y)=\varepsilon^2 \zeta_k^\star(\varepsilon x,\varepsilon^2 y)+ o(\varepsilon^2) \label{uniform KP approx}
\end{equation}
uniformly over $(x,y) \in {\mathbb R}^2$.
\end{theorem}
\begin{remark} \label{More lumps}
In fact Theorem \ref{Main result} generates a fully localised
solitary water wave from any symmetric lump solution $\zeta_k^\star$ of \eqref{Physical KP} which is nondegenerate in
the sense of Lemma \ref{LWY results}(iii), and a sketch of the proof of the nondegeneracy of
$\zeta_k^\star$ for $k \geq 3$ was given by Liu, Wei \& Yang \cite{LiuWeiYang24b}.
\end{remark}

\subsection{Reformulation}

We proceed using a recent formulation of \eqref{Beltrami 1}--\eqref{Beltrami 5} due to Groves \& Horn \cite{GrovesHorn20} which generalises the
Zakharov-Craig-Sulem formulation of the irrotational problem (Zakharov \cite{Zakharov68}, Craig \& Sulem \cite{CraigSulem93}).
For a vector field\linebreak $\bfF=(F_1,F_2,F_3)^\mathrm{T}$ let
$$\bfF_\parallel = \bfF_\mathrm{h}+F_3\nabla \eta\big|_{z=\eta},$$\pagebreak
where $\bfF_\mathrm{h}=(F_1,F_2)^\mathrm{T}$,
and write, according to the Hodge-Weyl decomposition for vector fields in two-dimensional free space,
$$
\bfv_\parallel = \nabla \Phi + \nabla^\perp \Psi,
$$
where $\Phi = \Delta^{-1} (\nabla \cdot \bfv_\parallel)$, $\Psi=\Delta^{-1}(\nabla^\perp\cdot\bfv_\parallel)
=-\Delta^{-1}(\nabla \cdot \bfv_\parallel^\perp)$ and $\Delta^{-1}$ is the two-dimensional Newtonian potential.
Define a \emph{generalised Dirichlet--Neumann operator} $H(\eta)$ by
$$
H(\eta)\Phi = \udl{\curl \bfA} \cdot \bfn=\nablac \bfA_\parallel^\perp,
$$
where $(f_1,f_2)^\perp=(f_2,-f_1)$,  the underscore denotes evaluation at $z=\eta$ and $\bfA$ is the unique solution of the boundary-value problem
\begin{alignat}{2}
\curl \curl \bfA &= \alpha  \curl \bfA & & \mbox{in $D_\eta$,}  \label{A BVP 1} \\
\Div \bfA &= 0 & & \mbox{in $D_\eta$,}  \label{A BVP 2} \\
\bfA \wedge \bfe_3 &= \bfzero & & \mbox{at $z=-1$,}  \label{A BVP 3} \\
\bfA \cdot \bfn &= 0 & & \mbox{at $z=\eta$,}  \label{A BVP 4} \\
(\curl \bfA)_\parallel &= \nabla\Phi - \alpha \nablap\Delta^{-1}(\nablac \bfA_\parallel^\perp)\qquad & & \mbox{at $z=\eta$.} \label{A BVP 5}
\end{alignat}
(Note that $\Psi=\Delta^{-1}(\nablapc (\curl \bfA)_{\parallel})$ is necessarily given by $\Psi = -\alpha\, \Delta^{-1}(\nablac\bfA_\parallel^\perp)$ because
\begin{equation}
\Psi \!=\! -\Delta^{-1}(\nablac\curl \bfA^{\perp}_\parallel)\!=\!-\Delta^{-1}(\udl{\curl\curl \bfA}\cdot\bfn)
\!=\!-\alpha\, \Delta^{-1}( \udl{\curl \bfA}\cdot\bfn)\!=\!- \alpha\, \Delta^{-1} (\nablac\bfA^{\!\perp}_\parallel),
\label{compatibility}
\end{equation}
and that $\bfv=\curl \bfA$ satisfies \eqref{Beltrami 1}--\eqref{Beltrami 3}.)

A straightforward calculation shows that
equations \eqref{Beltrami 4}--\eqref{Beltrami 5} are equivalent to
\begin{align}
& H(\eta)\Phi + \udl{\bfu}^\star \cdot \bfn=0, \label{GH1} \\
& \frac{1}{2} |\bfK(\eta)\Phi|^2\! -\! \frac{ ( H(\eta)\Phi + \bfK(\eta)\Phi \!\cdot\! \nabla\eta)^2 }{2(1+|\nabla\eta|^2)}  \nonumber \\
&\qquad\quad\mbox{} + \bfK(\eta)\Phi \cdot \udl{\bfu}_\mathrm{h}^\star+ \eta - \beta \left( \frac{\eta_x}{ (1+|\nabla\eta|^2)^{\frac{1}{2}} } \right)_{\!\!x}\!\!\! - \beta \left( \frac{\eta_y}{ (1+|\nabla\eta|^2)^{\frac{1}{2}} } \right)_{\!\!y}\!\!\!= 0, \label{GH2}\hspace{-0.35cm}
\end{align}
where
$$
\bfK(\eta)\Phi \coloneqq  \nabla\Phi - \alpha  \nablap \Delta^{-1}(H(\eta)\Phi),
$$
and these equations can in fact be reduced to a single equation for the variable $\eta$ (see Groves \& Horn \cite{GrovesHorn20} for full details and Oliveras \& Vashal \cite{OliverasVasan13} for a simpler version of this equation for irrotational waves). Equation \eqref{GH1} implies that
$\Phi=-H(\eta)^{-1}(\udl{\bfu}^\star \cdot \bfn)$, whereby \eqref{GH2} yields
\begin{equation}
\JJ(\eta)=0,  \label{J eqn}
\end{equation}
where
\begin{align}
\JJ(\eta)\coloneqq \frac{1}{2} |\bfT(\eta)|^2 &- \frac{ ( -\udl{\bfu}^\star \cdot \bfn + \bfT(\eta) \cdot \nabla\eta)^2 }{2(1+|\nabla\eta|^2)} \nonumber \\
&\mbox{}+ \bfT(\eta) \cdot \udl{\bfu}_\mathrm{h}^\star + \eta- \beta \left( \frac{\eta_x}{ (1+|\nabla\eta|^2)^{\frac{1}{2}} } \right)_{x} -\beta \left( \frac{\eta_y}{ (1+|\nabla\eta|^2)^{\frac{1}{2}} } \right)_{y} \label{Defn of J}
\end{align}
and
$$
\bfT(\eta) \coloneqq  - \nabla \left( H(\eta)^{-1} (\udl{\bfu}^\star \cdot \bfn) \right) + \alpha \, \nablap \Delta^{-1} (\udl{\bfu}^\star \cdot \bfn).
$$
Equation \eqref{J eqn} is invariant under the transformation $\eta(x,y) \mapsto \eta(-x,-y)$ (see the discussion beneath equations \eqref{Beltrami 1}--\eqref{Beltrami 5}), and in
this paper we show that \eqref{J eqn} has solutions $\eta_1^\star$, $\eta_2^\star \in H^3({\mathbb R}^2)$ which satisfy the estimate \eqref{uniform KP approx}
and are invariant under this transformation.

The operator $\bfT(\eta)$ can also be defined directly in terms of a boundary-value problem. Noting that\linebreak
$\udl{\bfu}^\star \cdot \bfn = \nablac \bfS(\eta)^\perp$ (and, for later use, that $\udl{\bfu}_\mathrm{h}^\star = \bfc+\alpha\bfS(\eta)$), where
$$
\bfS(\eta) =  \frac{c_1}{\alpha} 
\begin{pmatrix}
\cos (\alpha \, \eta) -1 \\
-\sin (\alpha \, \eta)
\end{pmatrix}
+\frac{c_2}{\alpha} 
\begin{pmatrix}
\sin (\alpha \, \eta) \\
\cos (\alpha \, \eta) -1
\end{pmatrix},
$$
we can define
$$
\bfT(\eta) \coloneqq  \bfM(\eta)\bfS(\eta), 
$$
where
$$
\bfM(\eta)\bfg \coloneqq  -(\curl \bfB)_{\parallel}, \label{eq:MOpDef}
$$
and $\bfB$ solves the boundary-value problem
\begin{alignat}{2}
\curl \curl \bfB &= \alpha\curl \bfB \qquad\quad& & \mbox{in $D_\eta$}, \label{B BVP 1} \\
\Div \bfB &= 0 & & \mbox{in $D_\eta$}, \label{B BVP 2} \\
\bfB \wedge \bfe_{3} &= \bfzero & & \mbox{at $z=-1$}, \label{B BVP 3} \\
\bfB \cdot \bfn &= 0 & & \mbox{at $z=\eta$,} \label{B BVP 4} \\
\nablac \bfB_{\parallel}^{\perp} &= \nablac \bfg^\perp & & \mbox{at $z=\eta$}. \label{B BVP 5}
\end{alignat}
Any solution to this boundary-value problem satisfies
$$(\curl \bfB)_\parallel = \nabla \Phi - \alpha \nablap \Delta^{-1}(\nablac \bfB_\parallel^\perp)$$
 for some $\Phi$ (see equation \eqref{compatibility}), so that $\Phi=H(\eta)^{-1}(\nablac \bfg^\perp)$
and
$$-(\curl \bfB)_{\parallel}=- \nabla (H(\eta)^{-1}(\nablac \bfg^\perp)) + \alpha \nablap \Delta^{-1}(\nablac \bfg^\perp).$$

In Section \ref{anal} we show that the solutions to the boundary-value problems \eqref{A BVP 1}--\eqref{A BVP 5} and \eqref{B BVP 1}--\eqref{B BVP 5}
depend analytically upon $\eta$ and use this fact to deduce that the same is true of $H(\eta)$ and $\bfM(\eta)$.
We proceed by `flattening' the fluid domain by means
of the transformation $\Sigma\colon D_0 \to D_\eta$ given by
$$
\Sigma\colon(\bfx,v) \mapsto (\bfx,v+\sigma(\bfx,v)), \qquad \sigma(\bfx,v)\coloneqq  \eta(\bfx)(1+v)
$$
and differential operators
$$
\partial_x^\sigma \coloneqq  \partial_x - \frac{ \partial_x\sigma }{1+\partial_v\sigma} \partial_v, \quad
\partial_y^\sigma \coloneqq  \partial_y - \frac{ \partial_y\sigma }{1+\partial_v\sigma} \partial_v, \quad
\partial_v^\sigma \coloneqq  \frac{ \partial_v }{1+\partial_v\sigma}
$$
(with corresponding notation for differential operators derived from these partial derivatives).
The boundary-value problems for $\bfA^\prime$ and $\bfB^\prime$ are transformed into equivalent problems
for $\bfA^\prime \coloneqq \bfA \circ \Sigma$ and
$\bfB^\prime\coloneqq  \bfB \circ \Sigma$ in the fixed domain $D_0$ (equations \eqref{pre flattened A BVP 1}--\eqref{pre flattened A BVP 5} and
\eqref{pre flattened B BVP 1}--\eqref{pre flattened B BVP 5} respectively), for which we establish the following results
(the function spaces $\ZZ$, $\dot{H}^s({\mathbb R}^2)$ and
$\check{H}^s({\mathbb R}^2)$ are defined in Sections \ref{Function spaces} and \ref{Anal BVPs} below).

\begin{theorem}
Suppose that $|\alpha| \in [0,\tfrac{1}{2}\pi)$. There exists a neighbourhood $V$ of the origin in $\ZZ$ for which the following statements hold.
\begin{itemize}
\item[(i)]
The boundary-value problem \eqref{pre flattened A BVP 1}--\eqref{pre flattened A BVP 5} 
has a unique solution $\bfA^\prime \in H^\frac{5}{2}(D_0)^3$ which depends analytically upon
$\eta \in V$ and
$\Phi \in \dot{H}^{\frac{5}{2}}({\mathbb R}^2)$; furthermore
$$\bfA^\prime(\eta,\Phi) = \Grad^\sigma \phi(\eta,\tilde\bfA) + \tilde\bfA(\eta,\Phi),$$
where $\phi \in H^\frac{7}{2}(D_0)$ depends analytically upon $\eta \in V$ and $\tilde\bfA \in H^3(D_0)^3$,
while $\tilde\bfA \in H^3(D_0)^3$ depends analytically upon $\eta \in V$ and $\Phi \in \dot{H}^{\frac{5}{2}}({\mathbb R}^2)$.
\item[(ii)]
The boundary-value problem \eqref{pre flattened B BVP 1}--\eqref{pre flattened B BVP 5} has a unique solution
$\bfB^\prime \in H^\frac{5}{2}(D_0)^3$ which depends analytically upon
$\eta \in V$ and $\bfg \in H^{\frac{5}{2}}({\mathbb R}^2)^2$; furthermore
$$\bfB^\prime(\eta,\Phi) = \Grad^\sigma \phi(\eta,\tilde\bfB) + \tilde\bfB(\eta,\Phi),$$
where $\phi \in H^\frac{7}{2}(D_0)$ depends analytically upon $\eta \in V$ and $\tilde\bfB \in H^3(D_0)^3$,
while $\tilde\bfB \in H^3(D_0)^3$ depends analytically upon $\eta \in V$ and $\bfg \in H^{\frac{5}{2}}({\mathbb R}^2)^2$.
\end{itemize}
\end{theorem}

The analyticity of $H$ and $\bfM$ follows from the above theorem and the facts that
$$
H(\eta)(\Phi)=\nablac \bfA_\parallel^{\prime\perp}=\nablac \tilde\bfA_\parallel^\perp,
\qquad
\bfM(\eta)(\bfg)=-(\curl^\sigma \bfB^\prime)_\parallel=-(\curl^\sigma \tilde\bfB)_\parallel.
$$

\begin{theorem}
Suppose that $|\alpha| \in [0,\tfrac{1}{2}\pi)$.
The mappings $\eta \mapsto H(\eta)$ and $\eta \mapsto \bfM(\eta)$ are analytic\linebreak
$V \to L(\dot{H}^{\frac{5}{2}}({\mathbb R}^2),\check{H}^{\frac{3}{2}}({\mathbb R}^2))$ and
$V \to L(H^{\frac{5}{2}}(\mathbb{R}^2)^2, H^{\frac{3}{2}}(\mathbb{R}^2)^2)$ respectively.
\end{theorem}

Our final result follows by noting that $H^3({\mathbb R}^2)$ is continuously embedded in $\ZZ$, so that
\begin{equation}
U\coloneqq \{\eta \in H^3({\mathbb R}^2)\colon \|\eta\|_\ZZ<M\} \label{Defn of U with M}
\end{equation}
is an open neighbourhood of the origin in $H^3({\mathbb R}^2)$.
\begin{corollary}
Suppose that $|\alpha| \in [0,\tfrac{1}{2}\pi)$.
The formula \eqref{Defn of J} defines an analytic function $\JJ\colon U \rightarrow H^1({\mathbb R}^2)$ for sufficiently small values of $M>0$.
\end{corollary}

\subsection{Reduction}

\begin{figure}[h]
\centering

\includegraphics[scale=0.8]{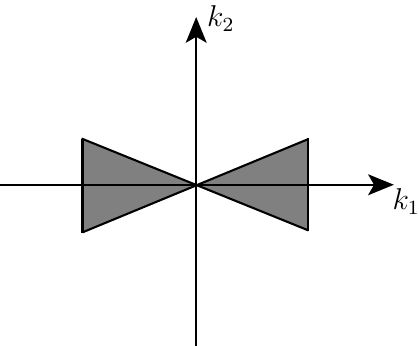}
{\it
\caption{The set $\displaystyle S=\left\{(k_1,k_2)\colon |k_1|, \Big|\frac{k_2}{k_1}\Big| \leq \delta\right\}$. \label{Splitting}}}
\end{figure}

The Ansatz \eqref{KP scaling} indicates that the Fourier transform of the surface-profile function
$\eta$ for a fully localised solitary wave is concentrated in the region
$$S=\left\{(k_1,k_2)\colon |k_1|, \Big|\frac{k_2}{k_1}\Big| \leq \delta\right\}$$
for some $0< \delta \ll 1$ (see Figure \ref{Splitting}). We therefore decompose $\eta \in H^3({\mathbb R}^2)$ into the sum of
$$\eta_1 =\chi(\bfD)\eta, \qquad
\eta_2 = \left(1-\chi(\bfD)\right)\eta,$$
where $\chi$ is the indicator function of the set $S$,
the Fourier transform $\hat{\eta}=\FF[\eta]$ of $\eta$ is defined by
$$\hat{\eta}(k_1,k_2)=\frac{1}{2\pi}\int_{{\mathbb R}^2}\eta(x,y)\ee^{-\ii(k_1x + k_2y)}\dx\dy,$$
and 
$\bfD=(-\i \partial_x, -\i \partial_y)^T$.
Setting
$$\bfc=(1-\varepsilon^2)\bfc_0,$$
choosing $\beta$ sufficiently large and writing equation \eqref{J eqn} as
\begin{align*}
\chi(\bfD)\JJ(\eta_1+\eta_2) & = 0,\\
\left(1-\chi(\bfD)\right)\JJ(\eta_1+\eta_2) & = 0,
\end{align*}
we find that the second equation is solvable for $\eta_2$ as a function of $\eta_1$ for sufficiently
small values of $\varepsilon$; the first therefore reduces to
$$
\chi(\bfD)\JJ(\eta_1+\eta_2(\eta_1)) = 0
$$
upon inserting $\eta_2=\eta_2(\eta_1)$. Finally, the scaling
$$
\eta_1(x,y)=\varepsilon^2 \zeta (X,Y), \qquad X=\varepsilon x,\ Y=\varepsilon^2 y
$$
transforms the reduced equation into a perturbation of the equation
\begin{equation}
\varepsilon^{-2}g_\varepsilon(\bfD)\zeta
+2\zeta
+ d_\alpha\chi_\varepsilon(\bfD)\zeta^2=0, \label{Red eq - intro}
\end{equation}
where $g_\varepsilon(k_1,k_2)=g(\varepsilon k_1, \varepsilon^2 k_2)$ and $\chi_\varepsilon(k_1,k_2)=\chi(\varepsilon k_1, \varepsilon^2 k_2)$,
such that $\chi_\varepsilon$ is the indicator function of the set
$$
S_\varepsilon=\left\{ \bfk \in {\mathbb R}^2 \colon |k_1| \leq \frac{\delta}{\varepsilon}, \Big|\frac{k_2}{k_1}\Big|\leq \frac{\delta}{\varepsilon}\right\}
$$
(see Sections \ref{reduction} and \ref{calculation of reduced eqn}; the reduced equation is stated precisely in equation \eqref{Final reduced eqn}).
Note that $\delta$ is a small, but fixed constant while $\varepsilon$ is a small parameter whose maximum value depends upon $\delta$.

Equation \eqref{Red eq - intro} is a \emph{full-dispersion}
version of the stationary KP-I equation \eqref{Physical KP} since it retains the linear part of the original equation
\eqref{J eqn};  noting that
$$\varepsilon^{-2}g_\varepsilon(k_1,k_2) =(\beta-\beta_0)k_1^2 +\sec^2\!\tfrac{1}{2}\alpha \frac{k_2^2}{k_1^2} + O(\varepsilon),$$
we recover the fully reduced model equation
in the formal limit $\varepsilon=0$.
In Section \ref{sec:existence} we demonstrate that equation \eqref{Red eq - intro} for $\zeta$ has solutions $\zeta_1^\varepsilon$, $\zeta_2^\varepsilon$
which satisfy $\zeta_k^\varepsilon \rightarrow \pm \zeta_k^\star$ as $\varepsilon \rightarrow 0$
in a suitable function space (see Theorem \ref{Final strong existence thm}). The key step is the nondegeneracy result given in Lemma \ref{LWY results}(iii)
which allows one to apply a suitable variant of the implicit-function theorem. For this purpose we
exploit the fact that the reduction procedure preserves the invariance
of equation \eqref{J eqn} under
$\eta(x,y) \mapsto \eta(-x,-y)$, so that equation \eqref{Red eq - intro} is
invariant under $\zeta(x,y) \mapsto \zeta(-x,-y)$; restricting to a space of functions with this invariance,
we find that the kernel of the appropriate linearisation is trivial since $\partial_x \zeta_k^\star$, $\partial_y \zeta_k^\star$ do not have this invariance.

The perturbation argument used in Section \ref{sec:existence} was developed by Groves \cite{Groves21} in the context of two-dimensional
irrotational solitary waves and applied to three-dimensional irrotational fully localised solitary waves on water of infinite depth by Buffoni, Groves \& Wahl\'{e}n
\cite{BuffoniGrovesWahlen22}. It has also been applied to the Whitham equation by Stefanov \& Wright \cite{StefanovWright20} and to a full dispersion KP-I equation (which differs from \eqref{Red eq - intro}) by Ehrnstr\"{o}m \& Groves \cite{EhrnstroemGroves26}.

\subsection{Function spaces} \label{Function spaces}

We work with the standard function spaces $H^n(D_0)$ for $n \in {\mathbb N}_0$ (and $H^s(D_0)$ for $s \geq 0$ defined by interpolation)
in the flattened fluid domain
together with $L^p({\mathbb R}^2)$ for $p \geq 1$, $W^{n,\infty}({\mathbb R}^2)$ for $n \in {\mathbb N}_0$
and
$$H^s({\mathbb R}^2) = \{\eta \in L^2({\mathbb R}^2)\colon (1+|\bfk|^2)^{\frac{1}{2}s} \hat{\eta} \in L^2({\mathbb R}^2)\},\qquad \|\eta\|_s^2 = \int_{{\mathbb R}^2} (1+|\bfk|^2)^s |\hat{\eta}(\bfk)|^2 \dbfk$$
for $s \geq 0$ in the plane; the definitions are extended componentwise to vector-valued functions.
The nonstandard spaces 
\begin{align*}
\dot{H}^s({\mathbb R}^2) &= \{\eta \in L^2_\mathrm{loc}({\mathbb R}^2)\colon \nabla \eta \in H^{s-1}({\mathbb R}^2)^2\}/{\mathbb R}, \quad \|\eta\|_{\dot{H}^s} \coloneqq  \|\nabla \eta\|_{s-1}
\cong\|\langle\bfk\rangle^{s-1}|\bfk|\hat{\eta}\|_0, \qquad s \geq 1,\\
\check{H}^s({\mathbb R}^2) &= \{\eta \in L^2({\mathbb R}^2)\colon \Delta^{-1} \eta \in \dot{H}^{s+2}({\mathbb R}^2)\}, \quad \|\eta\|_{\check{H}^s} \coloneqq  \|\Delta^{-1} \eta\|_{\dot{H}^{s+2}}
\cong\|\langle \bfk \rangle^{s+1} |\bfk|^{-1}\hat{\eta}\|_0, \qquad s \geq 0,
\end{align*}
where $L^2_\mathrm{loc}({\mathbb R}^2)$ denotes the space of locally square integrable functions in the plane, and the scale $\{Y_s,\|\cdot\|_s\}_{s \geq 0}$, where
$$Y_s = \left\{\eta \in L^2({\mathbb R}^2)\colon m^s\hat{\eta} \in L^2({\mathbb R}^2)\right\}, \quad
\|\eta\|_{Y_s} \coloneqq  \|m^s\hat{\eta}\|_0,
$$
and
$$m(\bfk) = \left( 1 + k_1^2 + \frac{k_2^2}{k_1^2} \right)^\frac{1}{2}$$
are also used. Note that $\Delta \colon \dot{H}^{s+2}({\mathbb R}^2) \to H^s({\mathbb R}^2)$ is injective so that $\check{H}^s({\mathbb R}^2)$ is well defined.

\begin{proposition} \label{Properties of Y} $ $
\begin{itemize}
\item[(i)] The space $Y_1$ is continuously embedded in $L^p({\mathbb R}^2)$ for $2 \leq p < 6$.
\item[(ii)] The space $Y_1$ is compactly embedded in $L^2(B_R(\bfzero))$ for each $R >0$.
\item[(iii)] The space $Y_s$ is continuously embedded in $C_\mathrm{b}({\mathbb R}^2)\coloneqq C({\mathbb R}^2)\cap L^\infty({\mathbb R}^2)$ for $s > \frac{3}{2}$.
\end{itemize}
\end{proposition}
{\bf Proof.}  $ $
\begin{itemize}
\item[(i)]
Suppose $p \geq 2$. Using the Hausdorff--Young and H\"{o}lder inequalities we find that
$$\|f\|_{L^p({\mathbb R}^2)} \leq \|\hat{f}\|_{L^q({\mathbb R}^2)}=\|m^{-1}m\hat{f}\|_{L^q({\mathbb R}^2)} \leq \|m^{-1}\|_{{L^r({\mathbb R}^2)}} \|m \hat{f}\|_{0} = I^\frac{1}{r}\|f\|_Y,$$
where $\frac{1}{p} + \frac{1}{q} = 1$, $\frac{1}{q} = \frac{1}{r} + \frac{1}{2}$ and
$$I=\int_{{\mathbb R}^2} m(\bfk)^{-\frac{1}{2}r}\dbfk = \int_{{\mathbb R}^2} (1+|\bfk|^2)^{-\frac{1}{2}r}|k_1| \dbfk < \infty$$
if $r>3$, that is $p<6$.

\item[(ii)]
Let $\chi_n$ be the indicator function of the set $B_n(\bfzero)$ and define bounded linear operators\linebreak
$T, T_n\colon Y_1 \rightarrow L^2(B_R(\bfzero))$ by
$$Tf = f|_{B_R(\bfzero)}, \qquad T_nf = T\chi_n(\bfD)f.$$
The calculation
$$
\|T_n \eta\|_{H^1(B_R(\bfzero))}^2
\leq \|\chi_n(\bfD)\eta\|_{1}^2
= \!\!\int_{|\bfk|<n} (1 + |\bfk|^2) |\hat{\eta}(\bfk)|^2 \dbfk
\leq (1 + n^2)\!\! \int_{{\mathbb R}^2} |\hat{\eta}(\bfk)|^2 \dbfk
\leq (1+n^2) \|\eta\|_{Y_1}^2
$$
shows that $T_n$ is a bounded linear operator $Y_1 \rightarrow H^1(B_R(\bfzero))$. Since
$H^1(B_R(\bfzero))$ is compactly embedded in $L^2(B_R(\bfzero))$ by the Rellich--Kondrachov
theorem we conclude that $T_n\colon Y_1 \rightarrow L^2(B_R(\bfzero))$ is compact.

Furthermore
\begin{align*}
\|T f - T_n f\|_{L^2(B_R(\bfzero))}^2 & \leq \|f - \chi_n(\bfD) f\|_{0}^2 \\
&= \int_{|\bfk| > n} |\hat{f}(\bfk)|^2 \dbfk \\
& = \int_{|\bfk| > n} \frac{1}{m(\bfk)^2} m(\bfk)^2 |\hat{f}(\bfk)|^2 \dbfk \\
& \leq \left( \sup_{|\bfk| > n} \frac{1}{m(\bfk)^2} \right) \int_{|\bfk| > n} m(\bfk)^2 |\hat{f}(\bfk)|^2 \dbfk \\
& \leq \left( \sup_{|\bfk| > n} \frac{1}{m(\bfk)^2} \right)\|f\|_{Y_1}^2,
\end{align*}
such that
$$\|T-T_n\|_{L(Y_1,L^2(B_R(\bfzero)))} \leq \sup_{|\bfk| > n} \frac{1}{m(\bfk)^2}\rightarrow 0$$
as $n \rightarrow \infty$ (because $m(\bfk) \rightarrow \infty$ as $|\bfk| \rightarrow \infty$; it follows that
$T$ is compact (as the uniform limit of the sequence $\{T_n\}$ of compact operators).
\item[(iii)]
Note that
$$
\|\eta\|_\infty \leq \|\hat{\eta}\|_{L^1({\mathbb R}^2)} = \| m^{-\frac{1}{2}s}m^{\frac{1}{2}s} \hat{\eta}\|_{L^1({\mathbb R}^2} 
\leq  I^\frac{1}{2}\|\eta\|_{Y_s},
$$
where
$$I=\int_{{\mathbb R}^2} m(\bfk)^{-s}\dbfk = \int_{{\mathbb R}^2} (1+|\bfk|^2)^{-s} |k_1|\dbfk < \infty$$
if $s>\frac{3}{2}$. The continuity of $\eta$ follows from a standard dominated convergence argument.\qed
\end{itemize}

Observe that the spaces $\chi(\bfD)H^s({\mathbb R}^2)$ and $\chi(\bfD)Y_s$, $s \geq 0$
of `truncated' functions all coincide and have equivalent norms. In Sections \ref{reduction} and \ref{calculation of reduced eqn} we 
identify in particular $\chi(\bfD)H^3({\mathbb R}^2)$ with $\chi(\bfD)Y_1$ and equip it with the scaled norm
\begin{equation}
\label{scaled norm}
\nn \eta \nn^2\coloneqq  \int_{{\mathbb R}^2} \left(1+\varepsilon^{-2}k_1^2+\varepsilon^{-2}\sdfrac{k_2^2}{k_1^2}\right)
|\hat \eta(\bfk)|^2\dbfk
\end{equation}
in anticipation of the KP scaling.

\begin{proposition}
\label{sup estimate}
The estimate 
$\|\hat \eta_1\|_{L^1({\mathbb R}^2)} \lesssim
\varepsilon \nn \eta_1\nn$
holds for each $\eta_1\in \chi(\bfD)Y_1$.
\end{proposition}
{\bf Proof.}
Observe that
$$
\int_{{\mathbb R}^2} |\hat \eta_1(\bfk)|\dbfk
= \int_{{\mathbb R}^2}\left(1+\varepsilon^{-2}k_1^2 + \varepsilon^{-2}\tfrac{k_2^2}{k_1^2}\right)^{-\frac{1}{2}}
\left(1+\varepsilon^{-2}k_1^2 + \varepsilon^{-2}\tfrac{k_2^2}{k_1^2}\right)^{\frac{1}{2}}
 |\hat \eta_1(\bfk)| \dbfk \lesssim \nn \eta_1 \nn I^{\frac{1}{2}},
$$
where
\begin{align*}
I &= \int_S \frac{1}{1+\varepsilon^{-2}k_1^2 + \varepsilon^{-2}\tfrac{k_2^2}{k_1^2}} \dbfk \\
&= 4\varepsilon^3 \int_0^{\delta/\varepsilon} \int_0^{\delta/\varepsilon} \frac{u_1}{1+u_2^2+u_1^2} \du_2\du_1 \\
&= 4\varepsilon^3 \int_0^{\delta/\varepsilon} \frac{u_1}{(1+u_1^2)^\frac{1}{2}}\arctan \frac{\delta/\varepsilon}{(1+u_1^2)^\frac{1}{2}}\du_1 \\
& \leq 2\pi\varepsilon^3 \int_0^{\delta/\varepsilon} \frac{u_1}{(1+u_1^2)^\frac{1}{2}}\du_1 \\
& = 2\pi\varepsilon^3 \left(\big(1+\tfrac{\delta^2}{\varepsilon^2}\big)^\frac{1}{2}-1\right) \\
& \lesssim \varepsilon^2. \alignqed
\end{align*}

\begin{corollary}
The estimate $\|\eta_1\|_{n,\infty} \lesssim \varepsilon \nn \eta_1 \nn$ holds for each $\eta_1 \in \chi(\bfD)Y_1$ and each $n \in {\mathbb N}_0$.
\end{corollary}
{\bf Proof.} The result follows from the calculation $\|\eta_1\|_{n,\infty} \lesssim \| |\bfk|^n \hat{\eta}\|_{L^1({\mathbb R}^2)} \lesssim  \|\hat \eta_1\|_{L^1({\mathbb R}^2)}$
and the previous proposition.\qed

Finally, we introduce the space $Y_s^\varepsilon=\chi_\varepsilon(\bfD)Y_s$ (with norm
$\|\cdot\|_{Y_s}$), noting the relationship
\[\nn \eta\nn^2=\varepsilon\|\zeta\|_{Y_1}^2,
\qquad  \eta(x,y)=\varepsilon^2 \zeta(\varepsilon x, \varepsilon^2 y)
\]
for $\zeta \in Y_1^\varepsilon$. Observe that $Y_s^\varepsilon$ coincides with
$\chi_\varepsilon(\bfD)H^s({\mathbb R}^2)$.

\section{Analyticity} \label{anal}

\subsection{The boundary-value problems} \label{Anal BVPs}

In this section we solve the boundary-value problems \eqref{A BVP 1}--\eqref{A BVP 5} and \eqref{B BVP 1}--\eqref{B BVP 5} and use these results to deduce that $H(\eta)$ and $\bfM(\eta)$ depend analytically upon $\eta \in \ZZ$,
where
$$\ZZ=\{\eta \in {\mathscr S}^\prime({\mathbb R}^2)\colon  \|\hat{\eta}_1\|_{L^1({\mathbb R}^2)} + \|\eta_2\|_3 < \infty\}$$
and
$$\eta_1 =\chi(\bfD)\eta, \qquad
\eta_2 = \left(1-\chi(\bfD)\right)\eta$$
(see Theorem \ref{main anal thm}(i) below for a precise statement).
We proceed by transforming \eqref{A BVP 1}--\eqref{A BVP 5} and \eqref{B BVP 1}--\eqref{B BVP 5} into equivalent boundary-value problems in the fixed domain $D_0$ by means
of the `flattening' transformation\linebreak
$\Sigma\colon D_0 \to D_\eta$ given by
$$
\Sigma\colon(x,y,v) \mapsto (x,y,v+\sigma(x,y,v)), \qquad \sigma(x,y,v)\coloneqq  \eta(x,y)(1+v)
$$
and differential operators
$$
\partial_x^\sigma \coloneqq  \partial_x - \frac{ \partial_x\sigma }{1+\partial_v\sigma} \partial_v, \quad
\partial_y^\sigma \coloneqq  \partial_y - \frac{ \partial_y\sigma }{1+\partial_v\sigma} \partial_v, \quad
\partial_v^\sigma \coloneqq  \frac{ \partial_v }{1+\partial_v\sigma}
$$
(with corresponding notation for differential operators derived from them).
Equations
\eqref{A BVP 1}--\eqref{A BVP 5} are equivalent to the flattened boundary-value problem\pagebreak
\begin{alignat}{2}
\curl^\sigma \curl^\sigma \bfA^\prime &= \alpha  \curl^\sigma \bfA^\prime  & &\mbox{in $D_0$}, \label{pre flattened A BVP 1} \\
\Div^\sigma \bfA^\prime &= 0 & & \mbox{in $D_0$}, \label{pre flattened A BVP 2} \\
\bfA^\prime \wedge \bfe_3 &= \bfzero & & \mbox{at $v=-1$}, \label{pre flattened A BVP 3} \\
\bfA^\prime \cdot \bfn &= 0 & & \mbox{at $v=0$}, \label{pre flattened A BVP 4} \\
(\curl^\sigma \bfA^\prime)_\parallel &= \nabla\Phi - \alpha \nablap\Delta^{-1}(\nablac \bfA_\parallel^{\prime\perp}) \qquad \label{pre flattened A BVP 5}
\end{alignat}
for $\bfA^\prime = \bfA \circ \Sigma$,
in terms of which
\begin{equation}
H(\eta)\Phi=\nablac \bfA_\parallel^{\prime\perp}, \label{Prime defn of H}
\end{equation}
while equations \eqref{B BVP 1}--\eqref{B BVP 5} are equivalent to the flattened boundary-value problem
\begin{alignat}{2}
\curl^\sigma \curl^\sigma \bfB^\prime &= \alpha  \curl^\sigma \bfB^\prime \qquad\quad& & \mbox{in $D_0$}, \label{pre flattened B BVP 1} \\
\Div^\sigma \bfB^\prime &= 0 & & \mbox{in $D_0$}, \label{pre flattened B BVP 2} \\
\bfB^\prime \wedge \bfe_3 &= \bfzero & & \mbox{at $v=-1$}, \label{pre flattened B BVP 3} \\
\bfB^\prime \cdot \bfn &= 0 & & \mbox{at $v=0$}, \label{pre flattened B BVP 4} \\
\nablac \bfB_{\parallel}^{\prime\perp} &= \nablac \bfg^\perp, \label{pre flattened B BVP 5}
\end{alignat}
for $\bfB^\prime = \bfB \circ \Sigma$,
in terms of which
$$
\bfM(\eta)\bfg=-(\curl^\sigma \bfB^\prime)_\parallel;
$$
note that the
orthogonal gradient part of $(\curl^\sigma \bfB^\prime)_\parallel$ is equal to
$-\alpha \, \nablap \Delta^{-1} (\nablac \bfB_\parallel^{\prime\perp})$ for any solution\linebreak
$\bfB^\prime \in H^2(D_0)^3$ of \eqref{Flattened B BVP 1}--\eqref{Flattened B BVP 5}.

It is in fact convenient to replace \eqref{A BVP 1}--\eqref{A BVP 5} with an equivalent boundary-value problem.
The following proposition was proved by Groves \& Horn \cite[Proposition 4.8]{GrovesHorn20} (under slightly different regularity assumptions on
$\eta$, the change in which does not affect the proof).

\begin{proposition}
Suppose that $\Phi \in \dot{H}^\frac{3}{2}({\mathbb R}^2)$ and $\eta$ lies in a sufficiently small neighbourhood of the origin in $\ZZ$. A function $\bfA \in H^2(D_\eta)^3$ solves
\eqref{A BVP 1}--\eqref{A BVP 5} if and only if it satisfies the boundary-value problem
\begin{alignat*}{2}
-\Delta \bfA &= \alpha  \curl \bfA & & \mbox{in $D_\eta$,}  \\
\bfA \wedge \bfe_3 &= \bfzero & & \mbox{at $z=-1$,}  \\
\partial_zA_3 & = 0 & & \mbox{at $z=-1$,} \\
\bfA \cdot \bfn &= 0 & & \mbox{at $z=\eta$,} \\
(\curl \bfA)_\parallel &= \nabla\Phi - \alpha \nablap\Delta^{-1}(\nablac \bfA_\parallel^\perp).\qquad
\end{alignat*}
\end{proposition}
\begin{corollary}
Suppose that $\Phi \in \dot{H}^\frac{3}{2}({\mathbb R}^2)$ and $\eta$ lies in a sufficiently small neighbourhood of the origin in $\ZZ$. A function $\bfA^\prime \in H^2(D_0)^3$ solves
\eqref{pre flattened A BVP 1}--\eqref{pre flattened A BVP 5} if and only if it satisfies the boundary-value problem
\begin{alignat}{2}
-\Delta^\sigma \bfA^\prime &= \alpha  \curl^\sigma \bfA^\prime & &\mbox{in $D_0$}, \label{pre alt flattened A BVP 1} \\
\bfA^\prime \wedge \bfe_3 &= \bfzero & & \mbox{at $v=-1$}, \label{pre alt flattened A BVP 2} \\
\partial_v A_3^\prime & = 0 & & \mbox{at $v=-1$,} \label{pre alt flattenedA BVP 3} \\
\bfA^\prime \cdot \bfn &= 0 & & \mbox{at $v=0$}, \label{pre alt flattened A BVP 4} \\
(\curl^\sigma \bfA^\prime)_\parallel &= \nabla\Phi - \alpha \nablap\Delta^{-1}(\nablac \bfA_\parallel^{\prime\perp}). \qquad \label{pre alt flattened A BVP 5}
\end{alignat}
\end{corollary}

We proceed by rewriting  \eqref{pre alt flattened A BVP 1}--\eqref{pre alt flattened A BVP 5} as
\begin{align}
& \parbox{8.25cm}{$-\Delta {\bfA^\prime} - \alpha\, \curl{\bfA^\prime}=\bfG^\sigma(\bfA^\prime)$}\mbox{in $D_0$,} \label{alt flat 1} \\
& \parbox{8.25cm}{${\bfA^\prime} \wedge \bfe_3 = \mathbf{0}$}\mbox{at $v=-1$,} \label{alt flat 2} \\
& \parbox{8.25cm}{$\partial_v A_3^\prime = 0$}\mbox{at $v=-1$,} \label{alt flat 3} \\
& \parbox{8.25cm}{$\bfA^\prime\cdot\bfe_3 = g^\sigma(\bfA^\prime)$}\mbox{at $v=0$,} \label{alt flat 4} \\
& \parbox{8.25cm}{$(\udl{\curl {\bfA^\prime}})_\mathrm{h} + \alpha \nabla^\perp \Delta^{-1} (\nabla\cdot\udl{\bfA^\prime}^{\!\perp}_\mathrm{h})= \bfg^\sigma(\bfA^\prime)+\nabla \Phi
$,} \label{alt flat 5}
\end{align}
where\pagebreak
\begin{align*}
\bfG^\sigma(\bfA^\prime) &= \Delta^\sigma \bfA^\prime + \alpha \curl^\sigma \bfA^\prime - \Delta \bfA^\prime - \alpha \curl \bfA^\prime, \\
\bfg^\sigma(\bfA^\prime) &= -(\udl{\curl^{\sigma}{\bfA^\prime}})_\mathrm{h}+(\udl{\curl \bfA^\prime})_\mathrm{h} 
-\nabla \eta (\udl{\curl^{\sigma}\bfA^\prime})_3
-\alpha\nabla^\perp\!\Delta^{-1}(\nabla \cdot(\nabla\eta^\perp \udl{A}_3^\prime)), \\
g^\sigma(\bfA^\prime) &=\nabla \eta \cdot \udl{\bfA^\prime}_\mathrm{h}
\end{align*}
and studying \eqref{alt flat 1}--\eqref{alt flat 5} in more generality.
(With a slight abuse of notation the underscore now denotes evaluation at $v=0$).
The following proposition treats an inhomogeneous linear version of \eqref{alt flat 1}--\eqref{alt flat 5}; it was proved by
Groves \& Horn \cite[Proposition 4.9]{GrovesHorn20}, who in particular give an explicit formula for the solution.

\begin{lemma} \label{Solve inhomo lin BVP}
Suppose that $s \geq 2$ and $\alpha^\star < \tfrac{1}{2}\pi$.
The boundary-value problem
\begin{align*}
& \parbox{8.25cm}{$-\Delta \bfC - \alpha\, \curl{\bfC}=\bfH$}\mbox{in $D_0$,} \\
& \parbox{8.25cm}{${\bfC} \wedge \bfe_3 = \mathbf{0}$}\mbox{at $v=-1$,} \\
& \parbox{8.25cm}{$\partial_v C_3 = 0$}\mbox{at $v=-1$,} \\
& \parbox{8.25cm}{$\bfC\cdot\bfe_3= h$}\mbox{at $v=0$,} \\
& \parbox{8.25cm}{$(\udl{\curl {\bfC}})_\mathrm{h} + \alpha \nabla^\perp \Delta^{-1} (\nabla\cdot\udl{\bfC}^{\!\perp}_\mathrm{h})= \bfh$}
\end{align*}
has a unique solution $\bfC \in H^s(D_0)^3$ for each $h \in H^{s-\frac{1}{2}}({\mathbb R}^2)$,
$\bfH \in H^{s-2}(D_0)^3$, $\bfh \in  H^{s-\frac{3}{2}}({\mathbb R}^2)^2$ and\linebreak
$|\alpha| \in [0,\alpha^\star]$. The solution operator defines a
mapping $H^{s-\frac{1}{2}}({\mathbb R}^2) \times H^{s-2}(D_0)^3 \times H^{s-\frac{3}{2}}({\mathbb R}^2)^2
\rightarrow H^s(D_0)^3$ which is bounded uniformly over $|\alpha| \in [0,\alpha^\star]$.
\end{lemma}

Lemma \ref{Solve inhomo lin BVP} can be used in particular to study the boundary-value problems
\begin{alignat}{4}
\curl \curl \bfA^{\prime 0} &= \alpha  \curl \bfA^{\prime 0} &  \curl \curl \bfB^{\prime 0} &= \alpha  \curl \bfB^{\prime 0} \qquad& &\mbox{in $D_0$}, \label{AB1} \\
\Div \bfA^{\prime 0} &= 0 &  \Div \bfB^{\prime 0} &= 0 & & \mbox{in $D_0$}, \label{AB2} \\
\bfA^{\prime 0} \wedge \bfe_3 &= \bfzero & \bfB^{\prime 0} \wedge \bfe_3 &= \bfzero & & \mbox{at $v=-1$}, \label{AB3} \\
\bfA^{\prime 0} \cdot \bfe_3 &= 0  &\bfB^{\prime 0} \cdot \bfe_3 &= 0 &  & \mbox{at $v=0$}, \label{AB4} \\
(\udl{\curl \bfA}^{\prime 0})_\mathrm{h} &= \nabla\Phi - \alpha \nablap\Delta^{-1}(\nablac (\udl{\bfA}^{\prime 0})_\mathrm{h}^\perp)\qquad\qquad
& \nablac (\udl{\bfB}^{\prime 0})_\mathrm{h}^\perp &= \nablac \bfg^\perp\hspace{5mm} \label{AB5}
\end{alignat}
for $\Phi \in \dot{H}^{s-\frac{1}{2}}({\mathbb R}^2)$ and $\bfg \in H^{s-\frac{1}{2}}({\mathbb R}^2)^2$ with $s \geq 2$ and $|\alpha| \in [0,\tfrac{1}{2}\pi)$.
The boundary-value problem for $\bfA^{\prime 0}$ has a unique solution $\bfA^{\prime 0}(\Phi) \in H^s(D_0)^3$, and
it follows from
$$H(0)\Phi = \nablac \udl{{\bfA}}^{\prime 0}(\Phi)_\mathrm{h}^\perp$$
and the explicit formula for $\bfA^{\prime 0}(\Phi)$ given by Groves \& Horn that
$$H(0)\Phi = D^2 \, \mathtt{t}(D^2),
\qquad
\mathtt{t}(\mu) =  
\begin{cases}
\dfrac{ \tan(\sqrt{\alpha^2-\mu}) }{\sqrt{\alpha^2-\mu}},   & \mbox{if $\mu < \alpha^2$,} \\[5mm]
\dfrac{ \tanh(\sqrt{\mu-\alpha^2}) }{\sqrt{\mu-\alpha^2}},  & \mbox{if $\mu \geq \alpha^2$}
\end{cases}
$$
and
$$\bfD=(D_1,D_2)^T=-\i\nabla, \qquad D=|\bfD|.$$
Note that
$H(0) \in L(\dot{H}^{s-\frac{1}{2}}({\mathbb R}^2),\check{H}^{s-\frac{3}{2}}({\mathbb R}^2))$ is an isomorphism because
$$
\| H(0)^{-1}\Psi \|_{\dot{H}^{s-\frac{1}{2}}}\!=\!\left\| \langle \bfk \rangle^{s-\frac{3}{2}} |\bfk| \frac{1}{|\bfk|^2 \mathtt{t}(|\bfk|^2)}
\hat{\Psi}\right\|_0\!\!=\!\left\| \langle \bfk\rangle^{s-\frac{1}{2}}|\bfk|^{-1}\frac{\langle \bfk \rangle^{-1}}{\mathtt{t}(|\bfk|^2)}\hat{\Psi}\right\|_0\!\!
\lesssim \|\langle \bfk\rangle^{s-\frac{1}{2}}|\bfk|^{-1}\hat{\Psi}\|_0
\!=\! \|\Psi\|_{\check{H}^{s-\frac{3}{2}}},
$$
where we have used the fact that $\langle \bfk \rangle^{-1}/\mathtt{t}(|\bfk|^2)$ is bounded.

Observe that
$\bfB^{\prime 0}(\bfg)\coloneqq \bfA^{\prime 0}(\Phi)$ with $\Phi=H(0)^{-1}(\nablac\bfg^\perp)$ solves the boundary-value problem for $\bfB^{\prime 0}$ because
$$\nablac \bfg^\perp = H(0)\Phi = \nablac \udl{\bfA}^{\prime 0}(\Phi)_\mathrm{h}^\perp= \nablac \udl{\bfB}^{\prime 0}(\bfg)_\mathrm{h}^\perp;$$
this solution is unique because any other solution $\bfB^{\prime 0}(\bfg)$ is equal to $\bfA^{\prime 0}(\Phi)$ with $\Phi=\Delta^{-1} (\nablac \udl{\curl \bfB}^{\prime 0}(\bfg)_\mathrm{h})$,
so that
$$H(0)\Phi = \nablac \udl{\bfA}^{\prime 0}(\Phi)_\mathrm{h}^\perp =  \nablac \udl{\bfB}^{\prime 0}(\bfg)_\mathrm{h}^\perp=\nablac \bfg^\perp.$$
It now follows from
$$\bfM(0)\bfg=-\udl{\curl \bfB}^{\prime 0}(\bfg)_\mathrm{h}=-\udl{\curl \bfA}^{\prime 0}(\Phi)_\mathrm{h}=-\nabla\Phi+\alpha\nablap \Delta^{-1}(\nablac \bfg^\perp)$$
that $\bfM(0) \in L(H^{s-\frac{1}{2}}({\mathbb R}^2)^2,H^{s-\frac{3}{2}}({\mathbb R}^2)^2)$ is given by
$$\bfM(0)\bfg =\frac{1}{D^2} \, \left( \alpha \, \bfD^{\perp} + \bfD \, \mathtt{c}(D^2) \right) \, \bfD \cdot \bfg^{\perp},
\qquad
\mathtt{c}(\mu) =
\begin{cases}
\sqrt{\alpha^2-\mu} \, \cot(\sqrt{\alpha^2-\mu}),   & \mbox{if $\mu < \alpha^2$,} \\[2mm]
\sqrt{\mu-\alpha^2} \, \coth(\sqrt{\mu-\alpha^2}), & \mbox{if $\mu \geq \alpha^2$.}
\end{cases}
$$

The following lemma is the key to solving the nonlinear problem \eqref{alt flat 1}--\eqref{alt flat 5}. We henceforth assume that\linebreak
$|\alpha| \in [0,\tfrac{1}{2}\pi)$ and replace $V$ with a smaller neighbourhood of the origin in $\ZZ$ if necessary.

\begin{lemma} \label{Solve more general BVPs}
Suppose that $s \geq 2$ and that $\bfH\colon \ZZ \times H^s(D_0)^3 \rightarrow H^{s-2}(D_0)^3$, $h\colon \ZZ \times H^s(D_0)^3 \rightarrow H^{s-\frac{1}{2}}({\mathbb R}^2)$
and $\bfh\colon\ZZ \times H^s(D_0)^3 \rightarrow H^{s-\frac{3}{2}}({\mathbb R}^2)^2$ are analytic at the origin, linear in their second argument and satisfy
$\bfH(0,\bfC)=\bfzero$, $h(0,\bfC)=0$, $\bfh(0,\bfC)=\bfzero$ for all $\bfC \in H^s(D_0)^3$.
There exists a neighbourhood $V$ of the origin in $\ZZ$ with the property that the boundary-value problem
\begin{align*}
& \parbox{8.25cm}{$-\Delta {\bfC} - \alpha\, \curl{\bfC}=\bfH(\eta,\bfC)$}\mbox{in $D_0$,} \\
& \parbox{8.25cm}{${\bfC} \wedge \bfe_3 = \mathbf{0}$}\mbox{at $v=-1$,} \\
& \parbox{8.25cm}{$\partial_v C_3 = 0$}\mbox{at $v=-1$,} \\
& \parbox{8.25cm}{$\bfC\cdot\bfe_3 = h(\eta,\bfC)$}\mbox{at $v=0$,} \\
& \parbox{8.25cm}{$(\udl{\curl {\bfC}})_\mathrm{h} + \alpha \nabla^\perp \Delta^{-1} (\nabla\cdot\udl{\bfC}^{\!\perp}_\mathrm{h})= \bfh(\eta,\bfC)+\nabla \Phi
$,}
\end{align*}
has a unique solution $\bfC=\bfC(\eta,\Phi)$ in $H^s(D_0)^3$ which depends analytically upon
$\eta \in V$ and
$\Phi \in \dot{H}^{s+\frac{1}{2}}({\mathbb R}^2)$ (and linearly upon $\Phi$).
\end{lemma}
{\bf Proof.} The formula
$${\mathcal H}(\bfC,\eta,\Phi)=\begin{pmatrix} -\Delta {\bfC} - \alpha\, \curl{\bfC}-\bfH(\eta,\bfC) \\
\udl{\bfC}\cdot\bfe_3- h(\eta,\bfC) \\
(\udl{\curl {\bfC}})_\mathrm{h} + \alpha \nabla^\perp \Delta^{-1} (\nabla\cdot\udl{\bfC}^{\!\perp}_\mathrm{h}) - \bfh(\eta,\bfC)-\nabla \Phi 
\end{pmatrix},$$
defines a mapping
$${\mathcal H}\colon S \times \ZZ \times \dot{H}^{s+\frac{1}{2}}({\mathbb R}) \rightarrow H^{s-2}(D_0)^3 \times H^{s-\frac{1}{2}}({\mathbb R}^2) \times
H^{s-\frac{3}{2}}({\mathbb R}^2)^2,$$
where $S=\{\bfC \in H^s(D_0)^3\colon \bfC \wedge \bfe_3\big|_{v=-1}={\mathbf 0}, \partial_v C_3\big|_{v=-1}=0\}$,
which is analytic at the origin.
Furthermore, ${\mathcal H}({\mathbf 0},0,0)=({\mathbf 0},0,{\mathbf 0})$, and the calculation
$$\mathrm{d}_1{\mathcal H}[{\mathbf 0},0,0](\bfC)=\begin{pmatrix}
-\Delta {\bfC} - \alpha\, \curl{\bfC} \\
\udl{\bfC}\cdot\bfe_3 \\
(\udl{\curl {\bfC}})_\mathrm{h} + \alpha \nabla^\perp \Delta^{-1} (\nabla\cdot\udl{\bfC}^{\!\perp}_\mathrm{h})
\end{pmatrix}$$
and Lemma \ref{Solve inhomo lin BVP} show that
$$\mathrm{d}_1{\mathcal H}[{\mathbf 0},0,\bfzero]\colon S \rightarrow H^{s-2}(D_0)^3 \times H^{s-\frac{1}{2}}({\mathbb R}^2) \times
H^{s-\frac{3}{2}}({\mathbb R}^2)^2$$
is an isomorphism.
The analytic implicit-function theorem (Buffoni \& Toland \cite[Theorem 4.5.3]{BuffoniToland})
asserts the existence of open neighbourhoods $V_1$ and $V_2$ of the origin in respectively
$\ZZ \times \dot{H}^{s+\frac{1}{2}}({\mathbb R})$ and $S$ such that
the equation
$${\mathcal H}(\bfC,\eta,\Phi)=({\mathbf 0},0,{\mathbf 0})$$
and hence the boundary-value problem \eqref{alt flat 1}--\eqref{alt flat 5}
has a unique solution $\bfC=\bfC(\eta,\Phi)$ in $V_2$ for each $(\eta,\Phi) \in V_1$; furthermore
$\bfC(\eta,\Phi)$ depends
analytically upon $\eta$ and $\Phi$. Since $\bfC(\eta,\Phi)$ depends linearly upon $\Phi$ one can without loss of generality take
$V_1=V \times \dot{H}^{s+\frac{1}{2}}({\mathbb R})$, and clearly $V_2=S$
(with $\Phi=0$ the construction yields a unique solution in a neighbourhood of the origin in $S$,
which is evidently the zero solution).\qed

To apply Lemma \ref{Solve more general BVPs} to the boundary-value problem \eqref{alt flat 1}--\eqref{alt flat 5} we examine
the mapping properties of $\bfG^\sigma$, $\bfg^\sigma$ and $g^\sigma$.

\begin{proposition} \label{The gs}
$ $
\begin{itemize}
\item[(i)]
The formula $(\eta,\bfC) \mapsto \bfG^\sigma(\bfC)$ defines an analytic mapping $V \times H^s(D_0)^3 \rightarrow H^{s-2}(D_0)^3$ for $s \in [2,3]$.
\item[(ii)]
The formula $(\eta,\bfC) \mapsto \bfg^\sigma(\bfC)$ defines an analytic mapping $V \times H^s(D_0)^3 \rightarrow H^{s-\frac{1}{2}}({\mathbb R}^2)^2$
for $s\in [2,3]$
\item[(iii)]
The formula
$(\eta,\bfC) \mapsto g^\sigma(\bfC)$ defines an analytic mapping $V \times H^s(D_0)^3 \rightarrow H^{s-\frac{1}{2}}({\mathbb R}^2)$ for $s \in [2,\frac{5}{2}]$ and
$V \times H^s(D_0)^3 \rightarrow H^2({\mathbb R}^2)$ for $s \in [\frac{5}{2},3]$.
\end{itemize}
\end{proposition}
{\bf Proof.} The analyticity of $(\eta,\bfC) \mapsto \bfG^\sigma(\bfC)$ at the origin as a
mapping $\ZZ \times H^s(D_0)^3 \rightarrow H^{s-2}(D_0)^3$ follows from the
explicit expression
\begin{align*}
\bfG^\sigma(\bfC) &=  - 2 \frac{1+v}{1+\eta}( \eta_x  \partial_{vx}^2\bfC + \eta_y  \partial_{vy}^2\bfC ) - \frac{1+v}{1+\eta}\Delta\eta\, \partial_v\bfC  \\
&\qquad\mbox{}+ 2 \frac{1+v}{(1+\eta)^2} |\nabla\eta|^2 \partial_v\bfC + \left( \frac{1+v}{1+\eta} \right)^{\!\!2}  |\nabla\eta|^2\partial_v^2\bfC - \frac{\eta^2+2\eta}{(1+\eta)^2} \partial_v^2\bfC\\
&\qquad\mbox{}-\alpha\frac{\eta}{1+\eta} (-\partial_v C_2,\partial_v C_1,0)^T
-\alpha \frac{1+v}{1+\eta} (\eta_y  \partial_vC_3, -\eta_x  \partial_vC_3,\eta_x  \partial_v C_2-\eta_y  \partial_v C_1)^T
\end{align*}
by writing
$$\frac{1}{1+\eta}=1-\frac{\eta}{1+\eta}, \qquad \frac{1}{(1+\eta)^2}=1-\frac{\eta^2+2\eta}{(1+\eta)^2}$$
and noting that
\begin{itemize}
\item
the bilinear mappings $(\eta,\tilde{f}) \mapsto \eta_x \tilde{f}$ and $(\eta,\tilde{f}) \mapsto \eta_y \tilde{f}$
are bounded
$\ZZ \times H^{s-2}(D_0)\rightarrow H^{s-2}(D_0)$
because
\begin{align*}
\|\eta_x\tilde{f}\|_{L^2(D_0)}&\lesssim (\|\eta_{1x}\|_{\infty}+\|\eta_{2x}\|_{\infty})\|\tilde{f}\|_{L^2(D_0)}\lesssim (\|\hat{\eta}_1\|_{L^1({\mathbb R}^2)} +\|\eta_2\|_3)\|\tilde{f}\|_{L^2(D_0)},\\
\|\eta_x\tilde{f}\|_{H^1(D_0)}&\lesssim (\|\eta_{1x}\|_{1,\infty}+\|\eta_{2x}\|_{1,\infty})\|\tilde{f}\|_{H^1(D_0)}\lesssim (\|\hat{\eta}_1\|_{L^1({\mathbb R}^2)} +\|\eta_2\|_3)\|\tilde{f}\|_{H^1(D_0)},\\[6pt]
\|\eta_y\tilde{f}\|_{L^2(D_0)}&\lesssim (\|\eta_{1y}\|_{\infty}+\|\eta_{2y}\|_{\infty})\|\tilde{f}\|_{L^2(D_0)}\lesssim (\|\hat{\eta}_1\|_{L^1({\mathbb R}^2)} +\|\eta_2\|_3)\|\tilde{f}\|_{L^2(D_0)},\\
\|\eta_y\tilde{f}\|_{H^1(D_0)}&\lesssim (\|\eta_{1y}\|_{1,\infty}+\|\eta_{2y}\|_{1,\infty})\|\tilde{f}\|_{H^1(D_0)}\lesssim (\|\hat{\eta}_1\|_{L^1({\mathbb R}^2)} +\|\eta_2\|_3)\|\tilde{f}\|_{H^1(D_0)}
\end{align*}
(interpolate between $L(L^2(D_0),L^2(D_0))$ and $L(H^1(D_0),H^1(D_0))$),
\item
the bilinear mapping $(\eta,\tilde{f}) \mapsto \Delta\eta\,\partial_v\tilde{f}$ is bounded
$\ZZ \times H^{s}(D_0)\rightarrow H^{s-2}(D_0)$ because
\begin{align*}
\|\Delta\eta\,\partial_v\tilde{f}\|_{L^2(D_0)} & \lesssim \|\Delta\eta_1\|_{\infty}\|\partial_v\tilde{f}\|_{L^2(D_0)}  + \|\Delta \eta_2\|_{L^4({\mathbb R}^2)}\|\partial_v\tilde{f}\|_{L^4(D_0)}\\
& \lesssim \|\eta_1\|_{2,\infty} \|\partial_v\tilde{f}\|_{H^1(D_0)}+\|\Delta\eta_2\|_1 \|\partial_v\tilde{f}\|_{H^1(D_0)}\\
& \lesssim (\|\hat{\eta}_1\|_{L^1({\mathbb R}^2)} +\|\eta_2\|_3)\|\tilde{f}\|_{H^2(D_0)}, \\[6pt]
\|\Delta\eta\,\partial_v\tilde{f}\|_{H^1(D_0)}
&\lesssim \|\Delta\eta_1\|_{1,\infty}\|\partial_v\tilde{f}\|_{H^1(D_0)} + \|\nabla (\Delta \eta_2) \partial_v\tilde{f}\|_{L^2(D_0)}+\|\Delta \eta_2 \nabla(\partial_v\tilde{f})\|_{L^2(D_0)}\| \\
& \lesssim \|\eta_1\|_{3,\infty} \|\tilde{f}\|_{H^2(D_0)}+\|\eta_2\|_3 \|\partial_v\tilde{f}\|_{L^\infty(D_0)}+\|\Delta \eta_2\|_{L^4({\mathbb R}^2)}\|\partial_v \tilde{f}\|_{W^{1,4}(D_0)}\\
& \lesssim (\|\hat{\eta}_1\|_{L^1({\mathbb R}^2)} +\|\eta_2\|_3)\|\tilde{f}\|_{H^3(D_0)},
\end{align*}
(interpolate between $L(H^2(D_0),L^2(D_0))$ and $L(H^3(D_0),H^1(D_0))$),
\item
the mapping $\tilde{f} \mapsto (1+v)\tilde{f}$ belongs to $L(H^{s-2}(D_0),H^{s-2}(D_0))$ for all $t \geq 0$,
\item
a function $f\colon{\mathbb R} \rightarrow {\mathbb R}$ which is analytic at the origin (in particular
$f(t)=t(1+t)^{-1}$ and\linebreak
$f(t)=(t^2+2t)(1+t)^{-1}$) induces a mapping $W^{1,\infty}({\mathbb R}^2)
\rightarrow W^{1,\infty}({\mathbb R}^2)$ and hence $\ZZ \mapsto W^{1,\infty}({\mathbb R}^2)$
which is analytic at the origin,
\item
the bilinear mapping $(\rho,\tilde{f}) \mapsto \rho\tilde{f}$ is bounded $W^{1,\infty}({\mathbb R}^2) \times L^2(D_0) \rightarrow L^2(D_0)$ and\linebreak
$W^{1,\infty}({\mathbb R}^2) \times H^1(D_0) \rightarrow H^1(D_0)$, and hence
$W^{1,\infty}({\mathbb R}^2) \times H^{s-2}(D_0) \rightarrow H^{s-2}(D_0)$.
\end{itemize}
The results for $(\eta,\bfC) \mapsto \bfg^\sigma(\bfC)$ and $(\eta,\bfC) \mapsto g^\sigma(\bfC)$ are obtained in a similar fashion.\qed
\begin{corollary} \label{Aprimeisanal}
Suppose that $s \in [2,\frac{5}{2}]$.
The boundary-value problem \eqref{alt flat 1}--\eqref{alt flat 5} has a unique solution\linebreak
 $\bfA^\prime=\bfA^\prime(\eta,\Phi) \in H^s(D_0)^3$ which depends analytically upon
$\eta \in V$ and
$\Phi \in \dot{H}^{s-\frac{1}{2}}({\mathbb R}^2)$ $\mbox{(and linearly upon $\Phi$).}$
\end{corollary}
{\bf Proof.} This result is obtained from Lemma \ref{Solve more general BVPs} with
$(\bfH(\eta,\bfC),\bfh(\eta,\bfC),h(\eta,\bfC)=(\bfG^\sigma(\bfC),\bfg^\sigma(\bfC),g^\sigma(\bfC))$ and Proposition \ref{The gs}.\qed

Using definition \eqref{Prime defn of H}, we find from Corollary \ref{Aprimeisanal} with $s=\frac{5}{2}$ that
$H(\eta)\Phi \in \check{H}^1({\mathbb R}^2)$ for $\Phi \in \dot{H}^2({\mathbb R}^2)$.
The regularity of $\bfA^\prime$, and hence of $H(\eta)\Phi$, cannot be improved by increasing the regularity of $\Phi$, the limiting factor being 
Proposition \ref{The gs}(iii). It is however possible to replace \eqref{pre flattened A BVP 4}  and \eqref{pre alt flattened A BVP 4} in
the boundary-value problems defining $H(\eta)$ with a simpler boundary condition and thus obtain $H(\eta)\Phi \in \check{H}^\frac{3}{2}({\mathbb R}^2)$
for $\Phi \in \dot{H}^\frac{5}{2}({\mathbb R}^2)$.

\begin{lemma} \label{Change to tildeA} $ $
\begin{itemize}
\item[(i)]
The boundary-value problem
\begin{alignat*}{2}
\Delta^\sigma \phi &=0 & & \mbox{in $D_0$}, \\
\phi &=0 & & \mbox{at $v=-1$}, \\
\Grad^\sigma \phi\cdot\bfn & =\nabla \eta \cdot \tilde\bfA_\mathrm{h}\qquad  & & \mbox{at $v=0$}.
\end{alignat*}
has a unique solution $\phi=\phi(\eta,\tilde\bfA) \in H^\frac{7}{2}(D_0)$ which depends analytically upon $\eta \in V$ and $\tilde\bfA \in H^\frac{5}{2}(D_0)^3$
(and linearly upon $\tilde\bfA$).
\item[(ii)] The mapping $\tilde\bfA \mapsto \tilde\bfA + \Grad^\sigma\phi(\eta,\tilde\bfA)$ is an isomorphism $H^\frac{5}{2}(D_0)^3 \rightarrow
H^\frac{5}{2}(D_0)^3$. Both it and its inverse depend analytically upon $\eta \in V$.
\item[(iii)]
Suppose that $\Phi \in \dot{H}^2({\mathbb R}^2)$ and $\eta \in V$. The function $\tilde\bfA \in H^\frac{5}{2}(D_0)^3$ solves
\begin{alignat}{2}
\curl^\sigma \curl^\sigma \tilde\bfA &= \alpha  \curl^\sigma \tilde\bfA  & &\mbox{in $D_0$}, \label{Flattened A BVP 1} \\
\Div^\sigma \tilde\bfA &= 0 & & \mbox{in $D_0$}, \label{Flattened A BVP 2} \\
\tilde\bfA \wedge \bfe_3 &= \bfzero & & \mbox{at $v=-1$}, \label{Flattened A BVP 3} \\
\tilde\bfA \cdot \bfe_3 &= 0 & & \mbox{at $v=0$}, \label{Flattened A BVP 4} \\
(\curl^\sigma \tilde\bfA)_\parallel &= \nabla\Phi - \alpha \nablap\Delta^{-1}(\nablac \tilde\bfA_\parallel^\perp), \qquad \label{Flattened A BVP 5}
\end{alignat}
and equivalently
\begin{alignat}{2}
-\Delta^\sigma \tilde\bfA &= \alpha  \curl^\sigma \tilde\bfA & &\mbox{in $D_0$}, \label{alt flattened A BVP 1} \\
\tilde\bfA \wedge \bfe_3 &= \bfzero & & \mbox{at $v=-1$}, \label{alt flattened A BVP 2} \\
\partial_v \tilde{A}_3 & = 0 & & \mbox{at $v=-1$,} \label{alt flattenedA BVP 3} \\
\tilde\bfA \cdot \bfe_3 &= 0 & & \mbox{at $v=0$}, \label{alt flattened A BVP 4} \\
(\curl^\sigma \tilde\bfA)_\parallel &= \nabla\Phi - \alpha \nablap\Delta^{-1}(\nablac \tilde{\bfA}_\parallel^\perp). \qquad \label{alt flattened A BVP 5}
\end{alignat}
if and only if $\bfA^\prime= \tilde\bfA + \Grad^\sigma\phi(\eta,\tilde\bfA)\in H^\frac{5}{2}(D_0)^3$ solves \eqref{pre flattened A BVP 1}--\eqref{pre flattened A BVP 5}
and equivalently \eqref{pre alt flattened A BVP 1}--\eqref{pre alt flattened A BVP 5}.
\item[(iv)] The boundary-value problems \eqref{Flattened A BVP 1}--\eqref{Flattened A BVP 5} and
\eqref{alt flattened A BVP 1}--\eqref{alt flattened A BVP 5} have a unique solution $\tilde\bfA =\tilde\bfA(\eta,\Phi) \in H^3(D_0)^3$ which 
depend analytically upon $\eta \in V$ and $\Phi \in \dot{H}^\frac{5}{2}({\mathbb R}^2)$
(and linearly upon $\Phi$).
\end{itemize}
\end{lemma}
{\bf Proof.}
\begin{itemize}
\item[(i)]
See Groves \emph{et al.} \cite[Proposition 2.13]{GrovesNilssonPasqualiWahlen24}.
\item[(ii)]
This statement follows from the analytic inverse-function theorem.
\item[(iii)]
A straightforward calculation shows that $\tilde\bfA$ satisfies \eqref{Flattened A BVP 1}--\eqref{Flattened A BVP 5}
if and only if $\bfA^\prime$ satisfies \eqref{pre flattened A BVP 1}--\eqref{pre flattened A BVP 5} and that
$\tilde\bfA$ satisfies \eqref{alt flattened A BVP 1}--\eqref{alt flattened A BVP 5}
if and only if $\bfA^\prime$ satisfies \eqref{pre alt flattened A BVP 1}--\eqref{pre alt flattened A BVP 5}.
The equivalence of \eqref{Flattened A BVP 1}--\eqref{Flattened A BVP 5} and \eqref{alt flattened A BVP 1}--\eqref{alt flattened A BVP 5}
follows from the equivalence of
\eqref{pre flattened A BVP 1}--\eqref{pre flattened A BVP 5} and \eqref{pre alt flattened A BVP 1}--\eqref{pre alt flattened A BVP 5}.
\item[(iv)]
This result is proved by applying Lemma \ref{Solve more general BVPs} with $(\bfH(\eta,\bfC),\bfh(\eta,\bfC),h(\eta,\bfC)=(\bfG^\sigma(\bfC),\bfg^\sigma(\bfC),0)$.\qed
\end{itemize}\pagebreak

Observe that $\nablac \bfA_\parallel^{\prime\perp}=\nablac \tilde\bfA_\parallel^\perp$, such that
\begin{equation}
H(\eta)\Phi=\nablac \tilde\bfA_\parallel^\perp \label{Defn of H tilde}
\end{equation}
for $\Phi \in \dot{H}^\frac{5}{2}({\mathbb R}^2)$, where $\tilde\bfA \in H^3(D_0)^3$ solves \eqref{Flattened A BVP 1}--\eqref{Flattened A BVP 5} or
\eqref{alt flattened A BVP 1}--\eqref{alt flattened A BVP 5}; it follows that $H(\eta)\Phi \in \check{H}^\frac{3}{2}({\mathbb R}^2)$.

We can apply the same procedure to the boundary-value problem for $\bfM(\eta)$.
Suppose that $\bfg \in H^2({\mathbb R}^2)^2$ and $\eta \in V$. The 
function $\tilde\bfB \in H^\frac{5}{2}(D_0)^3$ solves \begin{alignat}{2}
\curl^\sigma \curl^\sigma \tilde\bfB &= \alpha  \curl^\sigma \tilde\bfB \qquad\quad& & \mbox{in $D_0$}, \label{Flattened B BVP 1} \\
\Div^\sigma \tilde\bfB &= 0 & & \mbox{in $D_0$}, \label{Flattened B BVP 2} \\
\tilde\bfB \wedge \bfe_3 &= \bfzero & & \mbox{at $v=-1$}, \label{Flattened B BVP 3} \\
\tilde\bfB \cdot \bfe_3 &= 0 & & \mbox{at $v=0$}, \label{Flattened B BVP 4} \\
\nablac \tilde{\bfB}_{\parallel}^{\perp} &= \nablac \bfg^\perp. \label{Flattened B BVP 5}
\end{alignat}
if and only if $\bfB^\prime=\tilde\bfB + \Grad^\sigma\phi(\eta,\tilde\bfB) \in H^\frac{5}{2}(D_0)^3$ solves
\eqref{pre flattened B BVP 1}--\eqref{pre flattened B BVP 5}. In Theorem \ref{main anal thm}(ii) below we deduce from Lemma
\ref{Change to tildeA}(iv) that the boundary-value problem \eqref{Flattened B BVP 1}--\eqref{Flattened B BVP 5}
has a unique solution $\tilde\bfB=\tilde\bfB(\eta,\Phi) \in H^3(D_0)^3$ which 
depends analytically upon $\eta \in V$ and $\bfg \in H^\frac{5}{2}({\mathbb R}^2)^2$
(and linearly upon $\bfg$). Since $(\curl^\sigma \bfB^\prime)_\parallel=(\curl^\sigma \tilde\bfB)_\parallel$ we conclude that
\begin{equation}
\bfM(\eta)\bfg=-(\curl^\sigma \tilde\bfB)_\parallel, \label{Defn of M tilde}
\end{equation}
such that $\bfM(\eta)\bfg \in H^\frac{3}{2}({\mathbb R})^2$ for $\bfg \in H^\frac{5}{2}({\mathbb R}^2)^2$.

\begin{theorem} \label{main anal thm} \hspace{1cm}
\begin{itemize}
\item[(i)] The mappings $\eta \mapsto H(\eta)$ and $\eta \mapsto \bfM(\eta)$ are analytic $V \to L(\dot{H}^{\frac{5}{2}}({\mathbb R}^2),\check{H}^{\frac{3}{2}}({\mathbb R}^2))$ and\linebreak
$V \to L(H^{\frac{5}{2}}(\mathbb{R}^2)^2, H^{\frac{3}{2}}(\mathbb{R}^2)^2)$ respectively.
\item[(ii)] The boundary-value problem \eqref{Flattened B BVP 1}--\eqref{Flattened B BVP 5} has a unique solution
$\tilde{\bfB}=\tilde{\bfB}(\eta,\bfg)$ in $H^3(D_0)^3$ which depends analytically upon
$\eta \in V$ and $\bfg \in H^{\frac{5}{2}}({\mathbb R}^2)^2$ (and linearly upon $\bfg$).
\end{itemize}
\end{theorem}
{\bf Proof.} The analyticity of $H(\cdot)\colon V \to L(\dot{H}^{\frac{5}{2}}({\mathbb R}^2),\check{H}^{\frac{3}{2}}({\mathbb R}^2))$ follows from Lemma \ref{Change to tildeA}(iii) and equation \eqref{Defn of H tilde}. Since 
$H(0) \in L(\dot{H}^{\frac{5}{2}}({\mathbb R}^2),\check{H}^{\frac{3}{2}}({\mathbb R}^2))$ is an isomorphism we conclude
that $H(\eta) \in L(\dot{H}^{\frac{5}{2}}({\mathbb R}^2),\check{H}^{\frac{3}{2}}({\mathbb R}^2))$ is an isomorphism for each $\eta \in V$ and that $H(\eta)^{-1} \in L(\check{H}^{\frac{3}{2}}({\mathbb R}^2),\dot{H}^{\frac{5}{2}}({\mathbb R}^2))$
also depends analytically upon $\eta \in V$.

The next step is to note that $\tilde{\bfB}(\eta,\bfg)= \tilde\bfA(\eta,\Phi)$ with $\Phi=H(\eta)^{-1}(\nablac \bfg^\perp)$ 
depends analytically upon $\eta$ and $\bfg$, and solves \eqref{Flattened B BVP 1}--\eqref{Flattened B BVP 5} 
since by construction
$$\nablac \bfg^\perp = H(\eta)\Phi = \nablac \tilde\bfA(\eta,\Phi)_\parallel^\perp= \nablac \tilde\bfB(\eta,\Phi)_\parallel^\perp.$$
The uniqueness of this solution follows by noting that any other solution
$\tilde{\bfB}(\eta,\bfg)$ is equal to $\tilde\bfA(\eta,\Phi)$ with\linebreak
$\Phi=\Delta^{-1}\nablac (\curl^\sigma \tilde{\bfB})_\parallel$, so that 
$$H(\eta)\Phi=\nablac \tilde{\bfA}(\eta,\Phi)_\parallel^\perp=\nablac \tilde{\bfB}(\eta,\bfg)_\parallel^\perp = \nablac \bfg^\perp,$$
 that is $\Phi=H(\eta)^{-1}(\nablac \bfg^\perp)$. Finally, the analyticity of $\bfM$ follows from the calculation%
\begin{align*}
\bfM(\eta)\bfg &=-(\curl^\sigma \tilde{\bfB}(\eta,\bfg))_\parallel \\
& = -(\curl^\sigma \tilde{\bfA}(\eta,\Phi))_\parallel  \\
& = - \nabla \Phi + \alpha \, \nablap \Delta^{-1} (\nablac \bfg^\perp)
\end{align*}
with $\Phi=H(\eta)^{-1}(\nablac \bfg^\perp)$.\qed

We now choose $M>0$ sufficiently small and note that $H^3({\mathbb R}^2)$ is continuously embedded in $\ZZ$ and
$$U=\{\eta \in H^3({\mathbb R}^2)\colon \|\eta\|_\ZZ<M\}$$
is an open neighbourhood of the origin in $H^3({\mathbb R}^2)$.

\begin{proposition}
The mappings $\eta \mapsto \bfM(\eta)$ and $\eta \mapsto \bfT(\eta)$
are analytic $U \mapsto L(H^\frac{5}{2}({\mathbb R}^2)^2,H^\frac{3}{2}({\mathbb R}^2)^2)$ and\linebreak
$U \to H^{\frac{3}{2}}(\mathbb{R}^2)^2$ respectively.
\end{proposition}
{\bf Proof.} This result follows from Theorem \ref{main anal thm}(i), the formula $\bfT(\eta)=\bfM(\eta)\bfS(\eta)$ and the fact that
$\eta \mapsto \bfS(\eta)$ is an analytic mapping $U \to H^3(\mathbb{R}^2)^2$.\qed
\begin{corollary}
The formula \eqref{Defn of J}
defines an analytic function $\JJ\colon U \rightarrow H^1({\mathbb R}^2)$.
\end{corollary}
{\bf Proof.} We proceed by writing the formula as
\begin{align*}
\JJ(\eta)&=\frac{1}{2} |\bfT(\eta)|^2 -\frac{1}{2}( -\nabla \cdot \bfS(\eta)^\perp + \bfT(\eta) \cdot \nabla\eta)^2
+ \frac{ |\nabla \eta|^2( -\nabla \cdot \bfS(\eta)^\perp + \bfT(\eta) \cdot \nabla\eta)^2 }{2(1+|\nabla\eta|^2)} \nonumber \\
&\qquad\qquad\mbox{}+ \bfc\cdot\bfT(\eta) +\alpha \bfT(\eta)\cdot\bfS(\eta)+ \eta -\beta\Delta\eta \\
&\qquad\qquad\mbox{}
+\beta \left(\frac{|\nabla\eta|^2\eta_x}{(1+|\nabla\eta|^2)^{\frac{1}{2}}(1+(1+|\nabla\eta|^2)^{\frac{1}{2}})}\right)_x+\beta \left(\frac{|\nabla\eta|^2\eta_y}{(1+|\nabla\eta|^2)^{\frac{1}{2}}(1+(1+|\nabla\eta|^2)^{\frac{1}{2}})}\right)_y,
\end{align*}
from which the result follows because $\eta \mapsto \bfS(\eta)$, $\eta \mapsto \bfT(\eta)$ and 
$$\eta \mapsto \frac{|\nabla\eta|^2}{1+|\nabla\eta|^2}, \qquad \eta\mapsto \frac{|\nabla\eta|^2}{(1+|\nabla\eta|^2)^{\frac{1}{2}}(1+(1+|\nabla\eta|^2)^{\frac{1}{2}})}$$
are analytic mappings $U \rightarrow H^3({\mathbb R}^2)^2$, $U \rightarrow H^\frac{3}{2}({\mathbb R}^2)^2$
and $U \rightarrow H^2({\mathbb R}^2)$ respectively and $H^\frac{3}{2}({\mathbb R}^2)$ is a Banach algebra.\qed

\subsection{Taylor expansions}

The terms in the Taylor expansions
\begin{equation}
\tilde\bfB(\eta)=\sum_{k=0}^\infty \tilde\bfB^k(\eta), \qquad \tilde{\bfB}^k=\frac{1}{k!}\mathrm{d}^k\tilde\bfB[0](\eta^{(k)}), \label{Taylor B}
\end{equation}
of $\tilde\bfB$ and
\begin{equation}
\bfM(\eta)=\sum_{k=0}^\infty \bfM_k(\eta), \qquad \bfM_k = \frac{1}{k!}\mathrm{d}^k\bfM[0](\eta^{(k)}), \label{Taylor M}
\end{equation}
of $\bfM\colon U \mapsto L(H^3({\mathbb R}^2)^2,H^\frac{3}{2}({\mathbb R}^2)^2)$ can be determined recursively by substituting them into
\eqref{Flattened B BVP 1}--\eqref{Flattened B BVP 5} and \eqref{Defn of M tilde}.
It has already been established that
\begin{equation}
\bfM_0\bfg =-(\udl{\curl \tilde\bfB}^0)_\mathrm{h}=\frac{1}{D^2} \bfL\, \bfD \cdot \bfg^{\perp}, \qquad \bfL= \alpha  \bfD^{\perp} +  \mathtt{c}(D^2)\bfD, \label{Formula for M0}
\end{equation}
where $\tilde\bfB^0$ is the unique solution of \eqref{AB1}--\eqref{AB5}. Observing that $\bfM_0$ also defines an operator in\linebreak
$L(H^{s-\frac{1}{2}}({\mathbb R}^2)^2,H^{s-\frac{3}{2}}({\mathbb R}^2)^2)$ (with $\tilde\bfB^0 \in H^s(D_0)^3$) for $s \geq 2$, we can also obtain an explicit
expression for $\bfM_1(\eta)\bfg$.

\begin{lemma}
The formula
\begin{align}
\bfM_1(\eta)\bfg &= \bfM_0(\eta(\bfM_0\bfg)^\perp)-\nabla(\eta\nablac\bfg^\perp)+\alpha\eta(\bfM_0\bfg)^\perp \nonumber \\
&= -\frac{1}{D^2}\bfL\bfD\cdot\left(\eta\frac{1}{D^2}\bfL\bfD\cdot\bfg^\perp\right)+\bfD(\eta \bfD\cdot\bfg^\perp)
+\alpha\eta\frac{1}{D^2}\bfL^\perp \bfD\cdot\bfg^\perp \label{Formula for M1}
\end{align}
holds for each $\eta \in H^3({\mathbb R}^2)$ and each $\bfg \in H^3({\mathbb R}^2)^2$.
\end{lemma}
{\bf Proof.} 
Substituting the expansions \eqref{Taylor B}, \eqref{Taylor M} into \eqref{Flattened B BVP 1}--\eqref{Flattened B BVP 5} and \eqref{Defn of M tilde}
and equating constant terms shows that\linebreak
$\tilde\bfB^0 \in H^\frac{7}{2}(D_0)^3 \subseteq H^3(D_0)^3$ solves the
boundary-value problem \eqref{AB1}--\eqref{AB5}, while
equating terms which are linear in $\eta$ and making the Ansatz
$$\tilde\bfB^1=(v+1)\eta\partial_v\tilde\bfB^0+\tilde\bfC$$
leads to
$$\bfM_1(\eta)\tilde\bfg=-(\curl \tilde\bfC)_\mathrm{h}-(\curl (v+1)\eta\partial_v\tilde\bfB^0)_\mathrm{h}
-\eta(\partial_v\tilde\bfB^0)_\mathrm{h}^\perp+\nabla\eta^\perp\partial_v\tilde{B}^0_3-\nabla\eta\nablac(\tilde\bfB^0_\mathrm{h})^\perp\Big|_{v=0},$$
where
\begin{alignat*}{2}
 \curl \curl \tilde{\bfC} &= \alpha  \curl \tilde{\bfC} & & \mbox{in $D_0$}, \\
\Div \tilde{\bfC} &= 0 & & \mbox{in $D_0$}, \\
\tilde{\bfC} \wedge \bfe_3 &= \bfzero & & \mbox{at $v=-1$}, \\
\tilde\bfC\cdot\bfe_3 &= -\eta\partial_v\tilde{B}^0_3 & & \mbox{at $v=0$}, \\
\nablac\tilde{\bfC}_\mathrm{h}^\perp&=\nablac(-\eta(\partial_v\tilde\bfB^0)_\mathrm{h}-\nabla\eta \tilde{B}_3^0)^\perp \hspace{1cm} & & \mbox{at $v=0$}.
\end{alignat*}
Writing $\tilde\bfC=\tilde{\bfD} +\Grad\varphi$, where $\varphi \in H^3(D_0)$ is the unique solution of the boundary-value problem
\begin{alignat*}{2}
\Delta \varphi&=0 & & \mbox{in $D_0$},\\
\varphi &=0 & & \mbox{at $v=-1$}, \\
\varphi_v &=-\eta\partial_v\tilde{B}^0_3 \qquad\quad& & \mbox{at $v=0$},
\end{alignat*}
we find that $\tilde{\bfD} \in H^2(D_0)$ is the unique solution of the boundary-value problem
\begin{alignat*}{2}
 \curl \curl \tilde{\bfD} &= \alpha  \curl \tilde{\bfD} & & \mbox{in $D_0$}, \\
\Div \tilde{\bfD} &= 0 & & \mbox{in $D_0$}, \\
\tilde{\bfD} \wedge \bfe_3 &= \bfzero & & \mbox{at $v=-1$}, \\
\tilde{\bfD}\cdot\bfe_3 &= 0 & & \mbox{at $v=0$}, \\
\nablac(\tilde{\bfD})_\mathrm{h}^\perp&=\nablac(-\eta(\partial_v\tilde\bfB^0)_\mathrm{h}-\nabla\eta \tilde{B}_3^0)^\perp\hspace{1cm} & & \mbox{at $v=0$}
\end{alignat*}
and
\begin{equation}
\bfM_1(\eta)\tilde\bfg=-(\curl \tilde{\bfD})_\mathrm{h}-(\curl (v+1)\eta\partial_v\tilde\bfB^0)_\mathrm{h}
-\eta(\partial_v\tilde\bfB^0)_\mathrm{h}^\perp+\nabla\eta^\perp\partial_v\tilde{B}^0_3-\nabla\eta\nablac(\tilde\bfB^0_\mathrm{h})^\perp\Big|_{v=0}.
\label{pre-formula for M1}
\end{equation}

Comparing the boundary-value problem for $\tilde{\bfD}$ with \eqref{AB1}--\eqref{AB5}, we find that
\begin{align*}
(\udl{\curl \tilde\bfD})_\mathrm{h} &=\bfM_0\big(\eta(\udl{\partial_v\tilde\bfB}^0)_\mathrm{h}+\nabla\eta \udl{\tilde{B}}_3^0\big) \\
&=\bfM_0\big(\eta (\udl{\curl\tilde{\bfB}}^0)_\mathrm{h}^\perp+\nabla(\eta\udl{\tilde{B}}_3^0)\big)\\
&=-\bfM_0(\eta(\bfM_0\bfg)^\perp)
\end{align*}
because $\bfM \nabla (\cdot)=\bfzero$, and explicit calculations show that
$$-\nabla\eta \nablac(\tilde\bfB^0_\mathrm{h})^\perp\Big|_{v=0}=-\nabla(\eta\nablac(\tilde\bfB^0_\mathrm{h})^\perp)+\eta\nabla \big(\nablac(\tilde\bfB^0_\mathrm{h})^\perp\big)\Big|_{v=0}
=-\nabla(\eta\nablac\bfg^\perp)+\eta\nabla \big(\nablac(\tilde\bfB^0_\mathrm{h})^\perp\big)\Big|_{v=0}$$
and
$$-(\curl (v+1)\eta\tilde\bfB^0_v)_\mathrm{h}
-\eta(\tilde\bfB^0)_\mathrm{h}^\perp+\nabla\eta^\perp\tilde\bfB^0_{3v}+\eta\nabla \big(\nablac(\tilde\bfB^0_\mathrm{h})^\perp\big)\Big|_{v=0}
=\eta \big(\Delta \tilde\bfB^0_\mathrm{h}- \nabla \Div(\tilde{\bfB}^0)\big)^\perp\Big|_{v=0}=\eta (\Delta \tilde\bfB^0_\mathrm{h})^\perp\Big|_{v=0}.$$
The result follows by inserting these expressions into \eqref{pre-formula for M1} and noting that
$$\Delta \tilde\bfB^0_\mathrm{h}|_{v=0}=-\alpha(\udl{\curl\tilde{\bfB}}^0)_\mathrm{h}=\alpha (\bfM_0\bfg).\eqno{\Box}$$

\begin{remark}
This method leads to the loss of two derivatives in the individual terms in the formula for
$\bfM_1(\eta)$; the overall validity of the formula arises from subtle cancellations between the terms
(see Nicholls \& Reitich \cite[\S2.2]{NichollsReitich01a} for a discussion of this phenomenon in the context of the classical Dirichlet--Neumann operator). 
\end{remark}

Explicit expressions for the first few terms in the Taylor expansion
$$
\bfT(\eta) = \sum_{k=0}^\infty \bfT_k(\eta),\qquad \bfT_k(\eta)=\frac{1}{k!}\mathrm{d}\bfT^k[0](\eta^{(k)}),
$$
of $\bfT$ can be computed from the formula
$\bfT(\eta)=\bfM(\eta)\bfS(\eta)$ using \eqref{Formula for M0}, \eqref{Formula for M1}
and the corresponding expansion
$$
\bfS(\eta) = \sum_{k=1}^\infty \bfS_k(\eta)
$$
of $\bfS(\eta)$, where
$$
\bfS_k(\eta)\coloneqq \begin{cases}
\displaystyle(-1)^{\frac{k-1}{2}} \, \frac{\alpha^{k-1} \, \eta^{k}}{k!} \bfc^\perp,
& k=1,3,5,\ldots, \\[3mm]
\displaystyle(-1)^{\frac{k}{2}} \, \frac{\alpha^{k-1} \, \eta^{k}}{k!} \bfc,
& k=2,4,6,\ldots,
\end{cases}
$$
In particular, we find that
$$
\bfT_0 = \bfzero, \qquad \bfT_1(\eta)=\bfM_0 \bfS_1(\eta), \qquad \bfT_2(\eta)=\bfM_0 \bfS_2(\eta)+\bfM_1(\eta) \bfS_1(\eta),
$$
such that
\begin{align*}
\bfT_1(\eta) &= -\bfL \frac{\bfc\cdot \bfD}{D^2} \eta, \\
\bfT_2(\eta) &= \tfrac{1}{2}\alpha \bfL \frac{\bfc\cdot \bfD^\perp}{D^2}\eta^2
-\alpha\eta \bfL^\perp \frac{\bfc\cdot \bfD}{D^2}\eta +\bfL\frac{\bfD}{D^2}\cdot
\left(\eta \bfL\frac{\bfc\cdot\bfD}{D^2}\eta\right)-\bfD(\eta(\bfc\cdot\bfD)\eta).
\end{align*}

Turning to the Taylor expansion
$$
\JJ(\eta) = \sum_{k=0}^\infty \JJ_k(\eta),\qquad \JJ_k=\frac{1}{k!}\mathrm{d}^k\JJ[0](\eta^{(k)}),
$$
we conclude that
\begin{align}
\JJ_1(\eta) & = \bfT_1(\eta)\cdot \bfc + \eta - \beta \Delta \eta \nonumber \\
&=  \left(-\frac{1}{D^2}(\bfc\cdot\bfL)(\bfc\cdot \bfD) +1+\beta D^2\right)\eta, \nonumber \\[2mm]
\JJ_2(\eta) & = \tfrac{1}{2}|\bfT_1(\eta)|^2-\tfrac{1}{2}(\bfc \cdot \nabla\eta)^2+\bfT_2(\eta)\cdot \bfc+\alpha \eta\,\bfT_1(\eta)\cdot\bfc^\perp \nonumber \\
& = \tfrac{1}{2}\left|\bfL\frac{\bfc\cdot\bfD}{D^2}\eta\right|^2+\tfrac{1}{2}\alpha \frac{1}{D^2}(\bfc\cdot\bfL)(\bfc\cdot\bfD^\perp)\eta^2
+\frac{1}{D^2}(\bfc\cdot\bfL)\bfD\cdot\left(\eta \bfL \frac{\bfc\cdot\bfD}{D^2}\eta\right) \nonumber \\
& \qquad\qquad -\tfrac{1}{2}(\bfc\cdot\nabla\eta)^2+\bfc\cdot\nabla(\eta(\bfc\cdot\nabla\eta)) \label{Formula for J2}
\end{align}
and%
\begin{align*}
\JJ_{\ge 3}(\eta) & = \frac{1}{2}\big(2\bfT_1(\eta) + \bfT_{\ge 2}(\eta)\big) \cdot \bfT_{\ge 2}(\eta)
- \frac{ ( \alpha\bfS(\eta)\cdot\nabla\eta + \bfT(\eta) \cdot \nabla\eta)^2 }{2(1+|\nabla\eta|^2)}
-\frac{ \bfc\cdot\nabla\eta( \alpha\bfS(\eta)\cdot\nabla\eta + \bfT(\eta) \cdot \nabla\eta) }{1+|\nabla\eta|^2}\\
& \qquad\mbox{}+\frac{(\bfc\cdot\nabla\eta)^2|\nabla\eta|^2}{2(1+|\nabla\eta|^2)}
+\bfT_{\ge 3}(\eta)\cdot\bfc+\alpha\bfT_{\ge 2}(\eta)\cdot\bfS(\eta)+\alpha\bfT_1(\eta)\cdot\bfS_{\ge 2}(\eta)\\
& \qquad\mbox{}+
\beta \left(\frac{|\nabla\eta|^2\eta_x}{(1+|\nabla\eta|^2)^{\frac{1}{2}}(1+(1+|\nabla\eta|^2)^{\frac{1}{2}})}\right)_x+\beta \left(\frac{|\nabla\eta|^2\eta_y}{(1+|\nabla\eta|^2)^{\frac{1}{2}}(1+(1+|\nabla\eta|^2)^{\frac{1}{2}})}\right)_y,
\end{align*}
where
$$\JJ_{\ge 3}(\eta)=\sum_{k=3}^\infty \JJ_k(\eta), \quad \bfS_{\ge 2}(\eta)=\sum_{k=2}^\infty \bfS_k(\eta), \quad \bfT_{\ge 2}(\eta)=\sum_{k=2}^\infty \bfT_k(\eta), \quad \bfT_{\ge 3}(\eta)=\sum_{k=3}^\infty \bfT_k(\eta).$$

\section{The reduction procedure} \label{reduction}

In this section we reduce the equation
$$\JJ(\eta)=0$$
with
$$\bfc=(1-\varepsilon^2)\bfc_0$$
into a locally equivalent equation for $\eta_1$.
Clearly $\eta \in U$ solves this equation if and only if
\begin{align*}
\chi(\bfD)\JJ(\eta_1+\eta_2) & = 0,\\
\left(1-\chi(\bfD)\right)\JJ(\eta_1+\eta_2) & = 0,
\end{align*}
where $\eta_1=\chi(\bfD)\eta$, $\eta_2 = (1-\chi(\bfD))\eta$, which equations are given explicitly by
\begin{align}
g(\bfD)\eta_1+2&\varepsilon^2\frac{1}{D^2}(\bfc_0\cdot\bfD)(\bfc_0\cdot\bfL)\eta_1 \nonumber \\
& -\varepsilon^4\frac{1}{D^2}(\bfc_0\cdot\bfD)(\bfc_0\cdot\bfL)\eta_1
+\chi(\bfD)\big(\JJ_2(\eta_1+\eta_2) + \JJ_{\ge 3}(\eta_1+\eta_2)\big) =0, \label{eta1 eqn} \\
\hspace{-2mm}g(\bfD)\eta_2+2&\varepsilon^2\frac{1}{D^2}(\bfc_0\cdot\bfD)(\bfc_0\cdot\bfL)\eta_2 \nonumber \\
&-\varepsilon^4\frac{1}{D^2}(\bfc_0\cdot\bfD)(\bfc_0\cdot\bfL)\eta_2+\left(1-\chi(\bfD)\right)\big(\JJ_2(\eta_1+\eta_2) + \JJ_{\ge 3}(\eta_1+\eta_2)\big)=0, \label{eta2 eqn}
\end{align}
where
$$
g(\bfk)=-\frac{1}{|\bfk|^2}(\alpha(\bfc_0\cdot\bfk^\perp)(\bfc_0\cdot\bfk)+\mathtt{c}(|\bfk|^2)(\bfc_0\cdot\bfk)^2)+1+\beta |\bfk|^2;
$$
note that \eqref{eta1 eqn}, \eqref{eta2 eqn} hold in respectively $\chi(\bfD)H^1({\mathbb R}^2)$ and $(1-\chi(\bfD))H^1({\mathbb R}^2)$.
We proceed by solving \eqref{eta2 eqn} to determine $\eta_2$ as a function of $\eta_1$ and inserting
$\eta_2=\eta_2(\eta)$ into \eqref{eta1 eqn} to derive a reduced equation for $\eta_1$. To this end
we write \eqref{eta2 eqn} in the form
\begin{equation}
\eta_2 = \left(1-\chi(\bfD)\right)g(\bfD)^{-1}\AA(\eta_1,\eta_2),
\label{Rewritten eta2 eqn}
\end{equation}
where
\begin{equation}
\AA(\eta_1,\eta_2)=-2\varepsilon^2 \frac{1}{D^2}(\bfc_0\cdot\bfD)(\bfc_0\cdot\bfL)\eta_2+\varepsilon^4\frac{1}{D^2}(\bfc_0\cdot\bfD)(\bfc_0\cdot\bfL)\eta_2-
\left(1-\chi(\bfD)\right)\big(\JJ_2(\eta_1+\eta_2) + \JJ_{\ge 3}(\eta_1+\eta_2)\big).
\label{Formula for AA}
\end{equation}

\begin{proposition}
\label{prop:bounded mapping}
The mapping $\left(1-\chi(\bfD)\right)g(\bfD)^{-1}$ defines a bounded linear operator $H^1({\mathbb R}^2)\to H^3({\mathbb R}^2)$.
\end{proposition}
{\bf Proof.} This result follows from the facts that $(0,0)$ is a strict global minimum of $\tilde{g}(k_1,\tfrac{k_2}{k_1})$ with $\tilde{g}(0,0)=0$ and
that $g(\bfk) \gtrsim |\bfk|^2$ as $|\bfk| \rightarrow \infty$.\qed

The next step is to estimate the nonlinear terms on the right-hand side of equation \eqref{Rewritten eta2 eqn}. The requisite estimates for $\JJ_2(\eta)$ are obtained by examining the explicit formula
$$\JJ_2(\eta)=b(\eta,\eta)-2\varepsilon^2 b(\eta,\eta)+\varepsilon^4 b(\eta,\eta),$$
where
\begin{align}
b(v,w)
& = \tfrac{1}{2}\left(\bfL\frac{\bfc_0\cdot\bfD}{D^2}v\right)\!\cdot\!\left(\bfL\frac{\bfc_0\cdot\bfD}{D^2}w\right)+\tfrac{1}{2}\alpha \frac{1}{D^2}(\bfc_0\cdot\bfL)(\bfc_0\cdot\bfD^\perp)vw \nonumber \\
& \qquad\qquad 
+\frac{1}{2D^2}(\bfc_0\cdot\bfL)\bfD\cdot\left(v \bfL \frac{\bfc_0\cdot\bfD}{D^2}w\right)+\frac{1}{2D^2}(\bfc_0\cdot\bfL)\bfD\cdot\left(w \bfL \frac{\bfc_0\cdot\bfD}{D^2}v\right) \nonumber \\
& \qquad\qquad +\tfrac{1}{2}((\bfc_0\cdot\bfD)v)((\bfc_0\cdot\bfD)w)-\tfrac{1}{2}\bfc_0\cdot\bfD(v(\bfc_0\cdot\bfD)w)-\tfrac{1}{2}\bfc_0\cdot\bfD(w(\bfc_0\cdot\bfD)v)
\label{Formula for m}
\end{align}
(see equation \eqref{Formula for J2}).

\begin{lemma} \label{Estimate for m}
The estimate
$\|b(v,w)\|_1 \lesssim \|v\|_\ZZ \|w\|_3$
holds for each $v$, $w\in H^3({\mathbb R}^2)$.
\end{lemma}
{\bf Proof.} We estimate each of the terms in the formula for $m$, observing that
$$\left|\frac{1}{|\bfk|^2}(\alpha \bfk^\perp + \mathtt{c}(|\bfk|^2) \bfk) \bfc_0\cdot\bfk\right|,\ 
\left|\frac{1}{|\bfk|^2}(\alpha \bfc_0\cdot\bfk^\perp + \mathtt{c}(|\bfk|^2) \bfc_0\cdot\bfk) \bfc_0\cdot\bfk^\perp\right|,\ 
\left|\frac{1}{|\bfk|^2}(\alpha \bfc_0\cdot\bfk^\perp + \mathtt{c}(|\bfk|^2) \bfc_0\cdot\bfk) \bfk\right| \lesssim \langle \bfk \rangle$$
and that
$$\|f(\bfD)\hat{v}_1\|_{n,\infty} \leq \|f(\bfk) \langle \bfk\rangle^n \hat{v}_1\|_{L^1({\mathbb R})} \lesssim \|\hat{v}_1\|_{L^1({\mathbb R})}$$
for all bounded multipliers $f$ because $\hat{v}_1$ has compact support.
We find that
\begin{align*}
\qquad\left\| \left(\bfL\frac{\bfc_0\cdot\bfD}{D^2}v\right)\!\cdot\!\left(\bfL\frac{\bfc_0\cdot\bfD}{D^2}w\right) \right\|_1
& \lesssim \left\|\bfL\frac{\bfc_0\cdot\bfD}{D^2}v_1\right\|_{1,\infty}\left\|\bfL\frac{\bfc_0\cdot\bfD}{D^2}w\right\|_1
+\left\| \bfL\frac{\bfc_0\cdot\bfD}{D^2}v_2\right\|_2\left\|\bfL\frac{\bfc_0\cdot\bfD}{D^2}w\right\|_2 \\
& \leq (\|\hat{v}_1\|_{L^1({\mathbb R}^2)}+\|v_2\|_3)\|w\|_3,\\
\\
\left\| \frac{1}{D^2}(\bfc_0\cdot\bfL)(\bfc_0\cdot\bfD^\perp)vw \right\|_1 & \lesssim \|vw\|_2 \\
& \lesssim (\|v_1\|_{2,\infty}+ \|v_2\|_2)\|w\|_2 \\
& \leq (\|\hat{v}_1\|_{L^1({\mathbb R}^2)}+\|v_2\|_3)\|w\|_3,\\
\\
\left\|\frac{1}{D^2}(\bfc_0\cdot\bfL)\bfD\cdot\left(v \bfL \frac{\bfc_0\cdot\bfD}{D^2}w\right)\right\|_1
& \lesssim \left\|v \bfL \frac{\bfc_0\cdot\bfD}{D^2}w\right\|_2 \\
& \lesssim (\|v_1\|_{2,\infty}+ \|v_2\|_2)\left\|\bfL \frac{\bfc_0\cdot\bfD}{D^2}w\right\|_2 \\
& \lesssim (\|\hat{v}_1\|_{L^1({\mathbb R}^2)}+\|v_2\|_3)\|w\|_3,\\
\\
\left\|\frac{1}{D^2}(\bfc_0\cdot\bfL)\bfD\cdot\left(w \bfL \frac{\bfc_0\cdot\bfD}{D^2}v\right)\right\|_1
& \lesssim \left\|w \bfL \frac{\bfc_0\cdot\bfD}{D^2}v\right\|_2 \\
& \lesssim \|w\|_2 \left(\left\|\bfL \frac{\bfc_0\cdot\bfD}{D^2}v_1\right\|_{2,\infty} + \left\|\bfL \frac{\bfc_0\cdot\bfD}{D^2}v_2\right\|_2\right) \\
& \lesssim (\|\hat{v}_1\|_{L^1({\mathbb R}^2)}+\|v_2\|_3)\|w\|_3,\\
\\
\|(\bfc_0\cdot\bfD)v)((\bfc_0\cdot\bfD)w)\|_1
& \lesssim (\|(\bfc_0\cdot\bfD)v_1\|_{1,\infty}\|\bfc_0\cdot\bfD)w\|_1 + \|(\bfc_0\cdot\bfD)v\|_2\|(\bfc_0\cdot\bfD)w\|_2 \\
& \lesssim (\|\hat{v}_1\|_{L^1({\mathbb R}^2)}+\|v_2\|_3)\|w\|_3,\\
\\
\|\bfc_0\cdot\bfD(v(\bfc_0\cdot\bfD)w)\|_1
& \lesssim \|v(\bfc_0\cdot\bfD)w)\|_2 \\
& \lesssim (\|v_1\|_{2,\infty} + \|v_2\|_2) \|(\bfc_0\cdot\bfD)w\|_2 \\
& \lesssim (\|\hat{v}_1\|_{L^1({\mathbb R}^2)}+\|v_2\|_3)\|w\|_3,\\
\\
\|\bfc_0\cdot\bfD(w(\bfc_0\cdot\bfD)v)\|_1
& \lesssim \|w(\bfc_0\cdot\bfD)v\|_2 \\
& \lesssim \|w\|_2 (\|(\bfc_0\cdot\bfD)v_1\|_{2,\infty} + \|(\bfc_0\cdot\bfD)v_2\|_2) \\
& \lesssim (\|\hat{v}_1\|_{L^1({\mathbb R}^2)}+\|v_2\|_3)\|w\|_3.\hspace{2.5in}\Box
\end{align*}

\begin{corollary} \label{Estimate for J2}
The estimates
$$
\|\JJ_2(\eta)\|_1 \lesssim \|\eta\|_\ZZ\|\eta\|_3
$$
and
$$
\|\mathrm{d}\JJ_2[\eta](u)\|_1 \lesssim \|\eta\|_3\|u\|_\ZZ, \qquad \|\mathrm{d}\JJ_2[\eta](u)\|_1 \lesssim \|\eta\|_\ZZ\|u\|_3
$$
hold for all $\eta$, $u \in H^3({\mathbb R}^2)$.
\end{corollary}

\begin{lemma} \label{Estimate for Jc}
The estimates
\begin{align*}
\|\JJ_{\ge 3}(\eta)\|_1 & \lesssim \|\eta\|_\ZZ^2\|\eta\|_3,  \\[1mm]
\|\mathrm{d}\JJ_{\ge 3}[\eta](u)\|_1 & \lesssim \|\eta\|_\ZZ^2 \|u\|_3 + \|\eta\|_{\ZZ}  \|\eta\|_3\|u\|_\ZZ
\end{align*}
hold for all $\eta \in U$ and $u \in H^3({\mathbb R}^2)$.
\end{lemma}
{\bf Proof.} Writing
\begin{align*}
\bfT_{\ge 2}(\eta)&=\bfM_0\bfS_{\ge 2}(\eta)+\big(\bfM_1(\eta)+\bfM_{\ge 2}(\eta)\big)\bfS_1(\eta), \\
\bfT_{\ge 3}(\eta)&=\bfM_0\bfS_{\ge 3}(\eta)+\bfM_{\ge 2}(\eta)\bfS_1(\eta)+\big(\bfM_1(\eta)+\bfM_{\ge 2}(\eta)\big)\bfS_2(\eta),
\end{align*}
we find by Theorem \ref{main anal thm}(i) and the fact that $\bfM_0 \in L(H^{\frac{5}{2}}({\mathbb R}^2)^2,H^{\frac{3}{2}}({\mathbb R}^2)^2)$ that
\begin{align*}
\|\bfT_{\ge 2}(\eta)\|_\frac{3}{2} &\lesssim \left\|\eta\frac{\bfS_{\ge 2}(\eta)}{\eta}\right\|_\frac{5}{2}+\|\eta\|_\ZZ\|\bfS_1(\eta)\|_\frac{5}{2}\\
& \lesssim (\|\eta_1\|_{3,\infty}+\|\eta_2\|_3)\left\|\frac{\bfS_{\ge 2}(\eta)}{\eta}
\right\|_3+\|\eta\|_\ZZ\|\bfS_1(\eta)\|_3\\
& \lesssim \|\eta\|_\ZZ\|\eta\|_3,
\\
\|\bfT_{\ge 3}(\eta)\|_\frac{3}{2} &\lesssim \left\|\eta\frac{\bfS_{\ge 3}(\eta)}{\eta}\right\|_\frac{5}{2}+\|\eta\|_\ZZ^2\|\bfS_1(\eta)\|_\frac{5}{2}+\|\eta\|_\ZZ\left\|\eta\frac{S_2(\eta)}{\eta}\right\|_\frac{5}{2} \\
& \lesssim \|\eta\|_\ZZ^2\|\eta\|_3
\end{align*}
and hence that
\begin{align*}
\|\bfT_1(\eta)\cdot\bfT_{\ge 2}(\eta)\|_1
&\lesssim
\|\bfM_0\bfS_1(\eta_1)\|_{1,\infty}\|\bfT_{\ge 2}(\eta)\|_1+\|\bfM_0\bfS_1(\eta_2)\|_\frac{3}{2}\|\bfT_{\ge 2}(\eta)\|_\frac{3}{2}\\
&\lesssim
\|\hat{\eta}_1\|_{L^1({\mathbb R}^2)}\|\bfT_{\ge 2}(\eta)\|_\frac{3}{2}+\|\eta_2\|_\frac{5}{2}\|\bfT_{\ge 2}(\eta)\|_\frac{3}{2}\\
&\lesssim\|\eta\|_\ZZ^2\|\eta\|_3,
\\[2mm]
\|\bfT_{\ge 2}(\eta)\cdot\bfT_{\ge 2}(\eta)\|_1
&\lesssim
\|\bfT_{\ge 2}(\eta)\|_\frac{3}{2}^2\\
&\lesssim\|\eta\|_\ZZ^2\|\eta\|_3^2,
\\[2mm]
\|\bfT_{\ge 2}(\eta)\cdot\bfS(\eta)\|_1 &\lesssim \|\bfT_{\ge 2}(\eta)\|_1\|\bfS(\eta)\|_{1,\infty} \\
& \lesssim \|\eta\|_\ZZ\|\eta\|_3(\|\eta_1\|_{1,\infty} +\|\eta_2\|_{1,\infty}) \\
& \lesssim \|\eta\|_\ZZ\|\eta\|_3(\|\eta_1\|_{1,\infty} +\|\eta_2\|_3) \\
& \lesssim \|\eta\|_\ZZ^2\|\eta\|_3,
\\[2mm]
\ \|\bfT_1(\eta)\cdot\bfS_{\ge 2}(\eta)\|_1
& \lesssim \|\bfT_1(\eta)\|_1\|\bfS_{\ge 2}(\eta)\|_{1,\infty} \\
& \lesssim \|\eta\|_3(\|\eta_1\|_{1,\infty} +\|\eta_2\|_{1,\infty})^2 \\
& \lesssim \|\eta\|_3\|\eta\|_\ZZ^2.
\end{align*}%

Furthermore
\begin{align*}
\left\| \frac{(\bfc \cdot \nabla \eta)^2|\nabla \eta|^2}{1+|\nabla \eta|^2}\right\|_1
& \lesssim \|\nabla \eta\|_0 \left\| \frac{|\nabla \eta|^3}{1+|\nabla \eta|^2}\right\|_\infty + \left\| \nabla \left(\frac{(\bfc \cdot \nabla \eta)^2|\nabla \eta|^2}{1+|\nabla \eta|^2}\right)\right\|_0 \\
& \lesssim \|\nabla \eta\|_0 \|\nabla\eta\|_\infty^3 + \|\nabla \eta\|_\infty^3 (\|\eta_{xx}\|_0+\|\eta_{xy}\|_0+\|\eta_{yy}\|_0) \\
& \lesssim \|\eta\|_3 \|\eta\|_\ZZ^3, \\
\\
\end{align*}
and

\begin{align*}
\left\|\frac{ \bfc\cdot\nabla\eta( \alpha\bfS(\eta) \cdot \nabla\eta + \bfT(\eta) \cdot \nabla\eta) }{1+|\nabla\eta|^2}\right\|_1\hspace{-3.5cm}\\
& \lesssim (\|\bfS(\eta)\|_1+\|\bfT(\eta)\|_1)\left\|\frac{|\nabla\eta|^2}{1+|\nabla\eta|^2}\right\|_\infty+(\|\bfS(\eta)\|_\infty+\|\bfT(\eta)\|_\infty)
\left\|\nabla\left(\frac{|\nabla\eta|^2}{1+|\nabla\eta|^2}\right)\right\|_0\\
& \lesssim (\|\bfS(\eta)\|_1+\|\bfT(\eta)\|_1)\|\nabla\eta\|_\infty^2+(\|\bfS(\eta)\|_\infty+\|\bfT(\eta)\|_\infty)
\|\nabla \eta\|_\infty^2 (\|\eta_{xx}\|_0+\|\eta_{xy}\|_0+\|\eta_{yy}\|_0)\\
& \lesssim \|\eta\|_3 \|\eta\|_\ZZ^2, \\
\\
\left\|\frac{ ( \alpha\bfS(\eta) \cdot \nabla\eta + \bfT(\eta) \cdot \nabla\eta)^2 }{2(1+|\nabla\eta|^2)}\right\|_1 \hspace{-3cm}\\
& \lesssim (\|\bfS(\eta)\|_1+\|\bfT(\eta)\|_1)(\|\bfS(\eta)\|_\infty+\|\bfT(\eta)\|_\infty)\left\|\frac{|\nabla\eta|^2}{1+|\nabla\eta|^2}\right\|_\infty\\
& \qquad\quad\mbox{}
+(\|\bfS(\eta)\|_\infty+\|\bfT(\eta)\|_\infty)^2\left\|\nabla\left(\frac{|\nabla\eta|^2}{1+|\nabla\eta|^2}\right)\right\|_0\\
& \lesssim (\|\bfS(\eta)\|_1+\|\bfT(\eta)\|_1)(\|\bfS(\eta)\|_\infty+\|\bfT(\eta)\|_\infty)\|\nabla\eta\|_\infty^2\\
& \qquad\quad\mbox{}
+(\|\bfS(\eta)\|_\infty+\|\bfT(\eta)\|_\infty)^2
\|\nabla \eta\|_\infty^2 (\|\eta_{xx}\|_0+\|\eta_{xy}\|_0+\|\eta_{yy}\|_0)\\
& \lesssim  \|\eta\|_3\|\eta\|_\ZZ^3, 
\end{align*}
where we have used the estimates
$$
\|\bfS(\eta)\|_1,\ \|\bfT(\eta)\|_1 \leq \|\eta\|_3, \qquad
\|\bfS(\eta)\|_\infty \lesssim \|\eta\|_\infty \lesssim \|\eta\|_\ZZ
$$
and
\begin{align*}
\|\bfT(\eta)\|_\infty &\leq \|\bfT_1(\eta)\|_\infty + \|\bfT_{\geq 2}(\eta)\|_\infty \\
& \lesssim \|\bfM_0\bfS_1(\eta_1)\|_\infty + \|\bfM_0\bfS_1(\eta_2)\|_{\frac{3}{2}} +\|\bfT_{\geq 2}(\eta)\|_{\frac{3}{2}} \\
& \lesssim \|\hat{\eta}_1\|_{L^1({\mathbb R}^2)} + \|\eta_2\|_{\frac{5}{2}} + \|\eta\|_{\ZZ}\|\eta\|_3 \\
& \lesssim \|\eta\|_\ZZ.
\end{align*}

Finally, we note that
\begin{eqnarray*}
\lefteqn{\left(\frac{|\nabla\eta|^2\eta_x}{(1+|\nabla\eta|^2)^{\frac{1}{2}}(1+(1+|\nabla\eta|^2)^{\frac{1}{2}})}\right)_{\!\!x}+\left(\frac{|\nabla\eta|^2\eta_y}{(1+|\nabla\eta|^2)^{\frac{1}{2}}(1+(1+|\nabla\eta|^2)^{\frac{1}{2}})}\right)_{\!\!y}} \hspace{1cm}\\
&=f_1(\eta_x, \eta_y)\eta_{xx}+f_2(\eta_x, \eta_y)\eta_{xy}+f_3(\eta_x,\eta_y)\eta_{yy},\hspace{1in}
\end{eqnarray*}
where $f_1, f_2, f_3$ are analytic functions with zeros of order two at the origin. Estimating
 \begin{align*}
 \|f_1(\eta_x,\eta_y)\eta_{xx}\|_0 &\lesssim \|\nabla\eta\|_\infty^2\|\eta_{xx}\|_0  \lesssim \|\eta\|_\ZZ^2\|\eta\|_3,
 \\[2mm]
 \|\nabla (f_1(\eta_x,\eta_y)\eta_{xx})\|_0 
 & \leq \|\partial_1 f_1(\eta_x,\eta_y) \eta_{xx} \nabla \eta_x +\partial_2 f_1(\eta_x,\eta_y) \eta_{xx} \nabla \eta_y \|_0+ \|f_1(\eta_x,\eta_y)\nabla\eta_{xx})\|_0 \\
 & \lesssim \|\nabla\eta\|_\infty(\|\nabla \eta_x\|_{L^4({\mathbb R}^2)}+\|\nabla \eta_y\|_{L^4({\mathbb R}^2)}) \|\eta_{2xx}\|_{L^4({\mathbb R}^2)}\\
 & \qquad \mbox{}+\|\nabla\eta\|_\infty(\|\nabla \eta_x\|_0+\|\nabla \eta_y\|_0) \|\eta_{1xx}\|_\infty+\|\nabla\eta\|_\infty^2\|\nabla\eta_{xx}\|_0 \\
 &  \lesssim \|\eta\|_\ZZ^2\|\eta\|_3,
 \end{align*}
 in which the last fine follows by the continuous embedding $H^1({\mathbb R}^2) \subseteq L^4({\mathbb R}^2)$,
 and $f_2(\eta_x, \eta_y)\eta_{xy}$, $f_3(\eta_x,\eta_y)\eta_{yy}$ similarly, we conclude that
$$
\left\|\left(\frac{|\nabla\eta|^2\eta_x}{(1+|\nabla\eta|^2)^{\frac{1}{2}}(1+(1+|\nabla\eta|^2)^{\frac{1}{2}})}\right)_{\!\!x}+\left(\frac{|\nabla\eta|^2\eta_y}{(1+|\nabla\eta|^2)^{\frac{1}{2}}(1+(1+|\nabla\eta|^2)^{\frac{1}{2}})}\right)_{\!\!y}\right\|_1
\lesssim \|\eta\|_\ZZ^2\|\eta\|_3.
$$

The estimates for the derivatives are obtained in a similar fashion.\qed\pagebreak

We proceed by solving \eqref{Rewritten eta2 eqn} for $\eta_2$ as a function of $\eta_1$ using the following fixed-point theorem,
which is proved by a straightforward application of the contraction mapping principle.

\begin{theorem}
\label{thm:fixed-point}
Let $\XX_1$, $\XX_2$ be Banach spaces, $X_1$, $X_2$ be closed, convex sets in, respectively, $\XX_1$, $\XX_2$ containing the origin and $\GG\colon X_1\times X_2 \to \XX_2$ be a smooth mapping. Suppose there exists a continuous mapping $r\colon X_1\to [0,\infty)$ such that
$$
\|\GG(x_1,0)\|\le \tfrac{1}{2}r, \quad \|\mathrm{d}_2 \GG[x_1,x_2]\|\le \tfrac{1}{3}
$$ 
for each $x_2\in \overline B_r(0)\subseteq X_2$ and each $x_1\in X_1$.

Under these hypotheses there exists for each $x_1\in X_1$ a unique solution $x_2=x_2(x_1)$ of the fixed-point equation
$x_2=\GG(x_1,x_2)$
satisfying $x_2(x_1)\in \overline B_r(0)$. Moreover $x_2(x_1)$ is a smooth function of $x_1\in X_1$ with
$$
\|\mathrm{d} x_2[x_1]\|_{L(\XX_1,\XX_2)}\le 2\|\mathrm{d}_1 \GG[x_1, x_2(x_1)]\|.
$$
\end{theorem}

We apply Theorem \ref{thm:fixed-point} to equation \eqref{Rewritten eta2 eqn} with
$$\XX_1=\chi(\bfD)H^3({\mathbb R}^2),
\qquad
\XX_2=\left(1-\chi(\bfD)\right)H^3({\mathbb R}^2),$$
equipping
$\XX_1$ with the scaled norm $\nn \cdot \nn$ defined in \eqref{scaled norm}
and $\XX_2$ with the usual norm for $H^3({\mathbb R}^2)$, and taking
$$
X_1=\{\eta_1\in \XX_1 \colon \nn \eta_1\nn \le R_1\}, \qquad
X_2=\{\eta_2\in \XX_2 \colon \| \eta_2\|_3 \le R_2\};
$$
the function $\GG$ is given by the right-hand side of \eqref{Rewritten eta2 eqn}. Recall that $\JJ$ is an analytic function
$U \rightarrow H^1({\mathbb R}^2)$ (see equation \eqref{Defn of U with M}). Using Proposition \ref{sup estimate}
we can guarantee that $\|\hat{\eta}_1\|_{L^1({\mathbb R}^2)} < \tfrac{1}{2}M$ for all
$\eta_1 \in X_1$ for an arbitrarily large value
of $R_1$; the value of $R_2$
is then constrained by the requirement that $\|\eta_2\|_3 < \tfrac{1}{2}M$ for all $\eta_2 \in X_2$.

We proceed by estimating each term appearing in the formula \eqref{Formula for AA} for $\AA$
using Corollary \ref{Estimate for J2}, Lemma \ref{Estimate for Jc},
together with
Proposition \ref{sup estimate} and the estimates
$$\|\eta\|_\ZZ \lesssim \varepsilon \nn \eta_1 \nn + \|\eta_2\|_3, \qquad
\|\eta\|_3 \lesssim \nn \eta_1 \nn + \|\eta_2\|_3$$
for $\eta \in H^3({\mathbb R}^2)$.

\begin{lemma}
\label{lem:3 estimates}
The estimates
\begin{list}{(\roman{count})}{\usecounter{count}}
\item
$\|\AA(\eta_1,\eta_2)\|_1\lesssim \varepsilon \nn \eta_1\nn^2+\nn \eta_1\nn \|\eta_2\|_3 
+\|\eta_2\|_3^2+\varepsilon^2 \|\eta_2\|_3$,
\item
$\|\mathrm{d}_1\AA[\eta_1,\eta_2]\|_{L(\XX_1,H^1({\mathbb R}^2))}\lesssim \varepsilon \nn \eta_1\nn+ \|\eta_2\|_3$,
\item
$\|\mathrm{d}_2\AA[\eta_1,\eta_2]\|_{L(\XX_2,H^1({\mathbb R}^2))}\lesssim \varepsilon \nn \eta_1\nn+ \|\eta_2\|_3+\varepsilon^2$
\end{list}
hold for each $\eta_1\in X_1$ and $\eta_2\in X_2$.
\end{lemma}
\begin{theorem} \label{Estimate for eta2}
Equation \eqref{Rewritten eta2 eqn} has a unique solution $\eta_2 \in
X_2$ which depends smoothly upon $\eta_1 \in X_1$ and satisfies the estimates
$$
\|\eta_2(\eta_1)\|_3 \lesssim \varepsilon \nn \eta_1\nn^2,  \quad
\|\mathrm{d}\eta_2[\eta_1]\|_{L(\XX_1,\XX_2)} \lesssim \varepsilon \nn \eta_1\nn.
$$
\end{theorem}
{\bf Proof.}
Choosing $R_2$ and $\varepsilon$ sufficiently small, one finds $r>0$ such that
$\|\GG(\eta_1,0)\|_3 \leq \tfrac{1}{2}r$ and\linebreak $\|\mathrm{d}_2 \GG[\eta_1,\eta_2]\|_{L(\XX_2,\XX_2)} \leq \tfrac{1}{3}$
for $\eta_1 \in X_1$, $\eta_2 \in X_2$, and Theorem \ref{thm:fixed-point} asserts that equation \eqref{Rewritten eta2 eqn} has a unique solution $\eta_2 \in X_2$ which depends smoothly
upon $\eta_1 \in X_1$. More precise estimates are obtained by choosing $C>0$ so that
$\|\GG(\eta_1,0)\|_3 \leq C\varepsilon\nn \eta_1\nn^2$ for $\eta_1 \in X_1$
and writing $r(\eta_1)=2C\varepsilon \nn \eta_1\nn^2$, so that
$$
\|\mathrm{d}_1 \GG[\eta_1,\eta_2]\|_{L(\XX_1,\XX_2)}\lesssim \varepsilon \nn \eta_1 \nn, \quad
\|\mathrm{d}_2 \GG[\eta_1,\eta_2]\|_{L(\XX_2,\XX_2)}\lesssim 1
$$
for $\eta_1 \in X_1$, $\eta_2 \in \overline{B}_{r(\eta_1)}(0) \subseteq X_2$, and the stated estimates for $\eta_2(\eta_1)$ follow
from Theorem \ref{thm:fixed-point}.\qed\pagebreak

Inserting $\eta_2=\eta_2(\eta_1)$ into \eqref{eta1 eqn} yields the reduced equation
\begin{equation}
g(\bfD)\eta_1+2\varepsilon^2\frac{1}{D^2}(\bfc_0\cdot\bfD)(\bfc_0\cdot\bfL)\eta_1-\varepsilon^4\frac{1}{D^2}(\bfc_0\cdot\bfD)(\bfc_0\cdot\bfL)\eta_1
+\chi(\bfD)\big(\JJ_2(\eta_1+\eta_2(\eta_1)) + \JJ_{\ge 3}(\eta_1+\eta_2(\eta_1))\big)=0
\label{Red eqn 1}
\end{equation}
for $\eta_1$, which holds in $\chi(\bfD)H^1({\mathbb R}^2)$. This equation is invariant under the reflection $\eta_1(x,y) \mapsto \eta_1(-x,-y)$; a familiar argument shows that it is inherited from the
corresponding invariance $\eta_1(x,y) \mapsto \eta_1(-x,-y)$,\linebreak
$\eta_2(x,y) \mapsto \eta_2(-x,-y)$,
of \eqref{eta1 eqn}, \eqref{eta2 eqn} when applying Theorem \ref{thm:fixed-point}.

\section{Derivation of the reduced equation} \label{calculation of reduced eqn}
In this section we compute the leading-order terms in the reduced equation \eqref{Red eqn 1}.
The first step is to write
\begin{align*}
\JJ_2(\eta_1+\eta_2(\eta_1))&=b(\eta_1,\eta_1)-2\varepsilon^2b(\eta_1,\eta_1)+\varepsilon^4b(\eta_1,\eta_1) \\
& \qquad \mbox{}+b(\eta_1,\eta_2(\eta_1))-2\varepsilon^2b(\eta_1,\eta_2(\eta_1))+\varepsilon^4b(\eta_1,\eta_2(\eta_1)) \\
& \qquad \mbox{}+b(\eta_2(\eta_1),\eta_2(\eta_1))-2\varepsilon^2b(\eta_2(\eta_1),\eta_2(\eta_1))+\varepsilon^4b(\eta_2(\eta_1),\eta_2(\eta_1))
\end{align*}
and examine each of the terms on the right-hand side of this expression individually. The first term is handled by approximating the Fourier-multiplier operators by constants according to Lemma \ref{Multiplier to const} below. The order-of-magnitude estimates in this section are computed with
respect to the $L^2({\mathbb R}^2)$-norm, which is equivalent to the $H^s({\mathbb R}^2)$-norm on the space
$\chi(\bfD)H^s({\mathbb R}^2)$ for any $s \geq 0$.

\begin{lemma} \label{Multiplier to const}
The estimates
\begin{itemize}
\item[(i)] $\eta_{1x}=O(\varepsilon\nn\eta_1\nn)$,
\item[(ii)] $\eta_{1y}=O(\varepsilon\nn\eta_1\nn)$,
\item[(iii)] $\bfL\dfrac{\bfc_0\cdot\bfD}{D^2}\eta_1=\begin{pmatrix}\alpha \cot\alpha c_{0,1} \\ -\alpha c_{0,1} \end{pmatrix}\eta_1+O(\varepsilon\nn\eta_1\nn)$,
\item[(iv)] $\dfrac{1}{D^2}(\bfc_0\cdot\bfL)(\bfc_0\cdot\bfD)\eta_1=\underbrace{\alpha c_{0,1}(-c_{0,2}+c_{0,1}\cot\alpha)}_{\displaystyle=1} \eta_1+O(\varepsilon\nn\eta_1\nn)$,
\item[(v)] $\dfrac{1}{D^2}(\bfc_0\cdot\bfL)(\bfc_0\cdot\bfD^\perp)\eta_1\rho_1=\alpha c_{0,2}(c_{0,2}-c_{0,1}\cot\alpha)\eta_1\rho_1+O(\varepsilon^2\nn\eta_1\nn\nn\rho_1\nn)+B_1(\bfD)\eta_1\rho_1$,
\item[(vi)] $\dfrac{1}{D^2}(\bfc_0\cdot\bfL)\bfD\cdot(\eta_1\rho_1\bfw)=\alpha(-c_{0,2}+c_{0,1}\cot\alpha)\eta_1\rho_1w_1+O(\varepsilon^2\nn\eta_1\nn\nn\rho_1\nn)+B_2(\bfD)\eta_1\rho_1$,
\end{itemize}
hold for all $\eta_1$, $\rho_1 \in \XX_1$, where $\bfw$ is a vector-valued constant,
$$|B_j(\bfk)| \lesssim |\tfrac{k_2}{k_1}|(1+\tfrac{k_2^2}{k_1^2})^{-1}$$
and $c_{0,1}=c_0\cos\tfrac{1}{2}\alpha$, $c_{0,2} = -c_0\sin\tfrac{1}{2}\alpha$.
\end{lemma}
{\bf Proof.} Parts (i)--(iv) follow from the calculations
\begin{align*}
\|\eta_{1x}\|_0^2&=\|k_1\hat{\eta}_1\|_0^2 \\
&\leq \varepsilon^2\nn\eta_1\nn^2,\\
\\
\|\eta_{1z}\|_0^2&=\|k_2\hat{\eta}_1\|_0^2\\
&=\|k_1\tfrac{k_2}{k_1}\hat{\eta}_1\|_0^2\\
&\leq \delta^2\|k_1\hat{\eta}_1\|_0^2\\
&\lesssim\varepsilon^2\nn\eta_1\nn^2
\end{align*}
and
\begin{align*}
\Big\|\bfL&\dfrac{\bfc_0\cdot\bfD}{D^2}\eta_1-\begin{pmatrix}\alpha \cot\alpha\, c_{0,1} \\ -\alpha\, c_{0,1} \end{pmatrix}\eta_1\Big\|_0^2\\
&=
\bigg\|\underbrace{\bigg(\frac{1}{1+\frac{k_2^2}{k_1^2}}\left(\alpha\tfrac{k_2}{k_1}+\mathtt{c}(|\bfk|^2)\right)(c_{0,1}+c_{0,2}\tfrac{k_2}{k_1})-\alpha\cot\alpha \, c_{0,1}\bigg)}_{\displaystyle = O(|(k_1,\tfrac{k_2}{k_1}|)}\hat{\eta}_1\bigg\|_0^2\\
&\qquad\mbox{}+
\bigg\|\underbrace{\bigg(\frac{1}{1+\frac{k_2^2}{k_1^2}}\left(-\alpha+\mathtt{c}(|\bfk|^2)\tfrac{k_2}{k_1}\right)(c_{0,1}+c_{0,2}\tfrac{k_2}{k_1})+\alpha c_{0,1}\bigg)}_{\displaystyle = O(|(k_1,\tfrac{k_2}{k_1}|)}\hat{\eta}_1\bigg\|_0^2\\
&\lesssim\|k_1\hat{\eta}_1\|_0^2+\|\tfrac{k_2}{k_1}\hat{\eta}_1\|_0^2\\
&\lesssim\varepsilon^2\nn\eta_1\nn^2,\\
\\
\bigg\|\dfrac{1}{D^2}&(\bfc_0\cdot\bfL)(\bfc_0\cdot\bfD)\eta_1-\alpha c_{0,1}(-c_{0,2}+c_{0,1}\cot\alpha)\eta_1\Big\|_0^2\\
&=
\bigg\|\underbrace{\bigg(\frac{1}{1+\frac{k_2^2}{k_1^2}}\left(\alpha (c_{0,1}\tfrac{k_2}{k_1}-c_{0,2})+\mathtt{c}(|\bfk|^2)(c_{0,1}+c_{0,2}\tfrac{k_2}{k_1})\right)(c_{0,1}+c_{0,2}\tfrac{k_2}{k_1})
-\alpha c_{0,1}(-c_{0,2}+c_{0,1}\cot\alpha)\bigg)}_{\displaystyle = O(|(k_1,\tfrac{k_2}{k_1}|)}\hat{\eta}_1\bigg\|_0^2 \\
&\lesssim\|k_1\hat{\eta}_1\|_0^2+\|\tfrac{k_2}{k_1}\hat{\eta}_1\|_0^2\\
&\lesssim\varepsilon^2\nn\eta_1\nn^2.
\end{align*}
Turning to parts (v) and (vi), note that
\begin{align*}
&\Big\|\dfrac{1}{D^2}(\bfc_0\cdot\bfL)(\bfc_0\cdot\bfD^\perp)\eta_1\rho_1-\alpha c_{0,2}(c_{0,2}-c_{0,1}\cot\alpha)\eta_1\rho_1\Big\|_0^2\\
& =\!
\bigg\|\underbrace{\!\!\bigg(\frac{1}{1+\frac{k_2^2}{k_1^2}}\left(\!\alpha (c_{0,1}\tfrac{k_2}{k_1}-c_{0,2})\!+\!\mathtt{c}(|\bfk|^2)(c_{0,1}\!+\!c_{0,2}\tfrac{k_2}{k_1})\!\right)\!\!\!\left(c_{0,1}\tfrac{k_2}{k_1}-c_{0,2}\right)\!-\!\alpha c_{0,2}(c_{0,2}\!-\!c_{0,1}\cot\alpha)\!\!\bigg)}_{\displaystyle = O(|(k_1,\tfrac{k_2}{k_1})|(1+\tfrac{k_2^2}{k_1^2})^{-1})}\!\FF[\eta_1\rho_1]\bigg\|_0^2,\\
&\Big\|\dfrac{1}{D^2}(\bfc_0\cdot\bfL)\bfD\cdot(\eta_1\rho_1\bfw)
-\alpha (-c_{0,2}+c_{0,1}\cot\alpha)\eta_1\rho_1w_1\Big\|_0^2\\
& =
\bigg\|\underbrace{\bigg(\frac{1}{1+\frac{k_2^2}{k_1^2}}\left(\alpha (c_{0,1}\tfrac{k_2}{k_1}-c_{0,2})+\mathtt{c}(|\bfk|^2)(c_{0,1}+c_{0,2}\tfrac{k_2}{k_1})\right)
-\alpha(-c_{0,2}+c_{0,1}\cot\alpha)\bigg)}_{\displaystyle =O(|(k_1,\tfrac{k_2}{k_1})|(1+\tfrac{k_2^2}{k_1^2})^{-1})}\FF[\eta_1\rho_1]w_1\bigg\|_0^2
\\
& \qquad\qquad \mbox{}+ \bigg\|\underbrace{\bigg(\frac{1}{1+\frac{k_2^2}{k_1^2}}\left(\alpha (c_{0,1}\tfrac{k_2}{k_1}-c_{0,2})+\mathtt{c}(|\bfk|^2)(c_{0,1}+c_{0,2}\tfrac{k_2}{k_1})\right)\tfrac{k_2}{k_1}
\bigg)}_{\displaystyle =O(|(k_1,\tfrac{k_2}{k_1})|(1+\tfrac{k_2^2}{k_1^2})^{-1})}\FF[\eta_1\rho_1]w_2\bigg\|_0^2
\end{align*}
and
$$\bigg\|\frac{|k_1|}{1+\tfrac{k_2^2}{k_1^2}}\FF[\eta_1\rho_1]\bigg\|_0\leq \|k_1 \FF[\eta_1\rho_1]\|_0 = \|(\eta_1\rho_1)_x\|_0
\leq \|k_1 \hat{\eta}_1\|_0\|\rho_1\|_\infty + \|\eta_1\|_\infty \|k_1 \hat{\rho}_1\|_0 \leq \varepsilon^2 \nn \eta_1 \nn \nn \rho_1 \nn.\eqno{\Box}$$\pagebreak

\begin{corollary} \label{lot in m}
The estimate
$$b(\eta_1,\rho_1)=d_\alpha\eta_1\rho_1+B(\bfD)\eta_1\rho_1+O(\varepsilon^2\nn\eta_1\nn\nn\rho_1\nn),$$
where
$$|B(\bfk)| \lesssim |\tfrac{k_2}{k_1}|(1+\tfrac{k_2^2}{k_1^2})^{-1}$$
and $d_\alpha=\alpha \cosec \alpha + \tfrac{1}{2}\alpha \cot \alpha$,
holds for all $\eta_1$, $\rho_1 \in \XX_1$.
\end{corollary}
{\bf Proof.} This result is obtained by estimating each of the terms in the formula \eqref{Formula for J2} for $\JJ_2$ using Lemma \ref{Multiplier to const}.\qed

The remaining terms in the reduced equation are treated in the next lemma, which follows directly from Lemmata† \ref{Estimate for m},
\ref{Estimate for Jc}, Theorem \ref{Estimate for eta2} and Corollary \ref{lot in m}.

\begin{lemma}
The estimates
$$\varepsilon^2b(\eta_1,\eta_1) = \udl{O}(\varepsilon^2 \nn \eta_1 \nn^2), \qquad b(\eta_1,\eta_2(\eta_1))=\udl{O}(\varepsilon^2\nn \eta_1 \nn^3), \qquad
b(\eta_2(\eta_1),\eta_2(\eta_1))=\udl{O}(\varepsilon^2\nn \eta_1 \nn^4).$$
and
$$\JJ_{\ge 3}(\eta_1+\eta_2(\eta_1)) = \udl{O}(\varepsilon^2 \nn \eta_1 \nn^3)$$
hold for all $\eta_1 \in X_1$.
Here
the symbol $\udl{O}(\varepsilon^\gamma \nn \eta_1\nn^r)$ (with $\gamma \geq 0$, $r \geq 1$)
denotes a smooth function\linebreak
 $\RR^\varepsilon\colon X_1 \rightarrow H^1({\mathbb R}^2)$ which satisfies the estimates 
$$
\|\RR^\varepsilon(\eta_1)\|_1 \lesssim \varepsilon^\gamma  \nn \eta_1\nn^r, \quad
\|\mathrm{d}\RR^\varepsilon[\eta_1]\|_{L(\XX_1,H^1({\mathbb R}^2))}\lesssim \varepsilon^\gamma  \nn \eta_1 \nn^{r-1}$$
for each $\eta_1 \in X_1$.
\end{lemma}

Altogether we conclude that \eqref{Red eqn 1} can be written as
$$
g(\bfD)\eta_1 +2\varepsilon^2\frac{1}{D^2}(\bfc_0\cdot\bfD)(\bfc_0\cdot\bfL)\eta_1
-\varepsilon^4\frac{1}{D^2}(\bfc_0\cdot\bfD)(\bfc_0\cdot\bfL)\eta_1
+ \chi(\bfD)\Big(d_\alpha\eta_1^2+B(\bfD)\eta_1^2
+\udl{O}(\varepsilon^2\nn\eta_1\nn^2)\!\Big)=0,
$$
and applying Lemma \ref{Multiplier to const}(iv) one can further simplify it to
$$
g(\bfD)\eta_1 +2\varepsilon^2\eta_1
+ \chi(\bfD)\Big(d_\alpha\eta_1^2+B(\bfD)\eta_1^2
+\udl{O}(\varepsilon^3\nn\eta_1\nn)+\udl{O}(\varepsilon^2\nn\eta_1\nn^2)\!\Big)=0.
$$
The reduction is completed by introducing the KP scaling
$$\eta_1(x,y)=\varepsilon^2 \zeta (\varepsilon x, \varepsilon^2 y),$$
noting that $I\colon \eta_1 \rightarrow \zeta$ is an isomorphism $\XX_1 \rightarrow Y_1^\varepsilon$
and $\chi(\bfD)L^2({\mathbb R}^2) \rightarrow Y_0^\varepsilon$ and choosing $R>1$ large enough so that
$\zeta_k^\star \in B_R(\bfzero)$ (and $\varepsilon>0$ small enough so that 
$B_R(\bfzero) \subseteq Y_1^\varepsilon$  is contained in $I[X_1]$). Here we have replaced $(\chi(\bfD)H^1({\mathbb R}^2),\|\cdot\|_1)$
and $(\chi_\varepsilon(\bfD)L^2({\mathbb R}^2),\|\cdot\|_0)$,
where
$$\chi_\varepsilon(\bfk) = \chi(\varepsilon k_1, \varepsilon^2 k_2)$$
is the indicator function of the set
$$
S_\varepsilon=\left\{ \bfk \in {\mathbb R}^2 \colon |k_1| \leq \frac{\delta}{\varepsilon}, \Big|\frac{k_2}{k_1}\Big|\leq \frac{\delta}{\varepsilon}\right\},
$$
by the identical spaces
$(\chi(\bfD)L^2({\mathbb R}^2),\|\cdot\|_0)$ and $(Y_0^\varepsilon,\|\cdot\|_{Y_0})$ 
in order to work exclusively with the scales\linebreak
$\{Y_s,\|\cdot\|_{Y^s}\}_{s\geq 0}$ and $\{Y_s^\varepsilon,\|\cdot\|_{Y_s}\}_{s \geq 0}$ of function spaces.
We find that $\zeta \in B_R(\bfzero) \subseteq Y_1^\varepsilon$ satisfies the equation
\begin{equation}
\varepsilon^{-2}g_\varepsilon(\bfD)\zeta
+2\zeta
+ d_\alpha\chi_\varepsilon(\bfD)\zeta^2+\chi_\varepsilon(\bfD)B_\varepsilon(\bfD)\zeta^2
+\underline{\OO}_0^\varepsilon(\varepsilon^{\frac{1}{2}}\|\zeta\|_{Y_1})=0,
\label{Final reduced eqn}
\end{equation}
which now holds in $Y_0^\varepsilon$, where
\[g_\varepsilon(\bfk) = g(\varepsilon k_1, \varepsilon^2 k_2), \qquad B_\varepsilon(k_1,k_2) = B(\varepsilon k_1, \varepsilon^2 k_2)\]
and the symbol $\underline{\OO}^\varepsilon_n(\varepsilon^s\| \zeta \|_{Y_1}^r)$
denotes a smooth function
$\RR\colon B_R(\bfzero) \subseteq Y_1^\varepsilon \rightarrow Y_n^\varepsilon$
which satisfies the estimates
$$
\|\RR(\zeta)\|_{Y_n} \lesssim \varepsilon^s\| \zeta \|_{Y_1}^r, \quad
\|\mathrm{d}\RR[\zeta]\|_{L(Y_1,Y_n)}\lesssim \varepsilon^s\| u \|_{Y_1}^{r-1}
$$
for each $\zeta \in B_R(\bfzero)\subseteq Y_1^\varepsilon$  (with $r \geq 1$, $s$, $n \geq 0$).
Note that $\nn\eta_1\nn^2=\varepsilon \|\zeta\|_{Y_1}^2$ and that the change of variables
from $(x,y)$ to $(\varepsilon x, \varepsilon^2y)$ introduces a further factor of $\varepsilon^{\frac{3}{2}}$ in the
remainder term. The invariance of the reduced equation under $\eta_1(x,y)\mapsto\eta_1(-x,-y)$
is inherited by \eqref{Final reduced eqn}, which is invariant under the reflection $\zeta(x,y)\mapsto\zeta(-x,-y)$.

\section{Solution of the reduced equation} \label{sec:existence}

In this section we find solitary-wave solutions of the reduced equation \eqref{Final reduced eqn},
noting that in the formal limit $\varepsilon \rightarrow 0$ it reduces to the stationary KP-I equation
$$- (\beta-\beta_0) \zeta_{xx} +2 \zeta +\sec^2\!\tfrac{1}{2}\alpha \frac{D_2^2}{D_1^2}\zeta
+ d_\alpha \zeta^2=0,$$
which has explicit solitary-wave solutions $\zeta_k^\star$.
For this purpose we use a  perturbation argument, rewriting \eqref{Final reduced eqn}
as a fixed-point equation and applying the following version of the implicit-function theorem.

\begin{theorem} \label{IFT}
Let $\WW$ be a Banach space, $W_0$ and $\Lambda_0$ be open neighbourhoods of respectively $w^\star$ in $\WW$ and the origin in ${\mathbb R}$, and $\HH\colon  W_0 \times \Lambda_0 \rightarrow \WW$ be a function which is differentiable with respect to $w \in W_0$ for each $\lambda \in \Lambda_0$.
Furthermore, suppose that 
$\HH(w^\star,0)=0$, $\mathrm{d}_1\HH[w^\star,0]\colon \WW \rightarrow \WW$ is an isomorphism,
\[
\lim_{w \rightarrow w^\star}\|\mathrm{d}_1\HH[w, 0]-\mathrm{d}_1\HH[w^\star,0]\|_{L(\WW,\WW)}=0
\]
and
\[\lim_{\lambda \rightarrow 0} \|\HH(w,\lambda)-\HH(w,0)\|_{\WW}=0, \quad \lim_{\lambda \rightarrow 0} \
\|\mathrm{d}_1\HH[w,\lambda]-\mathrm{d}_1\HH[w,0]\|_{L(\WW,\WW)}=0\]
uniformly over $w \in W_0$.

There exist open neighbourhoods $W \subseteq W_0$ of $w^\star$ in $\WW$ and $\Lambda \subseteq \Lambda_0$ of the origin in ${\mathbb R}$, and a uniquely determined mapping
$h\colon \Lambda \rightarrow W$ with the properties that
\begin{list}{(\roman{count})}{\usecounter{count}}
\item
$h$ is continuous at the origin with $h(0)=w^\star$,
\item
$\HH(h(\lambda),\lambda)=0$ for all $\lambda \in \Lambda$,
\item
$w=h(\lambda)$ whenever $(w,\lambda) \in  W \times \Lambda$ satisfies $\HH(w,\lambda)=0$.
\end{list}
\end{theorem}

\begin{theorem} \label{Final strong existence thm}
For each sufficiently small value of $\varepsilon>0$ equation \eqref{Final reduced eqn} has a solution
$\zeta_k^\varepsilon$ in $Y_{1+\theta}^\varepsilon$ with\linebreak
$\zeta(x,y)=\zeta(-x,-y)$ for all $(x,y) \in {\mathbb R}^2$ and
$\|\zeta_k^\varepsilon-\zeta_k^\star\|_{Y_{1+\theta}} \rightarrow 0$ as $\varepsilon \rightarrow 0$.
\end{theorem}

The first step in the proof of Theorem \ref{Final strong existence thm} is to write \eqref{Final reduced eqn} as the fixed-point equation
\begin{equation}
\zeta+\varepsilon^2\left(g_\varepsilon(\bfD)+2\varepsilon^2\right)^{-1}
\Big(d_\alpha\chi_\varepsilon(\bfD)\zeta^2+\chi_\varepsilon(\bfD)B_\varepsilon(\bfD)\zeta^2
+\udl{\OO}_0^\varepsilon(\varepsilon^{\frac{1}{2}}\|\zeta\|_{Y_1})\!\Big)=0 \label{Final reduced eqn FP}
\end{equation}
and use the following result to `replace' the
nonlocal operator with the KP operator
$$L_\alpha=2-(\beta-\beta_0)\partial_x^2+\sec^2\!\tfrac{1}{2}\alpha\ \frac{D_2^2}{D_1^2}.$$

\begin{proposition} \label{Replace g with L}
Suppose that $\theta \in [0,1]$. The inequality
$$\left|\frac{\varepsilon^2}{2\varepsilon^2 + \tilde{g}(\varepsilon k_1,\varepsilon \frac{k_2}{k_1})}- \frac{1}{2+(\beta-\beta_0)k_1^2+\sec^2\!\tfrac{1}{2}\alpha\ \frac{k_2^2}{k_1^2}}\right| \lesssim \frac{\varepsilon^{1-\theta}}{(1+|(k_1,\frac{k_2}{k_1})|^2)^{\frac{1}{2}(1+\theta)}}$$
holds uniformly over $|k_1|, \frac{|k_2|}{|k_1|} < \delta/\varepsilon$.
\end{proposition}
{\bf Proof.} Clearly
\begin{align*}
\left|\frac{\varepsilon^2}{2\varepsilon^2 + \tilde{g}(\varepsilon k_1,\varepsilon \frac{k_2}{k_1})}- \frac{1}{2+(\beta-\beta_0)k_1^2+\sec^2\!\tfrac{1}{2}\alpha\ \frac{k_2^2}{k_1^2}}\right| \\
&\hspace{-2.25in}=\frac{\big|\tilde{g}(\varepsilon k_1,\varepsilon \frac{k_2}{k_1})-(\beta-\beta_0)\varepsilon^2k_1^2-\sec^2\!\tfrac{1}{2}\alpha\ \varepsilon^2\frac{k_2^2}{k_1^2}\big|}{\left(2\varepsilon^2 + \tilde{g}(\varepsilon k_1,\varepsilon \frac{k_2}{k_1})\right)\left(2+(\beta-\beta_0)k_1^2+\sec^2\!\tfrac{1}{2}\alpha\ \frac{k_2^2}{k_1^2}\right)};
\end{align*}
furthermore
$$\left|\tilde{g}\left(s_1,\frac{s_2}{s_1}\right)-(\beta-\beta_0)s_1^2 - \sec^2\!\tfrac{1}{2}\alpha\ \frac{s_2^2}{s_1^2}\right| \lesssim \left|\left(s_1,\frac{s_2}{s_1}\right)\right|^3,$$
and
$$\tilde{g}\left(s_1,\frac{s_2}{s_1}\right) \gtrsim \left|\left(s_1,\frac{s_2}{s_1}\right)\right|^2$$
for $|s_1|$, $\frac{|s_2|}{|s_1|} \leq \delta$ and sufficiently small $\delta$ (see Remark \ref{Near zero}).

It follows that
\begin{align*}
\left|\frac{\varepsilon^2}{2\varepsilon^2 + \tilde{g}(\varepsilon k_1,\varepsilon \frac{k_2}{k_1})}- \frac{1}{2+(\beta-\beta_0)k_1^2+\sec^2\!\tfrac{1}{2}\alpha\ \frac{k_2^2}{k_1^2}}\right|
& \lesssim
\frac{\varepsilon |(k_1,\frac{k_2}{k_1})|^3}{(1+|(k_1,\frac{k_2}{k_1})|^2)^2} \\
& \lesssim \frac{\varepsilon}{(1+|(k_1,\frac{k_2}{k_1})|^2)^{\frac12}}
\end{align*}
uniformly over $|k_1|, \frac{|k_2|}{|k_1|} < \delta/\varepsilon$, and the stated result follows from this inequality and
the observation that\linebreak
$\varepsilon \lesssim \delta(1+t^2)^{-\frac{1}{2}}$ when $|t|< \delta/\varepsilon$.\qed

\begin{lemma}
Suppose that $\theta \in [0,1]$. The estimates
\begin{align*}
\varepsilon^2\left(g_\varepsilon(\bfD)+2\varepsilon^2\right)^{-1}\udl{\OO}_0^\varepsilon(\varepsilon^{\frac{1}{2}}\|\zeta\|_{Y_1}) &=
\udl{\OO}_{1+\theta}^\varepsilon(\varepsilon^{\frac{1}{2}}\|\zeta\|_{Y_{1+\theta}}), \\
\varepsilon^2\left(g_\varepsilon(\bfD)+2\varepsilon^2\right)^{-1}B_\varepsilon(\bfD)\zeta^2 &= \udl{\OO}_{1+\theta}^\varepsilon(\varepsilon^{\frac{1}{2}-\frac{\theta}{2}}\|\zeta\|_{1+\theta}^2)
\end{align*}
and
$$\left(\varepsilon^2\left(g_\varepsilon(\bfD)+2\varepsilon^2\right)^{-1}-L_\alpha^{-1}\right)\chi_\varepsilon(\bfD)\zeta^2=\udl{\OO}_{1+\theta}^\varepsilon(\varepsilon^{1-\theta}\|\zeta\|_{1+\theta}^2)
$$
hold for all $\zeta \in Y_{1+\theta}^\varepsilon$.
\end{lemma}
{\bf Proof.} It follows from Proposition \ref{Replace g with L} (with $\theta=1$) that
$$
\frac{\varepsilon^2}{2\varepsilon^2 + \tilde{g}(\varepsilon k_1,\varepsilon \frac{k_2}{k_1})} \lesssim \left(1+k_1^2+\sdfrac{k_2^2}{k_1^2}\right)^{-1},
$$
from which the first estimate is an immediate consequence (note that $\|\zeta\|_1 \leq \|\zeta\|_{1+\theta}$). Furthermore%
\begin{align*}
\bigg(1+k_1^2+&\sdfrac{k_2^2}{k_1^2}\bigg)^{\! \frac{1}{2}+\frac{\theta}{2}}\frac{\varepsilon^2}{2\varepsilon^2 + \tilde{g}(\varepsilon k_1,\varepsilon \frac{k_2}{k_1})}B(\varepsilon k_1, \varepsilon^2 k_2)\\
& \lesssim \bigg(1+k_1^2+\sdfrac{k_2^2}{k_1^2}\bigg)^{\! -\frac{1}{2}+\frac{\theta}{2}}\frac{\varepsilon \big|\frac{k_2}{k_1}\big|}{1+\varepsilon^2 \frac{k_2^2}{k_1^2}}\\
& = \varepsilon^{\frac{1}{2}-\frac{\theta}{2}} \left(\frac{\big|\frac{k_2}{k_1}\big|}{1+k_1^2+\frac{k_2^2}{k_1^2}}\frac{1}{1+\varepsilon^2 \frac{k_2^2}{k_1^2}}\right)^{\!\!\frac{1}{2}-
\frac{\theta}{2}}
\left(\frac{\varepsilon \big|\frac{k_2}{k_1}\big|}{1+\varepsilon^2 \frac{k_2^2}{k_1^2}}\right)^{\!\!\frac{1}{2}+\frac{\theta}{2}} \\
& \lesssim \varepsilon^{\frac{1}{2}-\frac{\theta}{2}},
\end{align*}
such that
$$\left\|\varepsilon^2\left(g_\varepsilon(\bfD)+2\varepsilon^2\right)^{-1}B_\varepsilon(\bfD)\zeta\xi\right\|_{Y_{1+\theta}}
\lesssim \varepsilon^{\frac{1}{2}-\frac{\theta}{2}} \|\zeta\xi\|_0 \leq \varepsilon^{\frac{1}{2}-\frac{\theta}{2}}
\|\zeta\|_{L^4({\mathbb R}^2)} \|\xi\|_{L^4({\mathbb R}^2)} \lesssim \varepsilon^{\frac{1}{2}-\frac{\theta}{2}}\|\zeta\|_{Y_{1+\theta}}\|\xi\|_{Y_{1+\theta}}.$$
for all $\zeta$, $\xi \in Y_{1+\theta}^\varepsilon$ (see Proposition \ref{Properties of Y}(i)).

The final estimate follows from the observation that
$$\left\|\left(\varepsilon^2\left(g_\varepsilon(\bfD)+2\varepsilon^2\right)^{-1}-L_\alpha^{-1}\right)\chi_\varepsilon(\bfD)\zeta\xi\right\|_{Y_{1+\theta}}
\lesssim \varepsilon^{1-\theta}\|\zeta\xi\|_0 \lesssim \varepsilon^{1-\theta}\|\zeta\|_{Y_{1+\theta}}\|\xi\|_{Y_{1+\theta}}$$
for all $\zeta$, $\xi \in Y_{1+\theta}^\varepsilon$, in which the first inequality follows from Proposition \ref{Replace g with L}.\qed\pagebreak

Using the above lemma, one can write equation \eqref{Final reduced eqn FP} as
$$
\zeta+F_\varepsilon(\zeta)=0,
$$
in which
$$
F_\varepsilon(\zeta)=d_\alpha L_\alpha^{-1}\chi_\varepsilon(\bfD)\zeta^2
+\udl{\OO}_{1+\theta}^\varepsilon(\varepsilon^{\frac{1}{2}-\frac{\theta}{2}}\|\zeta\|_{1+\theta}).$$
It is convenient to replace this equation with
$$\zeta+\tilde{F}_\varepsilon(\zeta)=0,$$
where $\tilde{F}_\varepsilon(\zeta) = F_\varepsilon(\chi_\varepsilon(\bfD)\zeta)$ and study it in the fixed space $Y_{1+\theta}$
for $\theta\in(\frac{1}{2},1)$
(the solution sets of the two equations evidently coincide); we choose $\theta>\frac{1}{2}$ so that $Y_{1+\theta}$ is embedded
in $C_\mathrm{b}({\mathbb R}^2)$ and $\theta<1$ so that the remainder term in $\tilde{F}_\varepsilon(\zeta)$ vanishes at $\varepsilon=0$.

We establish Theorem \ref{Final strong existence thm} by applying Theorem \ref{IFT} with
$$\WW=Y_{1+\theta}^\mathrm{e}\coloneqq  \{\zeta \in Y_{1+\theta}\colon \mbox{$\zeta(x,y)=\zeta(-x,-y)$ for all $(x,y) \in {\mathbb R}^2$}\},$$
$W_0=B_R(\bfzero)\subseteq Y_{1+\theta}$, $\Lambda_0=(-\varepsilon_0,\varepsilon_0)$ for a sufficiently small value of $\varepsilon_0$,
and
$$\HH(\zeta,\varepsilon)\coloneqq \zeta+\tilde{F}_{|\varepsilon|}(\zeta)$$
(here $\varepsilon$ is replaced by $|\varepsilon|$ so that $\HH(\zeta,\varepsilon)$ is defined for $\varepsilon$ in
a full neighbourhood of the origin in ${\mathbb R}$). 

We begin by verifying that the functions \(\zeta_k^\star\) belong to $Y_{1+\theta}^\mathrm{e}$.

\begin{proposition}
Each lump solution \(\zeta_k^\star\) belongs to $Y_2$.
\end{proposition}
{\bf Proof.}
First note that $(\zeta_k^\star)^2$ belongs to $L^2({\mathbb R}^2)=Y_0$ because \(|\zeta_k^\star(x,y)| \lesssim (1+x^2+y^2)^{-1}\) for all $(x,y) \in {\mathbb R}^2$
(see Lemma \ref{LWY results}(i)). Since \(\zeta_k^\star\) satisfies
$$
 \zeta_k^\star+L_\alpha^{-1}(\zeta_k^\star)^2=0
$$
and $L_\alpha^{-1}$ is a regularising operator of order \(2\) for the scale $\{Y_r,\|\cdot\|_{Y_r}\}_{r \geq 0}$, one finds that \(\zeta_k^\star \in Y_2\). 
\qed

Observe that \(\HH(\cdot,\varepsilon)\) is a continuously differentiable function \(B_R(\bfzero) \subseteq Y_\mathrm{e}^{1+\theta} \to Y_\mathrm{e}^{1+\theta}\)
for each fixed $\varepsilon \geq 0$, so that
\[
\lim_{\zeta \rightarrow \zeta_k^\star}\|\mathrm{d}_1\HH[\zeta, 0]-\mathrm{d}_1\HH[\zeta_k^\star,0]\|_{L(Y_{1+\theta},Y_{1+\theta})}=0.
\]
The facts that
$$
 \lim_{\varepsilon \rightarrow 0} \|\HH(\zeta,\varepsilon)-\HH(\zeta,0)\|_{Y_{1+\theta}}=0, \quad \lim_{\varepsilon \rightarrow 0} \
\|\mathrm{d}_1\HH[\zeta,\varepsilon]-\mathrm{d}_1\HH[\zeta,0]\|_{L(Y_{1+\theta},Y_{1+\theta})}=0
$$
uniformly over $\zeta \in B_R(\bfzero)\subseteq Y_\mathrm{e}^{1+\theta}$ are obtained from the equation
\[
\HH(\zeta,\varepsilon)-\HH(\zeta,0)
= L_\alpha^{-1}
\left(\chi_\varepsilon(\mathbf{D})
\left(\chi_\varepsilon(\mathbf{D})\zeta\right)^2-\zeta^2\right)
 + \udl{\OO}_{1+\theta}^\varepsilon(\varepsilon^{\frac{1}{2}-\frac{\theta}{2}}\|\zeta\|_{1+\theta})
 \]
using Corollary~\ref{Difference of G 3} below, which is a consequence of the next two lemmas.

\begin{lemma} \label{Difference of G 1}
Fix $\theta > \frac{1}{2}$. The estimate
\[
\| L_\alpha^{-1} \chi_\varepsilon(\mathbf{D}) \big( (( \chi_\varepsilon(\mathbf{D}) + I) \zeta) ((\chi_\varepsilon(\mathbf{D}) - I) \xi)\big) \|_{Y_{1+\theta}} \lesssim \varepsilon \| \zeta \|_{Y_{1+\theta}} \| \xi \|_{Y_{1+\theta}}
\]
holds for all \(\zeta, \xi \in Y_{1+\theta}\).
\end{lemma}
{\bf Proof.}
Recall that \(L_\alpha^{-1}\) is a regularising operator of order 2 for the scale $\{Y_r, \|\cdot\|_{Y_r}\}_{r \geq 0}$ and that \(\chi_\varepsilon(\mathbf{D})\) is a bounded projection on all subspaces of $L^2({\mathbb R}^2)$. It follows that
\begin{align*}
\| L_\alpha^{-1} & \chi_\varepsilon(\mathbf{D}) \big( (( \chi_\varepsilon(\mathbf{D}) + I) \zeta) ((\chi_\varepsilon(\mathbf{D}) - I) \xi)\big) \|_{Y_{1+\theta}}\\
&\quad\leq \| \chi_\varepsilon(\mathbf{D}) \big( (( \chi_\varepsilon(\mathbf{D}) + I) \zeta) ((\chi_\varepsilon(\mathbf{D}) - I) \xi)\big)\|_0\\
&\quad \leq  \|(( \chi_\varepsilon(\mathbf{D}) + I) \zeta) ((\chi_\varepsilon(\mathbf{D}) - I) \xi) \|_0\\
&\quad \leq  \| (\chi_\varepsilon(\mathbf{D}) + I) \zeta\|_\infty \| (\chi_\varepsilon(\mathbf{D}) - I) \xi \|_0 \\
&\quad \lesssim \| (\chi_\varepsilon(\mathbf{D}) + I)\zeta\|_{Y_{1+\theta}} \| (\chi_\varepsilon(\mathbf{D}) - I) \xi \|_0 \\
&\quad \leq 2\| \zeta\|_{Y_{1+\theta}} \| (\chi_\varepsilon(\mathbf{D}) - I) \xi \|_0,
\end{align*}\pagebreak
where we have used the embedding $Y_{1+\theta}  \hookrightarrow C_\mathrm{b}({\mathbb R}^2)$.
To estimate \(\| (\chi_\varepsilon(\mathbf{D}) - I) \zeta\|_0\), note that
\[
\mathbb{R}^2 \setminus S_\varepsilon \subset \underbrace{\left\{ (k_1,k_2) \colon |k_1| > \frac{\delta}{\varepsilon}\right\}}_{\displaystyle= S_\varepsilon^1} \cup  \underbrace{\left\{ (k_1,k_2) \colon \left|\frac{k_2}{k_1}\right| > \frac{\delta}{\varepsilon}\right\}}_{\displaystyle=S_\varepsilon^2},
\]
so that
\begin{align*}
\hspace{1.5in}\| (\chi_\varepsilon(\mathbf{D}) - I) \zeta \|_0^2 &= \int_{\mathbb{R}^2 \setminus S_\varepsilon} |\hat\zeta|^2 \dbfk \\
& \leq \int_{S_\varepsilon^1} |\hat\zeta|^2 \dbfk + \int_{S_\varepsilon^2} |\hat\zeta|^2 \dbfk\\
&\leq \frac{\varepsilon^2}{\delta^2} \int_{S_\varepsilon^1} k_1^2 |\hat\zeta|^2 \dbfk +  \frac{\varepsilon^2}{\delta^2} \int_{S_\varepsilon^2} \frac{k_2^2}{k_1^2} |\hat\zeta|^2 \dbfk \\
& \leq \frac{2\varepsilon^2}{\delta^2} \| \zeta \|_{Y_1}^2.\hspace{3in}\Box
\end{align*}

\begin{lemma} \label{Difference of G 2}
Fix $\theta \in (0,1)$. The estimate
\[
\| L_\alpha^{-1} ( \chi_\varepsilon(\mathbf{D}) - I) (\zeta \xi) \|_{Y_{1 +\theta}} \lesssim \varepsilon^{1-\theta} \|\zeta\|_{Y_1} \|\xi\|_{Y_1} \leq \varepsilon^{\frac{1}{2}-\frac{\theta}{2}} \|\zeta\|_{Y_{1+\theta}} \|\xi\|_{Y_{1+\theta}},
\]
holds for all \(\zeta, \xi \in Y_{1+\theta}\).
\end{lemma}
{\bf Proof.}
For \(\nu \in \{ k_1,  \frac{k_2}{k_1}\}\) we find that
\begin{align*}
&\bigg( 1 + k_1^2 + \frac{k_2^2}{k_1^2} \bigg)^{1 + \theta} \bigg( 1 + k_1^2 + \frac{k_2^2}{k_1^2} \bigg)^{-2} |\nu|^{2-2\theta}
=\Bigg( \frac{\nu^2}{1 + k_1^2 +\frac{k_2^2}{k_1^2}} \Bigg)^{1-\theta} \leq 1,
\end{align*}
so that%
\begin{align*}
\| L_\alpha^{-1} &( \chi_\varepsilon(\mathbf{D}) - I) \zeta \xi\|_{Y_{1+\theta}}^2\\
&\lesssim \int_{S_\varepsilon^1 \cup S_\varepsilon^2} \left( 1 + k_1^2 + \frac{k_2^2}{k_1^2}\right)^{1+\theta} \left( 1 + k_1^2 + \frac{k_2^2}{k_1^2}\right)^{-2} |\FF[\zeta \xi]|^2 \dbfk\\
&\lesssim \left( \frac{\varepsilon}{\delta} \right)^{2-2\theta} \int_{S_\varepsilon^1} \left( 1 + k_1^2 + \frac{k_2^2}{k_1^2}\right)^{1+\theta} \left( 1 + k_1^2 + \frac{k_2^2}{k_1^2}\right)^{-2} |k_1|^{2-2\theta} |\FF[\zeta \xi]|^2 \dbfk\\
&\qquad \mbox{}+ \left( \frac{\varepsilon}{\delta} \right)^{2-2\theta} \int_{S_\varepsilon^2} \left( 1 + k_1^2 + \frac{k_2^2}{k_1^2}\right)^{1+\theta} \left( 1 + k_1^2 + \frac{k_2^2}{k_1^2}\right)^{-2} \left| \frac{k_2}{k_1} \right|^{2-2\theta} |\FF[\zeta \xi]|^2 \dbfk\\
& \leq \left( \frac{\varepsilon}{\delta} \right)^{2-2\theta} \|\zeta \xi\|_0^2  \\
& \lesssim \left( \frac{\varepsilon}{\delta} \right)^{2-2\theta} \|\zeta \|_{L^4({\mathbb R}^2)}^2 \|\xi \|_{L^4({\mathbb R}^2)}^2 \\
&  \lesssim \left( \frac{\varepsilon}{\delta} \right)^{2-2\theta} \|\zeta \|_{Y_1}^2 \|\xi \|_{Y_1}^2,
\end{align*}
where we have used Parseval's theorem, the Cauchy--Schwarz inequality and the embedding\linebreak
\(Y_1 \hookrightarrow L^4({\mathbb R}^2)\).
\qed

\begin{corollary} \label{Difference of G 3}
Fix $\theta \in (\frac{1}{2},1)$. The estimate
\[
\left\| L_\alpha^{-1} \Big( \chi_\varepsilon(\mathbf{D}) \big( ( \chi_\varepsilon(\mathbf{D}) \zeta) (\chi_\varepsilon(\mathbf{D}) \xi) \big) - \zeta \xi\Big)\right \|_{Y_{1+\theta}} \lesssim \varepsilon^{1-\theta} \|\zeta\|_{Y_{1+\theta}} \|\xi\|_{Y_{1+\theta}}
\]
holds for all $\zeta$, $\xi \in Y_{1+\theta}$.
\end{corollary}
{\bf Proof.}
This result is obtained by writing
\begin{align*}
L_\alpha^{-1} & \Big( \chi_\varepsilon(\mathbf{D}) \big( ( \chi_\varepsilon(\mathbf{D}) \zeta) (\chi_\varepsilon(\mathbf{D}) \xi) \big) - \zeta \xi\Big)\\
&=  \tfrac{1}{2} L_\alpha^{-1} \chi_\varepsilon(\mathbf{D}) \big( ((\chi_\varepsilon(\mathbf{D}) +1) \zeta) ( (\chi_\varepsilon(\mathbf{D}) -1) \xi) \big)\\
&\qquad \mbox{}+ \tfrac{1}{2} L_\alpha^{-1} \chi_\varepsilon(\mathbf{D}) \big( ((\chi_\varepsilon(\mathbf{D}) +1) \xi) ( (\chi_\varepsilon(\mathbf{D}) -1) \zeta) \big) + L_\alpha^{-1} (\chi_\varepsilon(\mathbf{D}) -1) (\zeta \xi),
\end{align*}
and applying Lemma \ref{Difference of G 1} to the first two terms on the right-hand side and Lemma \ref{Difference of G 2} to the third.
\qed

It thus remains to show that
$$\mathrm{d}_1\HH[\zeta_k^\star,0]=I+2d_\alpha L_\alpha^{-1}(\zeta_k^\star\cdot)$$
is an isomorphism; this fact follows from the following result.

\begin{lemma}
The operator $L_\alpha^{-1}(\zeta_k^\star\cdot)\colon Y_{1+\theta} \rightarrow Y_{1+\theta}$ is compact.
\end{lemma}
 {\bf Proof.} 
 Let $\{\zeta_j\}$ be a sequence which is bounded in $Y_1$.
We can find a subsequence of $\{\zeta_j\}$ (still denoted by $\{\zeta_j\}$)
which converges weakly in $L^2({\mathbb R}^2)$ (because $\{\zeta_j\}$ is bounded in $L^2({\mathbb R}^2)$)
and strongly in $L^2(|\mathbf{x}|<n)$ for each $n \in {\mathbb N}$ (by Proposition \ref{Properties of Y}(ii) and a `diagonal'
argument). Denote the limit by $\zeta_\infty$. Since
$$\|\zeta_k^\star\zeta_j-\zeta_k^\star\zeta_\infty\|_{L^2(|\mathbf{x}|<n)}\leq \|\zeta_k^\star\|_\infty \|\zeta_j-\zeta_\infty\|_{L^2(|\mathbf{x}|<n)} \rightarrow 0$$
as $j \rightarrow \infty$ for each $n \in {\mathbb N}$ and
$$\sup_j\|\zeta_k^\star\zeta_j\|_{L^2(|\mathbf{x}|>n)}\leq \sup_{|\mathbf{x}|>n}|\zeta_k^\star({\mathbf x})|\sup_j\|\zeta_j\|_0 \rightarrow 0$$
as $n \rightarrow \infty$ we conclude that $\{\zeta_k^\star\zeta_j\}$ converges to $\zeta_k^\star\zeta_\infty$ in $L^2({\mathbb R}^2)$
as $j \rightarrow \infty$. It follows that
$\zeta \mapsto \zeta_k^\star \zeta$ is compact\linebreak $Y_1 \rightarrow L^2({\mathbb R})$ and hence $Y_{1+\theta} \rightarrow L^2({\mathbb R})$;
the result follows from this fact and the observation that $L_\alpha^{-1}$ is continuous $L^2({\mathbb R}^2) \rightarrow Y_{1+\theta}$.\qed%

\begin{lemma}
The operator $I+2d_\alpha L_\alpha^{-1}(\zeta_k^\star\cdot)$ is an isomorphism $Y_{1+\theta} \rightarrow Y_{1+\theta}$.
\end{lemma}
{\bf Proof.} The previous result shows that $I+2d_\alpha L_\alpha^{-1}(\zeta_k^\star\cdot)\colon Y_{1+\theta} \rightarrow Y_{1+\theta}$ is Fredholm with index $0$; it therefore remains to show that it is injective.

Suppose that $\zeta \in Y_{1+\theta}$ satisfies
\begin{equation}
\zeta+2d_\alpha L_\alpha^{-1}(\zeta_k^\star\zeta)=0. \label{Integral form}
\end{equation}
It follows that
$$k_1\hat{\zeta}=\frac{-2d_\alpha k_1^3}{2k_1^2+(\beta-\beta_0)k_1^4+\sec^2\!\tfrac{1}{2}\alpha\ k_2^2}\FF[\zeta_k^\star\zeta],
\qquad
k_2\hat{\zeta}=\frac{-2d_\alpha k_1^2k_2}{2k_1^2+(\beta-\beta_0)k_1^4+\sec^2\!\tfrac{1}{2}\alpha\ k_2^2}\FF[\zeta_k^\star\zeta]$$
and hence $\zeta \in H^{n+1}({\mathbb R}^2)$ whenever $\zeta_k^\star\zeta \in H^n({\mathbb R}^2)$.
Since $\zeta \in L^2({\mathbb R}^2)$ and $\zeta \in H^b({\mathbb R}^2)$ implies $\zeta_k^\star\zeta \in H^b({\mathbb R}^2)$
we find by bootstrapping that $\zeta \in H^\infty({\mathbb R}^2)$.

Since $\zeta$ is smooth and satisfies \eqref{Integral form} it satisfies the linear equation
$$\big((\beta-\beta_0)\zeta_{xx}+2\zeta+2d_\alpha(\zeta_k^\star\zeta)\big)_{xx}-\sec^2\!\tfrac{1}{2}\alpha\ \zeta_{zz}=0,$$
and the only smooth solution to this equation with $\zeta(x,y)=\zeta(-x,-y)$ for all $(x,y) \in {\mathbb R}^2$ is the trivial solution 
(see Lemma \ref{LWY results}(iii)).\qed

To establish Theorem \ref{Main result} it remains to confirm that the formula
$$\eta = \eta_1 + \eta_2(\eta_1),
\quad \eta_1(x,y) = \varepsilon^2 \zeta(\varepsilon x,\varepsilon^2 y)$$
leads to the estimate
$$\eta(x,y)=\varepsilon^2 \zeta_k^\star(\varepsilon x,\varepsilon^2 y) + o(\varepsilon^2)$$
uniformly over $(x,y) \in {\mathbb R}^2$. This fact follows from the calculations
$$
\|\zeta_k^\varepsilon-\zeta_k^\star\|_\infty \lesssim \|\zeta_k^\varepsilon-\zeta_k^\star\|_{Y_{1+\theta}} = o(1),
$$
such that
\begin{align*}
\eta_1(x,y) & = \varepsilon^2 \zeta_k^\star(\varepsilon x,\varepsilon^2 y)
+\varepsilon^2 (\zeta_k^\varepsilon-\zeta_k^\star)(\varepsilon x,\varepsilon^2 y)  \\
& = \varepsilon^2 \zeta_k^\star(\varepsilon x,\varepsilon^2 y) + o(\varepsilon^2)
\end{align*}
uniformly in $(x,y)$,
and
$$\|\eta_2(\eta_1)\|_\infty \lesssim \|\eta_2(\eta_1)\|_3 \lesssim \varepsilon \nn \eta_1 \nn^2 \lesssim \varepsilon^3$$
by Theorem \ref{Estimate for eta2} and $\nn \eta_1 \nn = \varepsilon\|\zeta\|_{Y_1}$ with $\zeta \in B_R(\bfzero) \subseteq Y_{1+\theta}^\varepsilon$.

\begin{appendices}
\section{Dispersion relation} \label{disprel}
\end{appendices}

Recall the dispersion relation
\begin{equation}
g(\bfk)=0, \label{dispersion relation}
\end{equation}
where
$$
g(\bfk)=-\frac{1}{|\bfk|^2}(\alpha(\bfc_0\cdot\bfk^\perp)(\bfc_0\cdot\bfk)+\mathtt{c}(|\bfk|^2)(\bfc_0\cdot\bfk)^2)+1+\beta |\bfk|^2
$$
is an analytic function $\tilde{g}$ of $k_1$ and $\frac{k_2}{k_1}$ with $\tilde{g}(0,0)=0$ if
$$\bfc_0=\begin{pmatrix}c_0\cos\tfrac{1}{2}\alpha \\[1mm] -c_0\sin\tfrac{1}{2}\alpha\end{pmatrix}, \qquad c_0^2=\tfrac{2}{\alpha} \tan \tfrac{1}{2}\alpha.$$

Suppose that $\bfk \neq \bfzero$, so that \eqref{dispersion relation} is equivalent to
$$
\mathtt{c}(|\bfk|^2)=\kappa(|\bfk|^2,\theta),
$$
where
$$\kappa(\mu,\theta)=\frac{1+\beta\mu-\alpha c_0^2 \sin\theta\cos\theta}{c_0^2\cos^2\!\theta}$$
and $\theta$ is the angle between $\bfc_0$ and $\bfk$ (note that $g(\bfk)>0$ if $\cos\theta=0$).
The function $\kappa(\mu,\cdot)$ takes every value in $[\kappa_\mathrm{min}(\mu),\infty)$, where
$$\kappa_\mathrm{min}(\mu)=\frac{1+\beta\mu}{c_0^2}-\frac{c_0^2\alpha^2}{4(1+\beta\mu)},$$
and the minimum is attained at
$$\theta=\tan^{-1}\frac{\alpha c_0^2}{2(1+\beta \mu)}.$$
It follows that $g(\bfk)\neq 0$ for all $\bfk$ with given magnitude $|\bfk|$ if and only if $\mathtt{c}(|\bfk|^2)<\kappa_\mathrm{min}(|\bfk|^2)$.

The functions $\mathtt{c}$ and $\kappa_\mathrm{min}$ are both strictly increasing and concave on $[0,\infty)$ with\linebreak
$\mathtt{c}(0)=\kappa_\mathrm{min}(0)$.
Obviously $\mathtt{c}(\mu)<\kappa_\mathrm{min}(\mu)$ for $\mu \in (0,\infty)$ if $\mathtt{c}^\prime(\mu)<\kappa_\mathrm{min}^\prime(\mu)$
for $\mu \in [0,\infty)$, and since
$$\mathtt{c}^\prime(\mu) \leq \mathtt{c}^\prime(0)=\frac{1}{2}\left(-\frac{\cot\alpha}{\alpha}+\cosec^2\!\alpha\right)$$
(because $\mathtt{c}$ is concave) and
$$\kappa_\mathrm{min}^\prime(\mu) = \frac{1}{2}\alpha\beta\left(\frac{1}{(1+\beta\mu)^2}+\cot^2\!\tfrac{1}{2}\alpha\right)\tan\tfrac{1}{2}\alpha>\tfrac{1}{2}\alpha\beta\cot\tfrac{1}{2}\alpha$$
this condition is met if
$$\frac{1}{2}\left(-\frac{\cot\alpha}{\alpha}+\cosec^2\!\alpha\right)<\tfrac{1}{2}\alpha\beta\cot\tfrac{1}{2}\alpha,$$
that is if
$$\beta>\beta^\star\coloneqq \tfrac{1}{\alpha}\left(-\tfrac{1}{\alpha}\cot\alpha+\cosec^2\!\alpha\right)\tan\tfrac{1}{2}\alpha.$$

\begin{remark} \label{Near zero}
The calculation
$$\tilde{g}(k_1,\tfrac{k_2}{k_1})=\left(\beta+\frac{1}{2\alpha^2}(\cos\alpha-\alpha\cosec\alpha)\right)k_1^2
+\sec^2\!\tfrac{1}{2}\alpha\ \tfrac{k_2^2}{k_1^2}+\udl{O}(|(k_1,\tfrac{k_2}{k_1})|^3)$$
as $(k_1,\tfrac{k_2}{k_1})\rightarrow(0,0)$ shows that $(0,0)$ is a strict local minimum of $\tilde{g}$ if
$$\beta>\beta_0\coloneqq \frac{1}{2\alpha^2}(-\cos\alpha+\alpha\cosec\alpha).$$
Note that
$$\beta^\star-\beta_0=\frac{1}{\alpha^2}\cosec^3\alpha\sin^4\tfrac{1}{2}\alpha(2\alpha-\sin2\alpha)\geq0$$
with equality if and only if $\alpha=0$ (the common value is $\tfrac{1}{3}$).
\end{remark}

\noindent\\
{\bf Acknowledgements} The authors would like to thank Boris Buffoni and Evgeniy Lokharu for helpful discussions during the preparation of this article.\\

\noindent
{\bf Funding}
E. W. was supported by the European Research Council (ERC) under the European Union's Horizon 2020 research and innovation programme (grant agreement no. 678698) and the Swedish Research Council (grant no. 2020-00440).\\
\\
The project was partially funded by the Deutsche Forschungsgemeinschaft (DFG, German Research Foundation), project number 545145736.

\end{document}